\documentclass[11pt,a4paper,fleqn]{article}
\usepackage{a4wide,amsfonts,amsmath,latexsym,amssymb,euscript,graphicx,units,mathrsfs}

\usepackage{graphicx}
\usepackage{color}
\usepackage{amssymb}
\usepackage{amssymb}
\usepackage[T1]{fontenc}
\usepackage{latexsym}
\usepackage{xypic}
\usepackage{eufrak}
\usepackage{euscript}
\usepackage{amsfonts,amsmath}
\usepackage{verbatim}
\usepackage{fancyhdr}
\usepackage[english]{babel}
\usepackage{mathrsfs}
\usepackage{units}
\usepackage{url}

\newtheorem{prop}{Proposition}[section]
\newtheorem{cor}[prop]{Corollary}
\newtheorem{lemma}[prop]{Lemma}
\newtheorem{rem}[prop]{Remark}
\newtheorem{remark}[prop]{Remark}

\newtheorem{thm}[prop]{Theorem}
\newtheorem{defi}[prop]{Definition}

\newcommand{\poly}{Q}
\newcommand{\ff}{\varphi}

\renewcommand{\geq}{\geqslant}
\def\leq{\leqslant}

\newcommand{\N}{\mathbb{N}}
\newcommand{\Z}{\mathbb{Z}}
\newcommand{\R}{\mathbb{R}}

\def\e{\varepsilon}

\def\1{{\mathbf{1}}}

\def\1{{\mathbf{1}}}
\def\0.5{{\frac{1}{2}}}

\newcommand{\qed}{\nopagebreak\hspace*{\fill}
{\vrule width6pt height6ptdepth0pt}\par}

\title{Lectures on Gaussian approximations with Malliavin calculus}

\author{
  Ivan Nourdin\\
\small
  Universit\'e de Lorraine,
  Institut de Math\'ematiques \'Elie Cartan\\
\small
  B.P. 70239,
  54506 Vandoeuvre-l\`es-Nancy Cedex, France\\
\small
 {\tt inourdin@gmail.com}\\
  }

\begin{document}

\maketitle

{\bf Overview}.
In a seminal paper of 2005, Nualart and Peccati \cite{nunugio} discovered a surprising central limit theorem (called the ``Fourth Moment Theorem'' in the sequel) for sequences of multiple stochastic integrals of a fixed order: in this context, convergence in distribution to the standard normal law is  equivalent to convergence of just the fourth moment. Shortly afterwards, Peccati and Tudor \cite{PTu04} gave a multidimensional version of this characterization.

  Since the publication of these two beautiful papers, many improvements and developments on this theme have been considered. Among them is the work by Nualart and Ortiz-Latorre \cite{NOL}, giving a new proof only based on Malliavin calculus and the use of integration by parts on Wiener space. A second step is my joint paper \cite{NP-PTRF} (written in collaboration with Peccati) in which, by bringing together Stein's method with Malliavin calculus, we have been able (among other things) to associate quantitative bounds to the Fourth Moment Theorem.
  It turns out that Stein's method and Malliavin calculus fit together admirably well. Their interaction has led to some remarkable new results involving central and non-central limit theorems for functionals of infinite-dimensional Gaussian fields.

The current survey aims to introduce the main features of this recent theory.
It originates from a series of lectures I delivered\footnote{You may watch the videos of the lectures at
\url{http://www.sciencesmaths-paris.fr/index.php?page=175}.} at the Coll\`ege de France between January and March 2012, within the framework of the annual prize of the Fondation des Sciences Math\'ematiques de Paris.
It may be seen as a teaser for the book \cite{book-malliavinstein}, in which the interested reader will
find much more than in this short survey.

\bigskip

{\bf Acknowledgments}.
It is a pleasure to thank the Fondation des Sciences Math\'ematiques de Paris
for its generous support during the academic year 2011-12 and for giving me the opportunity to
speak about my recent research in the prestigious Coll\`ege de France.
I am grateful to all the participants of these lectures for their assiduity.
Also, I would like to warmly thank two anonymous referees for their very careful reading and for their valuable suggestions and remarks.
Finally, my last thank goes to Giovanni Peccati, not only for accepting to give a lecture (resulting to the material
developed in Section \ref{sec:gio}) but also (and especially!) for all the nice theorems we recently discovered together. I do hope
it will continue this way
as long as possible!

\eject

\tableofcontents

\bigskip

\section{Breuer-Major Theorem}\label{sectionbreuer}

The aim of this first section is to illustrate, through a guiding example, the power of the approach we will develop in this survey.

\medskip

Let $\{X_k\}_{k\geq 1}$ be a centered stationary Gaussian family. In this context, stationary just means that there
exists $\rho:\Z\to\R$ such that $E[X_kX_l]=\rho(k-l)$, $k,l\geq 1$. Assume further that $\rho(0)=1$, that is,
each $X_k$ is $\mathcal{N}(0,1)$ distributed.

Let $\varphi:\R\to\R$ be a measurable function satisfying
\begin{equation}\label{integr}
E[\varphi^2(X_1)]=\frac{1}{\sqrt{2\pi}}\int_\R \varphi^2(x)e^{-x^2/2}dx<\infty.
\end{equation}
Let $H_0,H_1,\ldots$ denote the sequence of Hermite polynomials. The first few Hermite polynomials are
$H_0=1$, $H_1=X$, $H_2=X^2-1$ and $H_3=X^3-3X$. More generally, the $q$th Hermite polynomial $H_q$
 is defined through the relation $XH_q = H_{q+1}+qH_{q-1}$. It is a well-known fact
 that, when it verifies (\ref{integr}), the function $\varphi$ may be expanded in $L^2(\R,e^{-x^2/2}dx)$ (in a unique way) in terms of
 Hermite polynomials as follows:
 \begin{equation}\label{hermite-decompo}
 \varphi(x)=\sum_{q=0}^\infty a_q H_q(x).
 \end{equation}
Let $d\geq 0$ be the first integer $q\geq 0$ such that $a_q\neq 0$ in (\ref{hermite-decompo}). It is called
the {\it Hermite rank} of $\varphi$; it will play a key role in our study.
Also, let us mention the following crucial property of
Hermite polynomials with respect to Gaussian elements.
For any integer $p,q\geq 0$ and any jointly Gaussian random variables $U,V\sim \mathcal{N}(0,1)$, we have
\begin{equation}\label{orthher}
E[H_p(U)H_q(V)]=\left\{
\begin{array}{cl}
0&\mbox{if $p\neq q$}\\
q!E[UV]^q&\mbox{if $p=q$}.
\end{array}
\right.
\end{equation}
In particular (choosing $p=0$) we have that $E[H_q(X_1)]=0$ for all $q\geq 1$, meaning that $a_0=E[\varphi(X_1)]$
in (\ref{hermite-decompo}). Also, combining the decomposition (\ref{hermite-decompo}) with (\ref{orthher}),
it is straightforward to check that
\begin{equation}\label{variance}
E[\varphi^2(X_1)]=\sum_{q=0}^\infty q!a_q^2.
\end{equation}

\medskip

We are now in position to state the celebrated Breuer-Major theorem.

\begin{thm}[Breuer, Major, 1983; see \cite{BrMa}]\label{BM}
Let $\{X_k\}_{k\geq 1}$ and $\varphi:\R\to\R$ be as above.
Assume further that $a_0=E[\varphi(X_1)]=0$ and that $\sum_{k\in\Z}|\rho(k)|^d<\infty$, where $\rho$ is the covariance function of $\{X_k\}_{k\geq 1}$
and $d$ is the Hermite rank of $\varphi$ (observe that $d\geq 1$).
Then, as $n\to\infty$,
\begin{equation}\label{BMcv}
V_n=\frac{1}{\sqrt{n}}\sum_{k=1}^n \varphi(X_k)\overset{\rm law}{\to} \mathcal{N}(0,\sigma^2),
\end{equation}
with $\sigma^2$ given by
\begin{equation}\label{sigma2}
\sigma^2=\sum_{q=d}^\infty q!a_q^2\sum_{k\in\Z}\rho(k)^q\in[0,\infty).
\end{equation}
(The fact that $\sigma^2\in[0,\infty)$ is part of the conclusion.)
\end{thm}

The proof of Theorem \ref{BM} is far from being obvious.
The original proof consisted to show that {\it all} the moments of $V_n$
converge to those of the Gaussian law $\mathcal{N}(0,\sigma^2)$. As anyone might guess, this required a high
ability and a lot of combinatorics. In the proof we will offer, the complexity is the same as checking that
the variance and the fourth moment of $V_n$ converges to $\sigma^2$ and $3\sigma^4$ respectively, which is
a drastic simplification with respect to the original proof. Before doing so, let us make some other comments.
\begin{remark}
{\rm
\begin{enumerate}
\item First, it is worthwhile noticing that Theorem \ref{BM} (strictly) contains the classical central limit theorem
(CLT), which is not an evident claim at first glance.
Indeed, let $\{Y_k\}_{k\geq 1}$ be a sequence of i.i.d. centered random variables
with common variance $\sigma^2>0$, and let $F_Y$ denote the
common cumulative distribution function.
Consider the pseudo-inverse $F_Y^{-1}$ of $F_Y$, defined as
\[
F_Y^{-1}(u)={\rm inf}\{y\in\R: \,u\leq F_Y(y)\},\quad u\in(0,1).
\]
When $U\sim \mathcal{U}_{[0,1]}$ is uniformly distributed, it is well-known that $F_Y^{-1}(U)\overset{\rm law}{=}Y_1$.
Observe also that $\frac{1}{\sqrt{2\pi}}\int_{-\infty}^{X_1} e^{-t^2/2}dt$ is $\mathcal{U}_{[0,1]}$
distributed. By combining these two facts, we get that $\varphi(X_1)\overset{\rm law}{=}Y_1$ with
\[
\varphi(x)=F_Y^{-1}\left(\frac{1}{\sqrt{2\pi}}\int_{-\infty}^{x} e^{-t^2/2}dt\right), \quad x\in\R.
\]
Assume now that $\rho(0)=1$ and $\rho(k)=0$ for $k\neq 0$, that is, assume that the sequence $\{X_k\}_{k\geq 1}$ is composed of i.i.d. $\mathcal{N}(0,1)$
random variables.
Theorem \ref{BM} yields that
\[
\frac{1}{\sqrt{n}}\sum_{k=1}^n Y_k\overset{\rm law}{=}\frac{1}{\sqrt{n}}\sum_{k=1}^n \varphi(X_k)\overset{\rm law}{\to} \mathcal{N}\left(0,\sum_{q=d}^\infty q!a_q^2\right),
\]
 thereby concluding the proof
of the CLT since $\sigma^2=E[\varphi^2(X_1)]=\sum_{q=d}^\infty q!a_q^2$, see (\ref{variance}).
Of course, such a proof of the CLT is like to crack a walnut with a sledgehammer.
This approach has nevertheless its merits: it shows that the independence assumption in the CLT is not crucial to allow a Gaussian limit. Indeed, this is rather the summability of a series which is responsible of this fact,
see also the second point of this remark.
\item Assume that $d\geq 2$ and that
$\rho(k)\sim |k|^{-D}$ as $|k|\to\infty$ for some $D\in(0,\frac1d)$. In this case,
it may be shown that
$
n^{dD/2-1}\sum_{k=1}^n \varphi(X_k)
$
converges in law to a \underline{non}-Gaussian (non degenerated) random variable. This shows in particular that, in the case where $\sum_{k\in\Z}|\rho(k)|^d = \infty$,
we can get a non-Gaussian limit.
In other words, the summability assumption in Theorem \ref{BM} is, roughly speaking, equivalent (when $d\geq 2$) to  the asymptotic normality.
\item There exists a functional version of Theorem \ref{BM}, in which the sum $\sum_{k=1}^n$ is replaced
by $\sum_{k=1}^{[nt]}$ for $t\geq 0$. It is actually not that much harder to prove and, unsurprisingly, the limiting process is then the standard Brownian motion multiplied by $\sigma$.
\end{enumerate}
}
\end{remark}
\bigskip

Let us now prove Theorem \ref{BM}.
We first compute the limiting variance, which will justify the formula (\ref{sigma2}) we claim for $\sigma^2$.
Thanks to (\ref{hermite-decompo}) and (\ref{orthher}), we can write
\begin{eqnarray*}
E[V_n^2]&=&\frac{1}{n}\,E\left[\left(
\sum_{q=d}^\infty a_q \sum_{k=1}^n H_q(X_k)
\right)^2\right]
=\frac1n\sum_{p,q=d}^\infty a_pa_q \sum_{k,l=1}^n E[H_p(X_k)H_q(X_l)]\\
&=&\frac1n \sum_{q=d}^\infty q!a_q^2\sum_{k,l=1}^n \rho(k-l)^q
=\sum_{q=d}^\infty q!a_q^2\sum_{r\in\mathbb{Z}}\rho(r)^q\big(1-\frac{|r|}{n}\big){\bf 1}_{\{|r|<n\}}.
\end{eqnarray*}
When $q\geq d$ and $r\in\Z$ are fixed, we have that
\[
q!a_q^2\rho(r)^q\big(1-\frac{|r|}{n}\big){\bf 1}_{\{|r|<n\}}\to q!a_q^2\rho(r)^q\quad\mbox{as $n\to\infty$}.
\]
On the other hand, using that $|\rho(k)|=|E[X_1X_{k+1}]|\leq \sqrt{E[X_1^2]E[X_{1+k}^2]}=1$, we have
\[
q!a_q^2|\rho(r)|^q\big(1-\frac{|r|}{n}\big){\bf 1}_{\{|r|<n\}}\leq q!a_q^2|\rho(r)|^q
\leq q!a_q^2 |\rho(r)|^d,
\]
with $\sum_{q=d}^\infty \sum_{r\in\Z}  q!a_q^2 |\rho(r)|^d = E[\varphi^2(X_1)]\times \sum_{r\in\Z} |\rho(r)|^d<\infty$, see (\ref{variance}).
By applying the dominated convergence theorem, we deduce that
$E[V_n^2]\to\sigma^2$ as $n\to\infty$, with $\sigma^2\in[0,\infty)$ given by (\ref{sigma2}).

\medskip

Let us next concentrate on the proof of (\ref{BMcv}).
We shall do it in three steps of increasing generality (but of decreasing complexity!):
\begin{enumerate}
\item[(i)] when $\varphi=H_q$ has the form of a Hermite polynomial (for some $q\geq 1$);
\item[(ii)] when $\varphi=P\in\R[X]$ is a real polynomial;
\item[(iii)] in the general case when $\varphi\in L^2(\R,e^{-x^2/2}dx)$.
\end{enumerate}

We first show that (ii) implies (iii).
That is, let us assume that Theorem \ref{BM} is shown for polynomial functions $\varphi$, and let us show that it holds true for any function  $\varphi\in L^2(\R,e^{-x^2/2}dx)$.
We proceed by approximation. Let $N\geq 1$ be a (large) integer (to be chosen later) and
write
\[
V_n
=\frac{1}{\sqrt{n}} \sum_{q=d}^N a_q \sum_{k=1}^n H_q(X_k)
+\frac{1}{\sqrt{n}} \sum_{q=N+1}^\infty a_q \sum_{k=1}^n H_q(X_k)
=: V_{n,N}+R_{n,N}.
\]
Similar computations as above lead to
\begin{equation}\label{conv1}
\sup_{n\geq 1}E[R_{n,N}^2] \leq \sum_{q=N+1}^\infty q!a_q^2 \times \sum_{r\in\Z}|\rho(r)|^d \to 0\mbox{ as $N\to\infty$}.
\end{equation}
(Recall from (\ref{variance}) that $E[\varphi^2(X_1)]=\sum_{q=d}^\infty q!a_q^2<\infty$.)
On the other hand, using (ii) we have that,  for fixed $N$ and as $n\to\infty$,
\begin{equation}\label{conv2}
V_{n,N}\overset{\rm law}{\to} \mathcal{N}\left(0,\sum_{q=d}^N q!a_q^2\sum_{k\in\Z}\rho(k)^q\right).
\end{equation}
It is then a routine exercise (details are left to the reader) to deduce from (\ref{conv1})-(\ref{conv2}) that
$V_n = V_{n,N}+R_{n,N}\overset{\rm law}{\to} \mathcal{N}(0,\sigma^2)$ as $n\to\infty$, that is, that
(iii) holds true.

\medskip

Next, let us prove (i), that is, (\ref{BMcv}) when $\varphi=H_q$ is the $q$th Hermite polynomial. We actually need to work with a specific realization of the sequence $\{X_k\}_{k\geq 1}$. The space
\[
\mathcal{H}:=\overline{{\rm span}\{X_1,X_2,\ldots\}}^{L^2(\Omega)}
\]
being a real separable Hilbert space, it is isometrically isomorphic to either $\R^N$ (with $N\geq 1$) or
$L^2(\R_+)$. Let us assume that $\mathcal{H}\simeq L^2(\R_+)$, the case where $\mathcal{H}\simeq \R^N$ being easier to handle. Let $\Phi:\mathcal{H}\to L^2(\R_+)$ be an isometry. Set $e_k=\Phi(X_k)$ for each $k\geq 1$. We have
\begin{equation}
\rho(k-l)=E[X_kX_l]=\int_0^\infty e_k(x)e_l(x)dx,\quad k,l\geq 1\label{ekrho}
\end{equation}
If $B=(B_t)_{t\geq 0}$ denotes a standard Brownian motion, we deduce that
\[
\{X_k\}_{k\geq 1} \overset{\rm law}{=} \left\{\int_0^\infty e_k(t)dB_t\right\}_{k\geq 1},
\]
these two families being indeed centered, Gaussian and having the same covariance structure
(by construction of the $e_k$'s). On the other hand, it is a well-known result of stochastic analysis (which
follows from an induction argument through the It\^o formula) that, for any function $e\in L^2(\R_+)$ such
that $\|e\|_{L^2(\R_+)}=1$, we have
\begin{equation}\label{linkhermite}
H_q\left(\int_0^\infty e(t)dB_t\right) = q!\int_0^\infty dB_{t_1} e(t_1)\int_0^{t_1}dB_{t_2}e(t_2)
\ldots \int_0^{t_{q-1}} dB_{t_q} e(t_q).
\end{equation}
(For instance, by It\^o's formula we can write
\begin{eqnarray*}
\left(\int_0^\infty e(t)dB_t\right)^2 &= &2\int_0^\infty dB_{t_1} e(t_1)\int_0^{t_1} dB_{t_2} e(t_2) + \int_0^\infty e(t)^2dt\\
&=&2\int_0^\infty dB_{t_1} e(t_1)\int_0^{t_1} dB_{t_2} e(t_2) + 1,
\end{eqnarray*}
which is nothing but (\ref{linkhermite}) for $q=2$, since $H_2=X^2-1$.)
At this stage, let us adopt the two following notational conventions:
\begin{enumerate}
\item[(a)] If $\varphi$ (resp. $\psi$)
is a function of $r$ (resp. $s$) arguments,  then
the tensor product $\varphi\otimes \psi$ is the function of $r+s$ arguments given by
$\varphi\otimes \psi(x_1,\ldots,x_{r+s})=\varphi(x_1,\ldots,x_r)\psi(x_{r+1},\ldots,x_{r+s})$.
Also, if $q\geq 1$ is an integer and $e$ is a function, the tensor product function $e^{\otimes q}$
is the function $e\otimes \ldots\otimes e$ where $e$ appears $q$ times.
\item[(b)] If $f\in L^2(\R_+^q)$ is symmetric (meaning that $f(x_1,\ldots,x_q)=f(x_{\sigma(1)},\ldots,x_{\sigma(q)})$
    for all permutation $\sigma\in\mathfrak{S}_q$ and almost all $x_1,\ldots,x_q\in\R_+$)
    then \[
    I^B_q(f) = \int_{\R_+^q} f(t_1,\ldots,t_q)dB_{t_1}\ldots dB_{t_q}:=
    q!\int_0^\infty dB_{t_1}\int_0^{t_1}dB_{t_2}\ldots \int_{0}^{t_{q-1}} dB_{t_q}
    f(t_1,\ldots,t_q).
    \]
\end{enumerate}
With these new notations at hand, observe that we can rephrase (\ref{linkhermite}) in a simple way as
\begin{equation}\label{linkhermite2}
H_q\left(\int_0^\infty e(t)dB_t\right)=I^B_q(e^{\otimes q}).
\end{equation}

It is now time to
introduce a very powerful tool, the so-called {\it Fourth Moment Theorem} of Nualart and Peccati.
This wonderful result lies at the heart of the approach we shall develop in these lecture notes.
We will prove it in Section \ref{stein+malliavin}.
\begin{thm}[Nualart, Peccati, 2005; see \cite{nunugio}]\label{NP}
Fix an integer $q\geq 2$, and let $\{f_n\}_{n\geq 1}$ be a sequence of symmetric functions of $L^2(\R_+^q)$.
Assume that $E[I^B_q(f_n)^2]=q!\|f_n\|^2_{L^2(\R_+^q)}\to \sigma^2$ as $n\to\infty$ for some $\sigma>0$.
Then, the following three assertions are equivalent as $n\to\infty$:
\begin{enumerate}
\item[(1)] $I^B_q(f_n)\overset{\rm law}{\to} \mathcal{N}(0,\sigma^2)$;
\item[(2)] $E[I^B_q(f_n)^4] \overset{\rm law}{\to} 3\sigma^4$;
\item[(3)] $\|f_n\otimes_r f_n\|_{L^2(\R_+^{2q-2r})}\to 0$ for each $r=1,\ldots,q-1$,
where $f_n\otimes_r f_n$ is the function of $L^2(\R_+^{2q-2r})$ defined by
\begin{eqnarray*}
&&f_n\otimes_r f_n(x_1,\ldots,x_{2q-2r})\\
&=&\int_{\R_+^r}f_n(x_1,\ldots,x_{q-r},y_1,\ldots,y_r)
f_n(x_{q-r+1},\ldots,x_{2q-2r},y_1,\ldots,y_r)
dy_1\ldots dy_r.
\end{eqnarray*}
\end{enumerate}
\end{thm}
\begin{remark}
{\rm
In other words, Theorem \ref{NP} states that the convergence in law of a normalized sequence of multiple Wiener-It\^o integrals $I^B_q(f_n)$ towards the Gaussian law $\mathcal{N}(0,\sigma^2)$ is equivalent to convergence of just the fourth moment to $3\sigma^4$. This surprising result has been the
starting point of a new line of research, and has quickly led to several applications, extensions and improvements. One of these improvements is the following quantitative bound associated to Theorem \ref{NP} that we shall prove in
 Section \ref{stein+malliavin} by combining Stein's method with the Malliavin calculus.
\begin{thm}[Nourdin, Peccati, 2009; see \cite{NP-PTRF}]\label{NP-PTRF}
If $q\geq 2$ is an integer and $f$ is a symmetric element of $L^2(\R_+^q)$ satisfying $E[I^B_q(f)^2]=q!\|f\|^2_{L^2(\R_+^q)}=1$, then
\[
\sup_{A\subset\mathcal{B}(\R)}\left|P[I^B_q(f)\in A] - \frac{1}{\sqrt{2\pi}}\int_A e^{-x^2/2}dx\right|\leq
2\sqrt{\frac{q-1}{3q}}\sqrt{\big|E[I^B_q(f)^4]-3\big|}.
\]
\end{thm}
}
\end{remark}

\medskip

Let us go back to the proof of (i), that is, to the proof of (\ref{BMcv}) for $\varphi=H_q$.
Recall that the sequence $\{e_k\}$ has be chosen for (\ref{ekrho}) to hold.
Using (\ref{linkhermite}) (see also (\ref{linkhermite2})), we can write
$V_n=I^B_q(f_n)$, with
\[
f_n=\frac{1}{\sqrt{n}}\sum_{k=1}^n e_k^{\otimes q}.
\]
We already showed that $E[V_n^2]\to \sigma^2$ as $n\to\infty$. So, according to Theorem \ref{NP}, to get
(i) it remains to check that $\|f_n\otimes_r f_n\|_{L^2(\R_+^{2q-2r})}\to 0$ for any $r=1,\ldots,q-1$.
We have
\begin{eqnarray*}
f_n\otimes_r f_n &=& \frac1n\,\sum_{k,l=1}^n e_k^{\otimes q}\otimes_r e_l^{\otimes q} =
\frac1n\,\sum_{k,l=1}^n \langle e_k,e_l\rangle^r_{L^2(\R_+)}\,e_k^{\otimes q-r}\otimes e_l^{\otimes q-r}\\
&=&\frac1n\,\sum_{k,l=1}^n \rho(k-l)^r\,e_k^{\otimes q-r}\otimes e_l^{\otimes q-r},
\end{eqnarray*}
implying in turn
\begin{eqnarray*}
&&\|f_n\otimes_r f_n\|_{L^2(\R_+^{2q-2r})}^2 \\
&=& \frac1{n^2}\,\sum_{i,j,k,l=1}^n \rho(i-j)^r\rho(k-l)^r\langle e_i^{\otimes q-r}\otimes e_j^{\otimes q-r},e_k^{\otimes q-r}\otimes e_l^{\otimes q-r}\rangle_{L^2(\R_+^{2q-2r})}\\
&=& \frac1{n^2}\,\sum_{i,j,k,l=1}^n \rho(i-j)^r\rho(k-l)^r\rho(i-k)^{q-r}\rho(j-l)^{q-r}.
\end{eqnarray*}
Observe that $|\rho(k-l)|^r|\rho(i-k)|^{q-r}\leq |\rho(k-l)|^q+|\rho(i-k)|^{q}$. This, together with other
obvious manipulations, leads to the bound
\begin{eqnarray*}
\|f_n\otimes_r f_n\|_{L^2(\R_+^{2q-2r})}^2
&\leq& \frac2{n}\,\sum_{k\in\Z} |\rho(k)|^q \sum_{|i|<n}|\rho(i)|^{r} \sum_{|j|<n} |\rho(j)|^{q-r}\\
&\leq& \frac2{n}\,\sum_{k\in\Z} |\rho(k)|^d \sum_{|i|<n}|\rho(i)|^{r} \sum_{|j|<n}|\rho(j)|^{q-r}\\
&=&2\sum_{k\in\Z} |\rho(k)|^d \times n^{-\frac{q-r}{q}}\sum_{|i|<n}|\rho(i)|^{r} \times n^{-\frac{r}{q}}\sum_{|j|<n}|\rho(j)|^{q-r}.
\end{eqnarray*}
Thus, to get that $\|f_n\otimes_r f_n\|_{L^2(\R_+^{2q-2r})}\to 0$ for any $r=1,\ldots,q-1$,
it suffices to show that
\[
s_n(r):=n^{-\frac{q-r}{q}}\sum_{|i|<n}|\rho(i)|^{r}\to 0\,\mbox{ for any $r=1,\ldots,q-1$}.
\]
Let $r=1,\ldots,q-1$. Fix $\delta\in(0,1)$ (to be chosen later) and let us decompose $s_n(r)$ into
\[
s_n(r)=n^{-\frac{q-r}{q}}\sum_{|i|<[n\delta]}|\rho(i)|^{r}  + n^{-\frac{q-r}{q}}\sum_{[n\delta]\leq |i|<n}|\rho(i)|^{r} =: s_{1,n}(\delta,r)+s_{2,n}(\delta,r).
\]
Using H\"older inequality, we get that
\[
s_{1,n}(\delta,r)\leq n^{-\frac{q-r}{r}}\left(\sum_{|i|<[n\delta]}|\rho(i)|^q\right)^{r/q}(1+2[n\delta])^{\frac{q-r}{q}}\leq {\rm cst}\times \delta^{1-r/q},
\]
as well as
\[
s_{2,n}(\delta,r)\leq n^{-\frac{q-r}{r}}\left(\sum_{[n\delta]\leq |i|<n}|\rho(i)|^q\right)^{r/q}(2n)^{\frac{q-r}{q}}\leq {\rm cst}\times \left(\sum_{|i|\geq [n\delta]}|\rho(i)|^q\right)^{r/q}.
\]
Since $1-r/q>0$, it is a routine exercise (details are left to the reader) to deduce that $s_n(r)\to 0$
as $n\to\infty$. Since this is true for any $r=1,\ldots,q-1$, this concludes the proof of (i).

\medskip

It remains to show (ii), that is, convergence in law (\ref{BMcv}) whenever $\varphi$ is a real polynomial. We shall use the multivariate counterpart of Theorem \ref{NP}, which
was obtained shortly afterwards by Peccati and Tudor. Since only a weak version
(where all the involved multiple Wiener-It\^o integrals have different orders) is needed here, we
state the result of Peccati and Tudor only in this situation.
We refer to Section \ref{s:smartpath} for a more general version and its proof.

\begin{thm}[Peccati, Tudor, 2005; see \cite{PTu04}]\label{PT}
Consider $l$ integers $q_1,\ldots,q_l\geq 1$, with $l\geq 2$. Assume that all the $q_i$'s are pairwise different.
For each $i=1,\ldots,l$, let $\{f^i_n\}_{n\geq 1}$ be a sequence of symmetric functions of $L^2(\R_+^{q_i})$
satisfying $E[I^B_{q_i}(f^i_n)^2]=q_i!\|f^i_n\|^2_{L^2(\R_+^{q_i})}\to \sigma_i^2$ as $n\to\infty$ for some $\sigma_i>0$.
Then, the following two assertions are equivalent as $n\to\infty$:
\begin{enumerate}
\item[(1)] $I^B_{q_i}(f^i_n)\overset{\rm law}{\to} \mathcal{N}(0,\sigma_i^2)$ for all $i=1,\ldots,l$;
\item[(2)] $\big(I^B_{q_1}(f^1_n),\ldots,I^B_{q_l}(f^l_n)\big)\overset{\rm law}{\to} \mathcal{N}\big(0,{\rm diag}(\sigma_1^2,\ldots,\sigma_l^2)\big)$.
\end{enumerate}
\end{thm}
In other words, Theorem \ref{PT} proves the surprising fact that, for such a sequence of vectors of multiple Wiener-It\^o integrals, componentwise convergence to Gaussian
always implies joint convergence. We shall combine Theorem \ref{PT} with (i) to prove (ii). Let $\varphi$ have the form of a real
polynomial. In particular, it admits a decomposition of the type $\varphi=\sum_{q=d}^N a_q H_q$ for some {\it finite} integer $N\geq d$.
Together with (i), Theorem \ref{PT} yields that
\[
\left(
\frac{1}{\sqrt{n}}\sum_{k=1}^n H_d(X_k),\ldots,\frac{1}{\sqrt{n}}\sum_{k=1}^n H_N(X_k)
\right)
\overset{\rm law}{\to}
\mathcal{N}\big(0,{\rm diag}(\sigma_d^2,\ldots,\sigma_N^2)\big),
\]
where $\sigma_q^2=q!\sum_{k\in \Z}\rho(k)^q$, $q=d,\ldots,N$.
We deduce that \[
V_n= \frac{1}{\sqrt{n}}\sum_{q=d}^N a_q \sum_{k=1}^n H_q(X_k) \overset{\rm law}{\to}
\mathcal{N}\left(0,\sum_{q=d}^N a_q^2 q! \sum_{k\in \Z}\rho(k)^q\right),
\]
which is the desired conclusion in (ii) and conclude the proof of Theorem \ref{BM}. \qed

\bigskip

{\bf To go further}.
In \cite{NoPePo}, one associates quantitative bounds to Theorem \ref{BM} by using a similar approach.

\noindent
\section{Universality of Wiener chaos}

Before developing the material which will be necessary for the proof of the Fourth Moment Theorem \ref{NP}
(as well as other related results), to motivate the reader let us study yet another consequence of this beautiful result.

For {\it any} sequence $X_1,X_2,\ldots$ of i.i.d. random variables
with mean 0 and variance 1, the central limit theorem asserts that
$V_n=(X_1+\ldots+X_n)/\sqrt{n}\overset{\rm law}{\to} \mathcal{N}(0,1)$ as $n\to\infty$.
It is a particular instance of what is commonly referred to as a `universality phenomenon'
in probability.
Indeed, we observe that the limit of the sequence $V_n$ does not rely on the specific law of the $X_i$'s, but only
of the fact that its first two moments are 0 and 1 respectively.

Another example that exhibits a universality phenomenon is given by Wigner's theorem in the  random matrix theory.
More precisely, let $\{X_{ij}\}_{j>i\geq 1}$ and $\{X_{ii}/\sqrt{2}\}_{i\geq 1}$ be two independent families composed
of i.i.d.\! random variables with mean 0, variance 1, and all the moments. Set $X_{ji}=X_{ij}$
and consider the $n\times n$ random matrix $M_n=(\frac{X_{ij}}{\sqrt{n}})_{1\leq i,j\leq n}$.
The matrix $M_n$ being symmetric,
its eigenvalues $\lambda_{1,n},\ldots,\lambda_{n,n}$ (possibly repeated with multiplicity) belong to $\R$.
Wigner's theorem then asserts that the spectral measure of $M_n$, that
is, the random probability measure defined as $\frac1n\sum_{k=1}^n \delta_{\lambda_{k,n}}$, converges almost surely to the semicircular law $\frac{1}{2\pi}\sqrt{4-x^2}{\bf 1}_{[-2,2]}(x)dx$, whatever the exact distribution of the
entries of $M_n$ are.

In this section, our aim is to prove yet another universality phenomenon, which is
in the spirit of the two afore-mentioned
results. To do so, we need to introduce the following two blocks of basic ingredients:
\begin{enumerate}
\item[(i)] Three sequences ${\bf X}=(X_1,X_2,\ldots)$,
${\bf G}=(G_1,G_2,\ldots)$ and ${\bf E}=(\varepsilon_1,\varepsilon_2,\ldots)$
 of i.i.d.\! random variables, all with mean 0, variance 1 and finite fourth moment.
 We are more specific with ${\bf G}$ and ${\bf E}$, by assuming further that $G_1\sim\mathcal{N}(0,1)$ and $P(\varepsilon_1=1)=P(\varepsilon_1=-1)=1/2$.
 (As we will see, ${\bf E}$ will actually play no role in the statement of Theorem \ref{NPR}; we will however use it to build a interesting counterexample, see Remark \ref{rk221}(1).)
 \item[(ii)] A fixed integer $d\geq 1$ as well as a sequence $g_n:\{1,\ldots,n\}^d\to \R$, $n\geq 1$
 of real functions, each $g_n$ satisfying in addition that, for all $i_1,\ldots,i_d=1,\ldots,n$,
 \begin{enumerate}
 \item[(a)] $g_n(i_1,\ldots,i_d)=g_n(i_{\sigma(1)},\ldots,i_{\sigma(d)})$ for all permutation $\sigma\in\mathfrak{S}_d$;
 \item[(b)] $g_n(i_1,\ldots,i_d)=0$ whenever $i_k=i_l$ for some $k\neq l$;
 \item[(c)] $d!\sum_{i_1,\ldots,i_d=1}^n g_n(i_1,\ldots,i_d)^2=1$.
 \end{enumerate}
 (Of course, conditions $(a)$ and $(b)$ are becoming immaterial when $d=1$.)
 If ${\bf x}=(x_1,x_2,\ldots)$ is a given real sequence, we also set
 \[
 Q_d(g_n,{\bf x})=
 \sum_{i_1,\ldots,i_d=1}^n g_n(i_1,\ldots,i_d)x_{i_1}\ldots x_{i_d}.
 \]
 Using (b) and (c), it is straightforward to check that, for any $n\geq 1$, we have
 $E[Q_d(g_n,{\bf X})]=0$ and $E[Q_d(g_n,{\bf X})^2]=1$.
\end{enumerate}

We are now in position to state our new universality phenomenon.

\begin{thm}[Nourdin, Peccati, Reinert, 2010; see \cite{NPR}]\label{NPR}
Assume that $d\geq 2$. Then, as $n\to\infty$, the following two assertions are equivalent:
\begin{enumerate}
\item[($\alpha$)] $Q_d(g_n,{\bf G})\overset{\rm law}{\to}\mathcal{N}(0,1)$;
\item[($\beta$)] $Q_d(g_n,{\bf X})\overset{\rm law}{\to}\mathcal{N}(0,1)$ for \underline{any} sequence ${\bf X}$ as given in (i).
\end{enumerate}
\end{thm}

Before proving Theorem \ref{NPR}, let us address some comments.
\begin{remark}\label{rk221}
{\rm
\begin{enumerate}
\item In reality, the universality phenomenon in Theorem \ref{NPR} is a bit more subtle than in the CLT or in Wigner's theorem.
To illustrate what we have in mind, let us consider an explicit situation (in the case $d=2$).
Let $g_n:\{1,\ldots,n\}^2\to\R$ be the function given by
\[
g_n(i,j)=\frac{1}{2\sqrt{n-1}}{\bf 1}_{\{i=1,j\geq 2\mbox{ or }j=1,i\geq 2\}}.
\]
It is easy to check that $g_n$ satisfies the three assumptions $(a)$-$(b)$-$(c)$ and also that
\[
Q_2(g_n,{\bf x})=x_1\times \frac{1}{\sqrt{n-1}}\sum_{k=2}^n x_k.
\]
The classical CLT then implies that
$Q_2(g_n,{\bf G})\overset{\rm law}{\to}G_1G_2$ and
$Q_2(g_n,{\bf E})\overset{\rm law}{\to}\varepsilon_1G_2$.
Moreover, it is a classical and easy exercise to check that $\varepsilon_1G_2$ is $\mathcal{N}(0,1)$ distributed.
Thus, what we just showed is that, although $Q_2(g_n,{\bf E})\overset{\rm law}{\to}\mathcal{N}(0,1)$ as $n\to\infty$, the assertion ($\beta$) in Theorem \ref{NPR}
fails when choosing ${\bf X}={\bf G}$ (indeed, the product of two independent $\mathcal{N}(0,1)$ random variables is not gaussian).
This means that, in Theorem \ref{NPR}, we cannot replace the sequence ${\bf G}$ in ($\alpha$) by another sequence
(at least, not by ${\bf E}$ !).
\item Theorem \ref{NPR} is completely false when $d=1$. For an explicit counterexample, consider for instance
$g_n(i)={\bf 1}_{\{i=1\}}$, $i=1,\ldots,n$. We then have $Q_1(g_n,{\bf x})=x_1$. Consequently, the assertion
($\alpha$) is trivially verified (it is even an equality in law!) but the assertion ($\beta$) is never true unless $X_1\sim\mathcal{N}(0,1)$.
\end{enumerate}
}
\end{remark}

\bigskip
\noindent
{\it Proof of Theorem \ref{NPR}}.
Of course, only the implication ($\alpha$)$\to$($\beta$)
must be shown. Let us divide its proof into three steps.

\medskip

\underline{\it Step 1}. Set $e_i={\bf 1}_{[i-1,i]}$, $i\geq 1$, and
let $f_n\in L^2(\R_+^d)$ be the symmetric function defined as
\[
f_n = \sum_{i_1,\ldots,i_d=1}^n g_n(i_1,\ldots,i_d)e_{i_1}\otimes\ldots\otimes e_{i_d}.
\]
By the very definition of $I^B_d(f_n)$, we have
\[
I^B_d(f_n)=d!\sum_{i_1,\ldots,i_d=1}^n g_n(i_1,\ldots,i_d)\int_0^\infty dB_{t_1}e_{i_1}(t_1)\int_0^{t_1}dB_{t_2}e_{i_2}(t_2)\ldots\int_0^{t_{d-1}}dB_{t_d}e_{i_d}(t_d).
\]
Observe that
\[
\int_0^\infty dB_{t_1}e_{i_1}(t_1)\int_0^{t_1}dB_{t_2}e_{i_2}(t_2)\ldots\int_0^{t_{d-1}}dB_{t_d}e_{i_d}(t_d)
\]
is {\it not} almost surely zero (if and) only if $i_d\leq i_{d-1}\leq\ldots\leq i_1$. By combining this fact with assumption (b), we deduce
that
\begin{eqnarray*}
I^B_d(f_n)&=&d!\sum_{1\leq i_d<\ldots<i_1\leq n} g_n(i_1,\ldots,i_d)\int_0^\infty dB_{t_1}e_{i_1}(t_1)\int_0^{t_1}dB_{t_2}e_{i_2}(t_2)\ldots\int_0^{t_{d-1}}dB_{t_d}e_{i_d}(t_d)\\
&=&d!\sum_{1\leq i_d<\ldots<i_1\leq n} g_n(i_1,\ldots,i_d)(B_{i_1}-B_{i_1-1})\ldots (B_{i_d}-B_{i_d-1})\\
&=&\sum_{i_1,\ldots,i_d=1}^n g_n(i_1,\ldots,i_d)(B_{i_1}-B_{i_1-1})\ldots (B_{i_d}-B_{i_d-1})
\overset{\rm law}{=}Q_d(g_n,{\bf G}).
\end{eqnarray*}
That is, the sequence $Q_d(g_n,{\bf G})$ in ($\alpha$) has actually the form of a multiple Wiener-It\^o integral.
On the other hand, going back to the definition of $f_n\otimes_{d-1}f_n$ and using that $\langle e_i,e_j\rangle_{L^2(\R_+)}=\delta_{ij}$
(Kronecker symbol), we get
\[
f_n\otimes_{d-1}f_n = \sum_{i,j=1}^n \left(\sum_{k_2,\ldots,k_d=1}^n g_n(i,k_2,\ldots,k_d)g_n(j,k_2,\ldots,k_d)\right)e_i\otimes
e_j,
\]
so that
\begin{eqnarray}
\left\|f_n\otimes_{d-1}f_n\right\|^2_{L^2(\R_+^2)}& =& \sum_{i,j=1}^n \left(\sum_{k_2,\ldots,k_d=1}^n g_n(i,k_2,\ldots,k_d)g_n(j,k_2,\ldots,k_d)\right)^2\notag\\
&\geq &\sum_{i=1}^n \left(\sum_{k_2,\ldots,k_d=1}^n g_n(i,k_2,\ldots,k_d)^2\right)^2\quad\mbox{(by summing only over
$i=j$)}\notag\\
&\geq&\max_{1\leq i\leq n} \left(\sum_{k_2,\ldots,k_d=1}^n g_n(i,k_2,\ldots,k_d)^2\right)^2=\tau_n^2,\label{crucial}
\end{eqnarray}
where
\begin{equation}\label{taun}
\tau_n:=\max_{1\leq i\leq n} \sum_{k_2,\ldots,k_d=1}^n g_n(i,k_2,\ldots,k_d)^2.
\end{equation}
Now, assume that ($\alpha$) holds. By Theorem \ref{NP} and because $Q_d(g_n,{\bf G})\overset{\rm law}{=}I^B_d(f_n)$,
 we have in particular that
$\|f_n\otimes_{d-1} f_n\|_{L^2(\R_+^2)}\to 0$ as $n\to\infty$.
Using the inequality (\ref{crucial}), we deduce that $\tau_n\to 0$ as $n\to\infty$.

\medskip

\underline{\it Step 2}. We claim that the following result (whose proof is given in Step 3) allows to conclude the proof of $(\alpha)\to(\beta)$.

\begin{thm}[Mossel, O'Donnel, Oleszkiewicz, 2010; see \cite{MOO}]\label{MOO}
Let ${\bf X}$ and ${\bf G}$ be given as in (i) and let $g_n:\{1,\ldots,n\}^d\to\R$ be a function
satisfying the three conditions (a)-(b)-(c).
Set $\gamma=\max\{3,E[X_1^4]\}\geq 1$
and let $\tau_n$ be the quantity given by (\ref{taun}). Then, for all function $\varphi:\R\to\R$
of class $\mathcal{C}^3$ with $\|\varphi'''\|_\infty<\infty$, we have
\[
\big|E[\varphi(Q_d(g_n,{\bf X}))] - E[\varphi(Q_d(g_n,{\bf G}))]
\big|
\leq
  \frac{\gamma}3(3+2\gamma)^{\frac32(d-1)}
 d^{3/2}\sqrt{d!}\,
 \|\varphi'''\|_\infty\, \sqrt{\tau_n}.
\]
\end{thm}

\medskip

Indeed, assume that ($\alpha$) holds. By Step 1, we have that $\tau_n\to 0$ as $n\to\infty$.
Next, Theorem \ref{MOO} together with ($\alpha$),
lead to ($\beta$) and therefore conclude the proof of Theorem \ref{NPR}.

\bigskip

\underline{\it Step 3}: Proof of Theorem \ref{MOO}. During the proof, we will need the following auxiliary lemma,
which is of independent interest.

\begin{lemma}[Hypercontractivity]\label{hyphyphyp...popotame!}
Let $n\geq d\geq 1$, and consider a multilinear polynomial $P\in\R[x_1,\ldots,x_n]$ of degree $d$, that is, $P$ is of the form
\[
P(x_1,\ldots,x_n)=
\sum_{\substack{S\subset\{1,\ldots,n\}\\|S|=d}}
a_S \prod_{i\in S} x_i.
\]
Let ${\bf X}$ be as in (i). Then,
\begin{equation}\label{ascona4}
E\left[P(X_1,\ldots,X_n)^4\right]\leq \left(3+2E[X_1^4]\right)^{2d}E\left[P(X_1,\ldots,X_n)^2\right]^2.
\end{equation}
\end{lemma}
{\it Proof}.
The proof follows ideas from \cite{MOO} is by induction on $n$.
The case $n=1$ is trivial. Indeed, in this case we have $d=1$ so that
$P(x_1)=ax_1$; the conclusion therefore asserts that (recall that $E[X_1^2]=1$, implying in turn
that $E[X_1^4]\geq E[X_1^2]^2=1$)
\[
a^4 E[X_1^4]\leq
a^4\left(3+2E[X_1^4]\right)^{2},
\]
which is evident.
Assume now that $n\geq 2$. We can write
\[
P(x_1,\ldots,x_n)=R(x_1,\ldots,x_{n-1})+x_nS(x_1,\ldots,x_{n-1}),
\]
where $R,S\in\R[x_1,\ldots,x_{n-1}]$ are multilinear polynomials of $n-1$ variables.
Observe that $R$ has degree $d$, while $S$ has degree $d-1$.
Now write $\mathbf{P}=P(X_1,\ldots,X_n)$,  $\mathbf{R}=R(X_1,\ldots,X_{n-1})$,
$\mathbf{S}=S(X_1,\ldots,X_{n-1})$ and $\alpha=E[X_1^4]$.
Clearly, $\mathbf{R}$ and $\mathbf{S}$ are independent of $X_n$.
We have, using $E[X_n]=0$ and $E[X_n^2]=1$:
\begin{eqnarray*}
E[\mathbf{P}^2]=E[(\mathbf{R}+\mathbf{S}X_n)^2]
&=&E[\mathbf{R}^2]+E[\mathbf{S}^2]\\
E[\mathbf{P}^4]=E[(\mathbf{R}+\mathbf{S}X_n)^4]
&=&E[\mathbf{R}^4]+6E[\mathbf{R}^2\mathbf{S}^2]+4E[X_n^3]E[\mathbf{R}\mathbf{S}^3]+E[X_n^4]E[\mathbf{S}^4].
\end{eqnarray*}
Observe that $E[\mathbf{R}^2\mathbf{S}^2]\leq \sqrt{E[{\bf R}^4]}\sqrt{E[{\bf S}^4]}$
and
\[
E[X_n^3]E[\mathbf{R}\mathbf{S}^3]\leq \alpha^{\frac34}
\big(E[{\bf R}^4]\big)^{\frac14}\big(E[{\bf S}^4]\big)^{\frac34}
\leq \alpha \sqrt{E[{\bf R}^4]}\sqrt{E[{\bf S}^4]} +\alpha E[{\bf S}^4],
\]
where the last inequality used both $x^{\frac14} y^{\frac34}\leq \sqrt{xy}+y$
(by considering $x<y$ and $x>y$) and $\alpha^\frac34\leq \alpha$ (because
$\alpha\geq E[X_n^4]\geq E[X_n^2]^2=1$).
Hence
\begin{eqnarray*}
E[{\bf P}^4]&\leq& E[{\bf R}^4]  +2(3+2\alpha)\sqrt{E[{\bf R}^4]}\sqrt{E[{\bf S}^4]}+5\alpha E[{\bf S}^4]\\
&\leq&E[{\bf R}^4] +2(3+2\alpha)\sqrt{E[{\bf R}^4]}\sqrt{E[{\bf S}^4]}+(3+2\alpha)^2 E[{\bf S}^4]\\
&=& \left( \sqrt{E[{\bf R}^4]}+(3+2\alpha)\sqrt{E[{\bf S}^4]}\right)^2.
\end{eqnarray*}
By induction, we have $ \sqrt{E[{\bf R}^4]} \leq (3+2\alpha)^d E[{\bf R}^2]$
and $ \sqrt{E[{\bf S}^4]} \leq (3+2\alpha)^{d-1} E[{\bf S}^2]$.
Therefore
\[
E[{\bf P}^4]\leq (3+2\alpha)^{2d}\left(E[{\bf R}^2]+E[{\bf S}^2]\right)^2 =  (3+2\alpha)^{2d} E[{\bf P}^2]^2,
\]
and the proof of the lemma is concluded.
\qed

\medskip

We are now in position to prove Theorem \ref{MOO}.
Following \cite{MOO}, we use the Lindeberg replacement trick.
Without loss of generality, we assume that ${\bf X}$ and ${\bf G}$ are stochastically
independent.
For $i=0,\ldots,n$, let ${\bf W}^{(i)}=(G_1,\ldots,G_i,X_{i+1},\ldots,X_n)$.
Fix a particular $i=1,\ldots,n$ and write
\begin{eqnarray*}
U_{i}&=&\sum_{
\stackrel{1\leq i_1,\ldots,i_{d}\leq n}{i_1\neq i,\ldots,i_{d}\neq i}
}
 g_n(i_1,\ldots,i_{d})W_{i_1}^{(i)}\ldots W_{i_{d}}^{(i)},\\
V_{i}&=&\sum_{
\stackrel{1\leq i_1,\ldots,i_{d}\leq n}{\exists j:\,i_j= i}
}
 g_n(i_1,\ldots,i_{d})W_{i_1}^{(i)}\ldots \widehat{W_{i}^{(i)}}\ldots W_{i_{d}}^{(i)}\\
 &=&d\sum_{i_2,\ldots,i_d=1}^n
 g_n(i,i_2,\ldots,i_{d})W_{i_2}^{(i)}\ldots W_{i_{d}}^{(i)},
\end{eqnarray*}
where $\widehat{W_i^{(i)}}$ means that this particular term is dropped
(observe that this notation bears no ambiguity: indeed, since $g_n$ vanishes on diagonals,
each string $i_1,\ldots,i_{d}$ contributing to the definition of $V_{i}$ contains the symbol
$i$ exactly once).
For each $i$, note that $U_i$ and $V_i$ are independent of
the variables $X_i$ and $G_i$, and that
\[
Q_d(g_n,{\bf W}^{(i-1)})=U_{i}+X_iV_{i}\quad\mbox{and}\quad
Q_d(g_n,{\bf W}^{(i)})=U_{i}+G_iV_{i}.
\]
By Taylor's theorem, using the independence of $X_i$ from $U_i$ and $V_i$, we have
\begin{eqnarray*}
&&\left| E\big[ \varphi(U_i+X_iV_i)\big]-
E\big[\varphi(U_i)\big]
-E\big[\varphi'(U_i)V_i\big]E[X_i]
-\frac12E\big[\varphi''(U_i)V_i^2\big]E[X_i^2]
\right|   \\
&\leq&
\frac16\|\varphi'''\|_\infty
E[|X_i|^3]E[|V_i|^3].
\end{eqnarray*}
Similarly,
\begin{eqnarray*}
&&\left| E\big[ \varphi(U_i+G_iV_i)\big]-
E\big[\varphi(U_i)\big]
-E\big[\varphi'(U_i)V_i\big]E[G_i]
-\frac12E\big[\varphi''(U_i)V_i^2\big]E[G_i^2]
\right|   \\
&\leq&
\frac16\|\varphi'''\|_\infty
E[|G_i|^3]E[|V_i|^3].
\end{eqnarray*}
Due to the matching moments up to second order on one hand, and using that
$E[|X_i|^3]\leq \gamma$
and $E[|G_i|^3]\leq \gamma$ on the other hand, we obtain that
\begin{eqnarray*}
&&\left|E\big[ \varphi(Q_d (g_n,{\bf W}^{(i-1)}))\big] - E\big[ \varphi( Q_d (g_n,{\bf W}^{(i)}))\big]\right|
=\left|E\big[ \varphi(U_i+G_i V_i)\big] - E\big[ \varphi( U_i+X_i V_i)\big]\right|\\
&\leq& \frac{\gamma}3
\|\varphi'''\|_\infty
   E[|V_i|^3].
\end{eqnarray*}
By Lemma \ref{hyphyphyp...popotame!},
we have
\[
E[ |V_{i}|^3]\leq E[V_{i}^4]^{\frac34}\leq
(3+2\gamma)^{\frac32(d-1)}E[ V_{i}^2]^\frac32.
\]
Using the independence between ${\bf X}$ and ${\bf G}$,
the properties of $g_n$ (which is symmetric and vanishes on diagonals)
as well as $E[X_i]=E[G_i]=0$ and
$E[X_i^2]=E[G_i^2]=1$, we get
\begin{eqnarray*}
E[V_{i}^2]^{3/2}&=& \bigg(
dd!\sum_{i_2,\ldots,i_d=1}^n
 g_n(i,i_2,\ldots,i_{d})^2
\bigg)^{3/2}\\
&\leq& (dd!)^{3/2}\sqrt{ \max_{1\leq j\leq n}
\sum_{j_2,\ldots,j_{d}=1}^n
 g_n(j,j_2,\ldots,j_{d})^2}\times \sum_{i_2,\ldots,i_{d}=1}^n
 g_n(i,i_2,\ldots,i_{d})^2,
\end{eqnarray*}
implying in turn that
\begin{eqnarray*}
\sum_{i=1}^n E[V_{i}^2]^{3/2} &\leq& (dd!)^{3/2}
\sqrt{ \max_{1\leq j\leq n}
\sum_{j_2,\ldots,j_{d}=1}^n
 g_n(j,j_2,\ldots,j_{d_k})^2}
 \times \sum_{i_1,\ldots,i_{d}=1}^n
 g_n(i_1,i_2,\ldots,i_{d})^2,
 \\
 &=&  d^{3/2}\sqrt{d!}\, \sqrt{\tau_n}.
\end{eqnarray*}
By collecting the previous bounds, we get
\begin{eqnarray*}
&& \left|
E[\varphi(Q_d(g_n,{\bf X}))]
-E[\varphi(Q_d(g_n,{\bf G}))]
\right|\\
&&\leq
\sum_{i=1}^n \left|
E\big[ \varphi(Q_d (g_n,{\bf W}^{(i-1)}))\big] - E\big[ \varphi( Q_d (g_n,{\bf W}^{(i)}))\big]
\right|\\
&&\leq \frac{\gamma}3
\|\varphi'''\|_\infty
   \sum_{i=1}^n E[|V_i|^3]
   \leq
    \frac{\gamma}3(3+2\gamma)^{\frac32(d-1)}
\|\varphi'''\|_\infty
   \sum_{i=1}^n
   E[ V_{i}^2]^\frac32
   \\
 & &\leq
    \frac{\gamma}3(3+2\gamma)^{\frac32(d-1)}
 d^{3/2}\sqrt{d!}\,
 \|\varphi'''\|_\infty
  \sqrt{\tau_n}.
\end{eqnarray*}
\qed

As a final remark, let us observe that Theorem \ref{MOO} contains the CLT as a special case.
Indeed, fix $d=1$ and let $g_n:\{1,\ldots,n\}\to\R$ be the function given by $g_n(i)=\frac1{\sqrt{n}}$.
We then have $\tau_n=1/n$. It is moreover clear that $Q_1(g_n,{\bf G})\sim \mathcal{N}(0,1)$.
Then, for any function  $\varphi:\R\to\R$ of class $\mathcal{C}^3$ with $\|\varphi'''\|_\infty<\infty$ and any sequence ${\bf X}$ as in (i), Theorem \ref{MOO} implies that
\[
\left|E\left[\varphi\left(\frac{X_1+\ldots+X_n}{\sqrt{n}}\right)\right]
-\frac{1}{\sqrt{2\pi}}\int_\R \varphi(y)e^{-y^2/2}dy
\right|\leq \max\{E[X_1^4]/3,1\}\|\varphi'''\|_\infty,
\]
from which it is straightforward to deduce the CLT.

\bigskip

{\bf To go further}.
In \cite{NPR}, Theorem \ref{NPR} is extended to the case where the target law is the centered Gamma law.
In \cite{peccatizheng}, there is a version of Theorem \ref{NPR} in which the sequence ${\bf G}$ is replaced
by ${\bf P}$, a sequence of i.i.d. Poisson random variables. Finally, let us mention that both
Theorems \ref{NPR} and \ref{MOO} have been extended to the free probability framework (see Section \ref{sec:free}) in the reference \cite{deyanourdin}.

\noindent
\section{Stein's method}\label{sec:stein}

In this section, we shall introduce some basic features of the so-called Stein method, which is the first step toward the proof
of the Fourth Moment Theorem \ref{NP}. Actually, we will not need the full force of this method, only a basic estimate.

A random variable $X$ is $\mathcal{N}(0,1)$ distributed if and only if
$E[e^{itX}]=e^{-t^2/2}$ for all $t\in\R$. This simple fact leads to the idea that a random variable
$X$ has a law which is {\it close} to $\mathcal{N}(0,1)$ if and only if $E[e^{itX}]$ is {\it approximately}
$e^{-t^2/2}$ for all $t\in\R$. This last claim is nothing but the usual criterion
for the convergence in law through the use of characteristic functions.

Stein's seminal idea is somehow similar. He noticed in \cite{Stein_orig} that $X$ is $\mathcal{N}(0,1)$ distributed
if and only if $E[f'(X)-Xf(X)]=0$ for all function $f$ belonging to a sufficiently rich class of functions
(for instance, the functions which are $\mathcal{C}^1$ and whose derivative grows at most polynomially).
He then wondered whether a suitable quantitative version of this identity may have fruitful consequences.
This is actually the case and, even for specialists (at least for me!), the reason why it works so well remains a bit mysterious.
Surprisingly, the simple following statement (due to Stein \cite{Stein_orig}) happens to contain all the elements of Stein's method that are needed for our discussion.
(For more details or extensions of the method, one can
consult the recent books \cite{goldstein,book-malliavinstein} and the references therein.)

\begin{lemma}[Stein, 1972; see \cite{Stein_orig}]\label{stein}
Let $N\sim \mathcal{N}(0,1)$ be a standard Gaussian random variable.
Let $h:\R\to [0,1]$ be any continuous function.
Define
$f:\mathbb{R}\to\mathbb{R}$ by
\begin{eqnarray}\label{fz}
f(x)&=&e^{\frac{x^2}{2}}\int_{-\infty}^x \big(h(a)-E[h(N)]\big)
e^{-\frac{a^2}{2}}da\\
&=&-e^{\frac{x^2}{2}}\int_{x}^\infty \big(h(a)-E[h(N)]\big)
e^{-\frac{a^2}{2}}da.\label{fz2}
\end{eqnarray}
Then $f$ is of class $\mathcal{C}^1$, and satisfies
$|f(x)|\leq \sqrt{\pi/2}$, $|f'(x)|\leq 2$ and
\begin{equation}\label{stein-equation}
f'(x)=xf(x)+h(x)-E[h(N)]
\end{equation}
for all $x\in\R$.
\end{lemma}
\noindent{\it Proof}:
The equality between (\ref{fz}) and (\ref{fz2}) comes from
\[
0=E\big[h(N)-E[h(N)]\big]
=\frac{1}{\sqrt{2\pi}}\int_{-\infty}^{+\infty}
\big(h(a)-E[h(N)]\big)
e^{-\frac{a^2}{2}}da.
\]
Using (\ref{fz2}) we have, for $x\geq 0$:
\begin{eqnarray*}
\big|xf(x)\big|
&=&\left|xe^{\frac{x^2}{2}}\int_x^{+\infty}\big(h(a)-E[h(N)]\big)
e^{-\frac{a^2}{2}}da\right|\\
&\leq&xe^{\frac{x^2}{2}}\int_x^{+\infty}
e^{-\frac{a^2}{2}}da\leq e^{\frac{x^2}{2}}\int_x^{+\infty}
ae^{-\frac{a^2}{2}}da=1.
\end{eqnarray*}
Using (\ref{fz}) we have, for $x\leq 0$:
\begin{eqnarray*}
\big|xf(x)\big|
&=&\left|xe^{\frac{x^2}{2}}\int_{-\infty}^x \big(h(a)-E[h(N)]\big)
e^{-\frac{a^2}{2}}da\right|\\
&\leq&|x|e^{\frac{x^2}{2}}\int_{|x|}^{+\infty}
e^{-\frac{a^2}{2}}da
\leq e^{\frac{x^2}{2}}\int_{|x|}^{+\infty}
ae^{-\frac{a^2}{2}}da=1.
\end{eqnarray*}
The identity (\ref{stein-equation}) is readily checked.
We deduce, in particular, that
\[
|f'(x)|\leq |xf(x)|+|h(x)-E[h(N)]|\leq 2
\]
for all $x\in\R$.
On the other hand, by (\ref{fz})-(\ref{fz2}), we have, for every $x\in\R$,
\[
|f (x)|\leq e^{x^2/2}\min\left(\int_{-\infty}^x e^{-y^2/2}dy, \int_{x}^\infty e^{-y^2/2}dy \right)
=  e^{x^2/2}\int_{|x|}^\infty e^{-y^2/2}dy \leq \sqrt{\frac{\pi}{2}},
\]
where the last inequality is obtained by observing that the function $s:\R_+\to\R$ given by
$
s(x)= e^{x^2/2}\!\int_{x}^\infty \! e^{-y^2/2}dy
$
attains its maximum at $x=0$ (indeed, we have
\[
s'(x)=xe^{x^2/2}\int_x^\infty e^{-y^2/2}dy-1\leq
e^{x^2/2}\int_x^\infty ye^{-y^2/2}dy-1=0
\]
so that $s$ is decreasing on $\R_+$) and that $s(0) = \sqrt{\pi/2}$.

The proof of the lemma is complete.
\qed

\bigskip

To illustrate how Stein's method is a powerful approach, we shall use it to prove
the celebrated Berry-Esseen theorem. (Our proof is based on an idea introduced by Ho and Chen in
\cite{HoChen}, see also Bolthausen \cite{bolthausen}.)

\begin{thm}[Berry, Esseen, 1956; see \cite{esseen1956}]\label{CLT+BE}
Let ${\bf X}=(X_1,X_2,\ldots)$ be a sequence of i.i.d. random variables with $E[X_1] = 0$, $E[X^2_1] =1$
and $E[|X_1|^3]<\infty$, and define
\[
V_n = \frac{1}{\sqrt{n}} \sum_{k=1}^n X_k, \quad n\geq 1,
\]
to be the associated sequence of normalized partial sums.
Then, for any $n\geq 1$, one has
\begin{equation}\label{e:berryesseen}
\sup_{x\in\R}\left|P\left(V_n\leq x\right)-\frac{1}{\sqrt{2\pi}}\int_{-\infty}^x e^{-u^2/2}du\right|\leq \frac{33\,E[|X_1|^3]}{\sqrt{n}}.
\end{equation}
\end{thm}

\begin{rem}{\rm
One may actually show that (\ref{e:berryesseen}) holds with the constant $0.4784$ instead
of $33$. This has been proved by Korolev and Shevtsova \cite{korolev} in 2010. (They do not use Stein's method.)
On the other hand, according to Esseen \cite{esseen1956} himself, it is impossible to expect a
universal constant smaller than $0.4097$.
}
\end{rem}

\noindent{\it Proof of (\ref{e:berryesseen}).}
For each $n\geq 2$, let $C_n>0$
be the best possible constant satisfying, for all i.i.d. random variables $X_1,\ldots,X_n$
with $E[|X_1|^3]<\infty$, $E[X_1^2]=1$ and $E[X_1]=0$, that
\begin{equation}\label{kn}
\sup_{x\in\R}\left|P(V_n\leq x)-\frac{1}{\sqrt{2\pi}}\int_{-\infty}^x e^{-u^2/2}du\right|\leq \frac{C_n\,E[|X_1|^3]}{\sqrt{n}}.
\end{equation}
As a first (rough) estimation, we first observe that, since $X_1$ is centered with $E[X_1^2]=1$, one has $E[|X_1|^3]\geq E[X_1^2]^{\frac32}=1$,
so that $C_n\leq \sqrt{n}$.
This is of course not enough to conclude, since we need to show that $C_n\leq 33$.

For any $x\in\R$ and $\e>0$, introduce the function
\[
h_{x,\e}(u)=\left\{
\begin{array}{ll}
1&\quad\mbox{if $u\leq x-\e$}\\
\mbox{linear}&\quad\mbox{if $x-\e<u<x+\e$}\\
0&\quad\mbox{if $u\geq x+\e$}
\end{array}
\right..
\]
It is immediately checked that, for all $n\geq 2$, $\e>0$ and $x\in\R$, we have
\[
E[h_{x-\e,\e}(V_n)]\leq P(V_n\leq x)\leq E[h_{x+\e,\e}(V_n)].
\]
Moreover, for $N\sim\mathcal{N}(0,1)$, $\e>0$ and $x\in\R$, we have, using that the density of $N$ is bounded by
$\frac{1}{\sqrt{2\pi}}$,
\begin{eqnarray*}
E[h_{x+\e,\e}(N)]-\frac{4\e}{\sqrt{2\pi}}\leq E[h_{x-\e,\e}(N)]&\leq& P(N\leq x)\\
&\leq& E[h_{x+\e,\e}(N)]\leq E[h_{x-\e,\e}(N)]+\frac{4\e}{\sqrt{2\pi}}.
\end{eqnarray*}
Therefore, for all $n\geq 2$ and $\e>0$, we have
\[
\sup_{x\in\R}\left|P(V_n\leq x)-\frac{1}{\sqrt{2\pi}}\int_{-\infty}^x e^{-u^2/2}du\right|\leq \sup_{x\in\R}\big|E[h_{x,\e}(V_n)]-E[h_{x,\e}(N)]\big| +\frac{4\e}{\sqrt{2\pi}}.
\]
Assume for the time being that,
for all $\e>0$,
\begin{equation}\label{lmoip}
\sup_{x\in\R}\left|E[h_{x,\e}(V_n)]-E[h_{x,\e}(N)]\right|
\leq
\frac{6\,E[|X_1|^3]}{\sqrt{n}}+\frac{3\,C_{n-1}\,E[|X_1|^3]^2}{\e\,n}.
\end{equation}
We deduce that, for all $\e>0$,
\begin{eqnarray*}
\sup_{x\in\R}\left|P(V_n\leq x)-\frac{1}{\sqrt{2\pi}}\int_{-\infty}^x e^{-u^2/2}du\right|
\leq
\frac{6\,E[|X_1|^3]}{\sqrt{n}}+\frac{3\,C_{n-1}\,E[|X_1|^3]^2}{\e\,n}+
\frac{4\e}{\sqrt{2\pi}}.
\end{eqnarray*}
By choosing $\e=\sqrt{\frac{C_{n-1}}{n}}E[|X_1|^3]$, we get that
\begin{eqnarray*}
\sup_{x\in\R}\left|P(V_n\leq x)-\frac{1}{\sqrt{2\pi}}\int_{-\infty}^x e^{-u^2/2}du\right|
\leq
\frac{E[|X_1|^3]}{\sqrt{n}}\left[6+\left(3
+\frac{4}{\sqrt{2\pi}}\right)\sqrt{C_{n-1}}\right],
\end{eqnarray*}
so that $C_n\leq 6+\left(3
+\frac{4}{\sqrt{2\pi}}\right)\sqrt{C_{n-1}}$.
It follows by induction that $C_n\leq 33$ (recall that $C_n\leq\sqrt{n}$ so that $C_2\leq 33$ in particular), which is the desired conclusion.

We shall now use Stein's Lemma \ref{stein} to prove that (\ref{lmoip}) holds.
Fix $x\in\R$ and $\e>0$, and let $f$ denote the Stein solution associated with $h=h_{x,\e}$, that is,
$f$ satisfies (\ref{fz}).
Observe that $h$ is continuous, and therefore $f$ is $\mathcal{C}^1$.
Recall from Lemma \ref{stein} that
$\|f\|_\infty\leq \sqrt{\frac{\pi}{2}}$ and $\|f'\|_\infty\leq 2$.
Set also $\widetilde{f}(x)=xf(x)$, $x\in\R$. We then have
\begin{equation}\label{cordero-star}
\big|\widetilde{f}(x)-\widetilde{f}(y)\big|=\big|f(x)(x-y)+(f(x)-f(y))y\big|
\leq \left( \sqrt{\frac{\pi}{2}} + 2|y|\right)|x-y|.
\end{equation}
On the other hand, set
\[
V^i_n=V_n-\frac{X_i}{\sqrt{n}},\quad i=1,\ldots,n.
\]
Observe that $V^i_n$ and $X_i$ are independent by construction.
One can thus write
\begin{eqnarray*}
&&E[h(V_n)]-E[h(N)]
=E[f'(V_n)-V_nf(V_n)]\\
&=&\sum_{i=1}^n E\left[
f'(V_n)\frac1n - f(V_n)\,\frac{X_i}{\sqrt{n}}\right]\\
&=&\sum_{i=1}^n E\left[f'(V_n)\frac1n - \big(f(V_n)-f(V^i_n)\big)\frac{X_i}{\sqrt{n}}\right]\,\,\,\mbox{because $E[f(V^i_n)X_i]
=E[f(V^i_n)]E[X_i]=0$}\\
&=&\sum_{i=1}^n E\left[f'(V_n)\frac1n - f'\left(V^i_n+\theta\frac{X_i}{\sqrt{n}}\right)\frac{X_i^2}{n}\right]\,\,\mbox{with
$\theta\sim\mathscr{U}_{[0,1]}$ independent of $X_1,\ldots,X_n$}.
\end{eqnarray*}
We have $f'(x)=\widetilde{f}(x) + h(x) - E[h(N)]$, so that
\begin{equation}\label{aqwzsxedc}
E[h(V_n)]-E[h(N)] = \sum_{i=1}^n \big( a_i(\widetilde{f}) - b_i(\widetilde{f}) + a_i(h) - b_i(h) \big),
\end{equation}
where
\[
a_i(g)=E[g(V_n)-g(V^i_n)]\frac1n\quad\mbox{and}\quad
b_i(g)=E\left[\left(g\left(V^i_n+\theta\frac{X_i}{\sqrt{n}}\right)-g(V^i_n)\right)X_i^2\right]\frac1n.
\]
(Here again, we have used that $V^i_n$ and $X_i$  are independent.)
Hence, to prove that (\ref{lmoip}) holds true, we must bound four terms.
\\
{\it 1st term}. One has, using (\ref{cordero-star}) as well as $E[|X_1|]\leq E[X_1^2]^{\frac12}=1$ and $E[|V^i_n|]\leq E[(V^i_n)^2]^\frac12 \leq 1$,
\begin{eqnarray*}
\big|a_i(\widetilde{f})\big|&\leq &
\frac1{n\sqrt{n}}\left(
E[|X_1|]\sqrt{\frac{\pi}{2}} + 2E[|X_1|]E[|V^i_n|]
\right)\leq \left(\sqrt{\frac{\pi}2} + 2\right)\frac{1}{n\sqrt{n}}.
\end{eqnarray*}
{\it 2nd term}.
Similarly and because $E[\theta]=\frac12$, one has
\begin{eqnarray*}
\big|b_i(\widetilde{f})\big|&\leq &
\frac1{n\sqrt{n}}\left(
E[\theta]E[|X_1|^3]\sqrt{\frac{\pi}{2}} + 2E[\theta]E[|X_1|^3]E[|V^i_n|]
\right)\leq \left(\frac12\sqrt{\frac{\pi}2}+1\right)\frac{E[|X_1|^3]}{n\sqrt{n}}.
\end{eqnarray*}
{\it 3rd term}. By definition of $h$, we have
\[
h(v)-h(u) =(v-u)\int_0^1 h'(u+s(v-u))ds =-\frac{v-u}{2\e}\,E\left[
{\bf 1}_{[x-\e,x+\e]}(u+\widehat{\theta}(v-u))
\right],
\]
with $\widehat{\theta}\sim\mathscr{U}_{[0,1]}$ independent of $\theta$ and $X_1,\ldots,X_n$, so that
\begin{eqnarray*}
\big|a_i(h)\big|
&\leq&\frac{1}{2\e\,n\sqrt{n}}
E\left[
|X_i|{\bf 1}_{[x-\e,x+\e]}\left(V^i_n+\widehat{\theta}\frac{X_i}{\sqrt{n}}\right)
\right]\\
&=&\frac{1}{2\e\,n\sqrt{n}}
E\left[
|X_i|\,
P\left(x-\frac{y}{\sqrt{n}}-\e\leq V^i_n\leq x-\frac{y}{\sqrt{n}}+\e\right)
\bigg|_{y=\widehat{\theta}X_i}
\right]\\
&\leq&
\frac{1}{2\e\,n\sqrt{n}}
\,\,\sup_{y\in\R}\,\,
P\left(x-\frac{y}{\sqrt{n}}-\e\leq V^i_n\leq x-\frac{y}{\sqrt{n}}+\e\right).
\end{eqnarray*}
We are thus left to bound $P(a\leq V_n^i\leq b)$ for all $a,b\in\R$ with $a\leq b$.
For that, set $\widetilde{V}^i_n=\frac{1}{\sqrt{n-1}}\sum_{j\neq i}X_j$, so that $V^i_n=\sqrt{1-\frac1n}\,\widetilde{V}^i_n$.
We then have, using in particular (\ref{kn}) (with $n-1$ instead of $n$) and the fact that the standard Gaussian density is bounded by $\frac{1}{\sqrt{2\pi}}$,
\begin{eqnarray*}
P(a\leq V^i_n\leq b)&=&P\left(\frac{a}{\sqrt{1-\frac1n}}\leq \widetilde{V}^i_n \leq \frac{b}{\sqrt{1-\frac1n}}\right)\\
&=&P\left(\frac{a}{\sqrt{1-\frac1n}}\leq N \leq \frac{b}{\sqrt{1-\frac1n}}\right)\\
&&+P\left(\frac{a}{\sqrt{1-\frac1n}}\leq \widetilde{V}^i_n \leq \frac{b}{\sqrt{1-\frac1n}}\right)
-P\left(\frac{a}{\sqrt{1-\frac1n}}\leq N \leq \frac{b}{\sqrt{1-\frac1n}}\right)\\
&\leq& \frac{b-a}{\sqrt{2\pi}\sqrt{1-\frac1n}}
+\frac{2\,C_{n-1}\,E[|X_1|^3]}{\sqrt{n-1}}.
\end{eqnarray*}
We deduce that
\[
\big|a_i(h)\big|\leq \frac{1}{\sqrt{2\pi}n\sqrt{n-1}} + \frac{C_{n-1}\,E[|X_1|^3]}{n\sqrt{n}\sqrt{n-1}\,\e}.
\]
{\it 4th term}. Similarly, we have
\begin{eqnarray*}
\big|b_i(h)\big|
&=&\frac{1}{2n\sqrt{n}\e}\left|E\left[X_i^3\,\theta\,{\bf 1}_{[x-\e,x+\e]}\left(V^i_n+\widehat{\theta}\,
\theta\frac{X_i}{\sqrt{n}}\right)\right]\right|\\
&\leq&\frac{E[|X_1|^3]}{4n\sqrt{n}\e}\,\,\,\sup_{y\in\R}
P\left(x-\frac{y}{\sqrt{n}}-\e\leq V^i_n\leq x-\frac{y}{\sqrt{n}}+\e\right)\\
&\leq&\frac{E[|X_1|^3]}{2\sqrt{2\pi}n\sqrt{n-1}}+\frac{C_{n-1}\,E[|X_1|^3]^2}{2n\sqrt{n}\sqrt{n-1}\,\e}.
\end{eqnarray*}
Plugging these four estimates into (\ref{aqwzsxedc}) and by using the fact that $n\geq 2$ (and therefore $n-1\geq\frac{n}2$)
and $E[|X_1|^3]\geq 1$,
we deduce the desired conclusion.
\qed

\bigskip

{\bf To go further}.
Stein's method has developed considerably since its first appearance in 1972.
A comprehensive and very nice reference to go further is the book \cite{goldstein} by Chen, Goldstein and Shao, in which several applications of Stein's
method
are carefully developed.

\noindent
\section{Malliavin calculus in a nutshell}\label{malliavin}

The second ingredient for the proof of the Fourth Moment Theorem \ref{NP} is the Malliavin calculus (the first one being Stein's method, as
developed in the previous section).
So, let us introduce the reader to
the basic operators of Malliavin calculus.
For the sake of simplicity and to avoid technicalities that would be useless in this survey, we will only consider the case where
the underlying Gaussian process (fixed once for all throughout the sequel) is a classical Brownian motion $B=(B_t)_{t\geq 0}$
defined on some probability space $(\Omega,\mathcal{F},P)$; we further assume
that the $\sigma$-field $\mathcal{F}$ is generated by $B$.

For a detailed exposition of Malliavin calculus (in a more general context) and for missing proofs,
we refer the reader to the textbooks \cite{book-malliavinstein,nunubook}.

\bigskip

{\bf Dimension one}.
In this first section, we would like
to introduce the basic operators of Malliavin calculus in the simplest situation (where only {\it one} Gaussian random variable
is involved). While easy, it is a sufficiently rich context to encapsulate all the essence of this theory.
We first need to recall some useful properties of Hermite polynomials.
\begin{prop}\label{thm-hermite}
The family $(H_q)_{q\in\mathbb{N}}\subset\R[X]$ of Hermite polynomials has the following properties.
\begin{enumerate}
\item[(a)] $H'_q=qH_{q-1}$ and $H_{q+1}=XH_q-qH_{q-1}$ for all $q\in\mathbb{N}$.
\item[(b)] The family $\left(\frac{1}{\sqrt{q!}}H_q\right)_{q\in\mathbb{N}}$
is an orthonormal basis of $L^2(\R,\frac{1}{\sqrt{2\pi}}e^{-x^2/2}dx)$.
\item[(c)] Let $(U,V)$ be a Gaussian vector with $U,V\sim \mathcal{N}(0,1)$. Then, for all $k,l\in\mathbb{N}$,
\[
E[H_p(U)H_q(V)]=\left\{
\begin{array}{ll}
q! E[UV]^q&\quad\mbox{if $p=q$}\\
0&\quad\mbox{otherwise}.
\end{array}\right.
\]
\end{enumerate}
\end{prop}
{\it Proof}. This is well-known. For a proof,  see, e.g., \cite[Proposition 1.4.2]{book-malliavinstein}.
\qed
\bigskip

Let $\varphi:\R\to\R$ be an element of $L^2(\R,\frac{1}{\sqrt{2\pi}}e^{-x^2/2}dx)$.
Proposition \ref{thm-hermite}(b) implies that $\varphi$ may be expanded (in a unique way) in terms of
 Hermite polynomials as follows:
 \begin{equation}\label{hermite-decompo2}
 \varphi=\sum_{q=0}^\infty a_q H_q.
 \end{equation}
When $\varphi$ is such that $\sum qq!a_q^2<\infty$, let us define
\begin{equation}\label{D1}
D\varphi=\sum_{q=0}^\infty qa_q H_{q-1}.
\end{equation}
Since the Hermite polynomials satisfy $H_q'=qH_{q-1}$ (Proposition \ref{thm-hermite}(a)), observe that
\[
D\varphi=\varphi'
\]
(in the sense of distributions).
Let us now define the {\it Ornstein-Uhlenbeck semigroup} $(P_t)_{t\geq 0}$ by
\begin{equation}\label{pt}
P_t\varphi=\sum_{q=0}^\infty e^{-qt}a_q H_{q}.
\end{equation}
Plainly, $P_0 =Id$, $P_t P_s = P_{t+s}$ ($s,t\geq 0$) and
\begin{equation}\label{DP_t1}
DP_t = e^{-t} P_t D.
\end{equation}
Since $(P_t)_{t\geq 0}$ is a semigroup, it admits a generator $L$ defined as
\[
L=\frac{d}{dt}|_{t=0}P_t.
\]
Of course,
for any $t\geq 0$ one has that
\[
\frac{d}{dt}P_t = \lim_{h\to 0}\frac{P_{t+h}-P_t}{h}
= \lim_{h\to 0}P_t\frac{P_{h}-Id}{h}
= P_t\lim_{h\to 0}\frac{P_{h}-Id}{h}
=P_t\, \frac{d}{dh}\bigg|_{h=0}{P_h}
= P_tL,
\]
and, similarly, $\frac{d}{dt}P_t = LP_t$. Moreover, going back to the definition of $(P_t)_{t\geq 0}$, it
is clear that the domain of $L$ is the set of functions $\varphi\in L^2(\R,\frac{1}{\sqrt{2\pi}}e^{-x^2/2}dx)$ such that $\sum q^2q!a_q^2<\infty$ and that, in this
case,
\[
L\varphi=-\sum_{q=0}^\infty q a_q H_q.
\]
We have the following integration by parts formula, whose proof is straightforward (start with the case $\varphi=H_p$ and $\psi=H_q$, and then use bilinearity and approximation to conclude in the general case) and left to the reader.
\begin{prop}\label{P:Libp}
Let $\varphi$ be in the domain of $L$ and $\psi$ be in the domain of $D$.
Then
\begin{equation}\label{Libp}
\int_\R L\varphi(x)\psi(x)\frac{e^{-x^2/2}}{\sqrt{2\pi}}dx =
-\int_\R D\varphi(x)D\psi(x)\frac{e^{-x^2/2}}{\sqrt{2\pi}}dx.
\end{equation}
\end{prop}

\bigskip

We shall now extend all the previous operators in a situation where, instead of dealing with a random variable of the form $F=\varphi(N)$ (that involves only {\it one} Gaussian
random variable $N$), we deal more generally with a random variable $F$ that is measurable with respect to the Brownian motion $(B_t)_{t\geq 0}$.

\bigskip

{\bf Wiener integral}.
For any adapted\footnote{Any
adapted process $u$
that is either c\`adl\`ag or c\`agl\`ad admits a progressively measurable version.
We will always assume that we are dealing with it.} and square integrable stochastic process $u=(u_t)_{t\geq 0}$, let us denote by $\int_0^\infty u_tdB_t$
its It\^o integral.
Recall from any standard textbook of stochastic analysis that the It\^o integral is a linear functional that takes its
values on $L^2(\Omega)$ and has the following basic features, coming mainly from the independence property of the increments of $B$:
\begin{eqnarray}
E\left[\int_0^\infty u_sdB_s\right]&=&0\label{zeroito}\\
E\left[\int_0^\infty u_sdB_s\times\int_0^\infty v_sdB_s\right]&=&E\left[\int_0^\infty u_sv_sds\right]\label{isomito}.
\end{eqnarray}
In the particular case where $u=f\in L^2(\R_+)$ is {\it deterministic}, we say that $\int_0^\infty f(s)dB_s$ is the {\it Wiener integral of $f$}; it is then easy to show that
\begin{equation}
\int_0^\infty f(s)dB_s\sim \mathcal{N}\left(0,\int_0^\infty f^2(s)ds\right)\label{gaugaugau}.
\end{equation}

\bigskip

{\bf Multiple Wiener-It\^o integrals and Wiener chaoses}.
Let $f\in L^2(\R_+^q)$.
Let us see how one could give a `natural' meaning to the $q$-fold multiple integral
\[
I^B_q(f)=\int_{\R_+^q} f(s_1,\ldots,s_q)dB_{s_1}\ldots dB_{s_q}.
\]
To achieve this goal, we shall use an iterated It\^o integral; the following heuristic `calculations'
are thus natural within this framework:
\begin{eqnarray}
&&\int_{\R_+^q} f(s_1,\ldots,s_q)dB_{s_1}\ldots dB_{s_q}\notag\\
&=&\sum_{\sigma\in\mathfrak{S}_q} \int_{\R_+^q} f(s_1,\ldots,s_q){\bf 1}_{\{s_{\sigma(1)}>\ldots>s_{\sigma(q)}\}}dB_{s_1}\ldots dB_{s_q}\notag\\
&=&\sum_{\sigma\in\mathfrak{S}_q} \int_{0}^\infty dB_{s_{\sigma(1)}}\int_{0}^{s_{\sigma(1)}}
dB_{s_{\sigma(2)}} \ldots\int_{0}^{s_{\sigma(q-1)}}dB_{s_{\sigma(q)}} f(s_1,\ldots,s_q)\notag\\
&=&\sum_{\sigma\in\mathfrak{S}_q} \int_{0}^\infty dB_{t_1}\int_{0}^{t_1}
dB_{t_2} \ldots\int_{0}^{t_{q-1}}dB_{t_q} f(t_{\sigma^{-1}(1)},\ldots,t_{\sigma^{-1}(q)})\notag\\
&=&\sum_{\sigma\in\mathfrak{S}_q} \int_{0}^\infty dB_{t_1}\int_{0}^{t_1}
dB_{t_2} \ldots\int_{0}^{t_{q-1}}dB_{t_q} f(t_{\sigma(1)},\ldots,t_{\sigma(q)}).\label{definat}
\end{eqnarray}

Now, we can use (\ref{definat}) as a natural candidate for being $I^B_q(f)$.
\begin{defi}
Let $q\geq 1$ be an integer.\\
1. When $f\in L^2(\R_+^q)$, we set
\begin{equation}\label{lili3}
I^B_q(f)=\sum_{\sigma\in\mathfrak{S}_q} \int_{0}^\infty dB_{t_1}\int_{0}^{t_1}
dB_{t_2} \ldots\int_{0}^{t_{q-1}}dB_{t_q} f(t_{\sigma(1)},\ldots,t_{\sigma(q)}).
\end{equation}
The random variable $I^B_q(f)$ is called the $q$th multiple Wiener-It\^o integral of $f$.\\
2. The set $\mathcal{H}^B_{q}$ of random variables of the form $I_q^B(f)$, $f\in L^2(\R_+^q)$, is called
the $q$th Wiener chaos of $B$.
We also use the convention $\mathcal{H}^B_{0} = \mathbb{R}$.
\end{defi}

The following properties are readily checked.

\begin{prop}
Let $q\geq 1$ be an integer and let $f\in L^2(\R_+^q)$.\\
1. If $f$ is symmetric (meaning that $f(t_1,\ldots,t_q)=f(t_{\sigma(1)},\ldots,t_{\sigma(q)})$ for any $t\in\R_+^q$ and any permutation
$\sigma\in\mathfrak{S}_q$), then
\begin{equation}\label{lili}
I^B_q(f)=q!\int_{0}^\infty dB_{t_1}\int_{0}^{t_1}dB_{t_2}\ldots\int_{0}^{t_{q-1}}dB_{t_q}\,f(t_1,\ldots,t_q).
\end{equation}
2. We have
\begin{equation}\label{lili2}
I^B_q(f)=I^B_q(\widetilde{f}),
\end{equation}
where $\widetilde{f}$ stands for the symmetrization of $f$ given by
\begin{equation}\label{symq}
\widetilde{f}(t_1,\ldots,t_q)=\frac{1}{q!}\sum_{\sigma\in\mathfrak{S}_q}f(t_{\sigma(1)},\ldots,t_{\sigma(q)}).
\end{equation}
3. For any $p,q\geq 1$, $f\in L^2(\R_+^p)$ and $g\in L^2(\R_+^q)$,
\begin{eqnarray}
E[I^B_q(f)]&=&0\label{isom1}\\
E[I^B_p(f)I^B_q(g)]&=&p!\langle \widetilde{f},\widetilde{g}\rangle_{L^2(\R_+^p)}\quad\mbox{if $p=q$}\label{isom}\\
\label{isom-dif-standard}
E[I^B_p(f)I^B_q(g)]&=&0\quad\mbox{if $p\neq q$}.
\end{eqnarray}
\end{prop}

\medskip

The space $L^2(\Omega)$ can be decomposed into the infinite orthogonal sum of the spaces $\mathcal{H}^B_{q}$.
(It is a statement which is analogous to the content of Proposition \ref{thm-hermite}(b), and it is precisely here that we need to assume that the $\sigma$-field $\mathcal{F}$ is generated by $B$.)
It follows
that any square-integrable random variable
$F\in L^2(\Omega)$ admits the following chaotic expansion:
\begin{equation}
F=E[F]+\sum_{q=1}^{\infty}I_{q}^B(f_{q}),  \label{E}
\end{equation}
where the functions $f_{q}\in L^2(\R_+^q)$ are symmetric and
uniquely determined by $F$. In practice and when $F$ is `smooth' enough, one may rely on Stroock's formula (see \cite{stroock} or \cite[Exercise 1.2.6]{nunubook}) to compute
the  functions $f_q$ explicitely.

\bigskip

The following result contains a very useful property of multiple Wiener-It\^o integrals.
It is in the same spirit as Lemma \ref{hyphyphyp...popotame!}.

\begin{thm}[Nelson, 1973; see \cite{Nelson1973}]\label{hyper-thm}
Let $f\in L^2(\R_+^q)$ with $q\geq 1$. Then, for all $r\geq 2$,
\begin{equation}\label{hypercontractivity}
E\big[|I^B_q(f)|^r\big]\leq [(r-1)^{q}q!]^{r/2}\|f\|^r_{L^2(\R_+^q)}<\infty.
\end{equation}
\end{thm}
\noindent{\it Proof}.
See, e.g., \cite[Corollary 2.8.14]{book-malliavinstein}. (The proof uses the hypercontractivity property of $(P_t)_{t\geq 0}$ defined as (\ref{Ptinfini}).)
\qed

\bigskip

Multiple Wiener-It\^o integrals are linear by construction. Let us see how they behave with respect to multiplication.
To this aim, we need to introduce the concept of {\it contractions}.

\begin{defi}\label{def-con}
When $r\in\{1,\ldots,p\wedge q\}$, $f\in L^2(\R_+^p)$ and
$g\in L^2(\R_+^q)$,
we write $f\otimes_{r} g$ to indicate the $r$th contraction of $f$ and $g$, defined as
being the element of $L^2(\R_+^{p+q-2r})$ given by
\begin{eqnarray}
&&(f\otimes_r g)(t_1,\ldots,t_{p+q-2r})\label{otherform}\\
&=&\int_{\R_+^r}f(t_1,\ldots,t_{p-r},x_1,\ldots,x_r)g(t_{p-r+1},\ldots,t_{p+q-2r},x_1,\ldots,x_r)dx_1\ldots dx_r.\notag
\end{eqnarray}
By convention, we set $f\otimes_0 g= f\otimes g$ as
being the tensor product of $f$ and $g$, that is,
\[
(f\otimes_0 g)(t_1,\ldots,t_{p+q})=
f(t_1,\ldots,t_p)g(t_{p+1},\ldots,t_{p+q}).
\]
\end{defi}
Observe that
\begin{equation}\label{cauchy-standard}
\|f\otimes_r g\|_{L^2(\R_+^{p+q-2r})}\leq \|f\|_{L^2(\R_+^p)}\|g\|_{L^2(\R_+^q)},\quad r=0,\ldots,p\wedge q
\end{equation}
by Cauchy-Schwarz, and that
$f\otimes_p g=\langle f,g\rangle_{L^2(\R_+^p)}$ when $p=q$.
The next result is the fundamental {\it product formula} between two multiple Wiener-It\^o integrals.
\begin{thm}\label{thm-multi}
Let $p,q\geq 1$ and let $f\in L^2(\R_+^p)$ and $g\in L^2(\R_+^q)$ be two symmetric functions. Then
\begin{equation}\label{multiplication}
I^B_p(f)I^B_q(g)=\sum_{r=0}^{p\wedge q} r!\binom{p}{r}\binom{q}{r}I^B_{p+q-2r}(f\widetilde{\otimes}_{r}g),
\end{equation}
where $f\otimes_{r}g$ stands for the contraction (\ref{otherform}).
\end{thm}
\noindent{\it Proof}. 
Theorem \ref{thm-multi} can be established by at least two routes, namely by induction (see, e.g., \cite[page 12]{nunubook}) or by using the concept of diagonal measure in the context of the Engel-Rota-Wallstrom theory (see \cite{PeTa}). Let us proceed to a heuristic proof following this latter strategy. Going back to the very definition of $I^B_p(f)$, we see that the diagonals are avoided. That is, $I^B_p(f)$ can be seen as
\[
I_p^B(f)=\int_{\R_+^p} f(s_1,\ldots,s_p){\bf 1}_{\{s_i\neq s_j,\,i\neq j\}}dB_{s_1}\ldots dB_{s_p}
\]
The same holds for $I^B_q(g)$.
Then we have (just as through Fubini)
\[
I_p^B(f)I_q^B(g)
=\int_{\R_+^{p+q}} f(s_1,\ldots,s_p){\bf 1}_{\{s_i\neq s_j,\,i\neq j\}}
g(t_1,\ldots,t_q){\bf 1}_{\{t_i\neq t_j,\,i\neq j\}}
dB_{s_1}\ldots dB_{s_p}dB_{t_1}\ldots dB_{t_q}.
\]
While there is no diagonals in the first and second blocks,
there are all possible mixed diagonals in the joint writing.
Hence we need to take into account all these diagonals
(whence the combinatorial coefficients in the statement,
which count all possible diagonal sets of size $r$)
and then integrate out (using the rule $(dB_t)^2=dt$). We thus obtain
\[
I_p^B(f)I_q^B(g)=
\sum_{r=0}^{p\wedge q} r!\binom{p}{r}\binom{q}{r}\int_{R_+^{p+q-2r}}(f\otimes_{r}g)
(x_1,\ldots,x_{p+q-2r})dB_{x_1}\ldots dB_{x_{p+q-2r}}
\]
which is exactly the claim (\ref{multiplication}).
\qed

\bigskip

{\bf Malliavin derivatives}.
We shall extend the operator $D$ introduced in (\ref{D1}).
Let $F\in L^2(\Omega)$ and consider its chaotic expansion (\ref{E}).
\begin{defi}
1. When $m\geq 1$ is an integer, we say that $F$ belongs to the Sobolev-Watanabe space $\mathbb{D}^{m,2}$
if
\begin{equation}\label{dm2}
\sum_{q=1}^\infty q^mq!\|f_q\|_{L^2(\R_+^q)}^2<\infty.
\end{equation}
2. When (\ref{dm2}) holds with $m=1$, the Malliavin derivative $DF=(D_tF)_{t\geq 0}$
of $F$ is
the element of $L^2(\Omega\times\R_+)$ given by
\begin{equation}
D_{t}F=\sum_{q=1}^{\infty }qI^B_{q-1}\left( f_{q}(\cdot ,t)\right).  \label{dtf}
\end{equation}
3. More generally, when (\ref{dm2}) holds with an $m$ bigger than or equal to $2$ we define
the $m$th Malliavin derivative $D^mF=(D_{t_1,\ldots,t_m}F)_{t_1,\ldots,t_m\geq 0}$
of $F$ as the element of $L^2(\Omega\times\R_+^m)$ given by
\begin{equation}
D_{t_1,\ldots,t_m}F=\sum_{q=m}^{\infty }q(q-1)\ldots (q-m+1)I^B_{q-m}\left( f_{q}(\cdot ,t_1,\ldots,t_m)\right).  \label{dtfm}
\end{equation}
\end{defi}
The power 2 in the notation  $\mathbb{D}^{m,2}$ is because it is related to the space $L^2(\Omega)$. (There exists a space 
 $\mathbb{D}^{m,p}$ related to $L^p(\Omega)$ but we
 will not use it in this survey.)
On the other hand, it is clear by construction that $D$ is a linear operator.
Also, using (\ref{isom})-(\ref{isom-dif-standard}) it is easy to compute the $L^2$-norm of $DF$ in terms of the kernels $f_q$ appearing in
the chaotic expansion (\ref{E}) of $F$:
\begin{prop}
Let $F\in\mathbb{D}^{1,2}$. We have
\[
E\left[ \Vert DF\Vert _{L^2(\R_+)}^{2}\right]
=\sum_{q=1}^{\infty }qq!\|f_q\|^2_{L^2(\R_+^q)}.
\]
\end{prop}
\noindent{\it Proof}.
By (\ref{dtf}), we can write
\begin{eqnarray*}
E\left[ \Vert DF\Vert _{L^2(\R_+)}^{2}\right]
&=&\int_{\R_+}E\left[\left(\sum_{q=1}^{\infty }qI^B_{q-1}\left( f_{q}(\cdot ,t)\right)\right)^2\right]dt\\
&=&\sum_{p,q=1}^{\infty }pq
\int_{\R_+}E\left[I^B_{p-1}\left( f_{p}(\cdot ,t)\right)I^B_{q-1}\left( f_{q}(\cdot ,t)\right)\right]dt.
\end{eqnarray*}
Using (\ref{isom-dif-standard}), we deduce that
\begin{eqnarray*}
E\left[ \Vert DF\Vert _{L^2(\R_+)}^{2}\right]
=
\sum_{q=1}^{\infty }q^2
\int_{\R_+}E\left[I^B_{q-1}\left( f_{q}(\cdot ,t)\right)^2\right]dt.
\end{eqnarray*}
Finally, using (\ref{isom}), we get that
\begin{eqnarray*}
E\left[ \Vert DF\Vert _{L^2(\R_+)}^{2}\right]
=
\sum_{q=1}^{\infty }q^2(q-1)!
\int_{\R_+}\left\|f_{q}(\cdot ,t)\right\|^2_{L^2(\R_+^{q-1})}dt
=
\sum_{q=1}^{\infty }qq!
\left\|f_{q}\right\|^2_{L^2(\R_+^{q})}.
\end{eqnarray*}
\qed

\bigskip

Let $H_q$ be the $q$th Hermite polynomial (for some $q\geq 1$)
and let $e\in L^2(\R_+)$ have norm 1.
Recall (\ref{linkhermite}) and Proposition \ref{thm-hermite}(a). We deduce that, for any $t\geq 0$,
\begin{eqnarray*}
&&D_t\left(H_q\left(\int_0^\infty e(s)dW_s\right)\right)= D_t(I^B_q(e^{\otimes q}))=qI^B_{q-1}(e^{\otimes q-1})e(t)\\
&=&qH_{q-1}\left(\int_0^\infty e(s)dB_s\right)e(t)
=H'_q\left(\int_0^\infty e(s)dB_s\right)D_t\left(\int_0^\infty e(s)dB_s\right).
\end{eqnarray*}
More generally,
the Malliavin derivative $D$ verifies the {\it chain rule}:

\begin{thm}
Let $\varphi :\mathbb{R}\rightarrow \mathbb{R}$ be both of class $\mathcal{C}^1$ and Lipschitz, and let $F\in{\mathbb{D}}^{1,2}$.
Then, $\varphi(F)\in {\mathbb{D}}^{1,2}$ and
\begin{equation}\label{chainrule}
D_t\varphi (F)=\varphi'(F)D_tF,\quad t\geq 0.
\end{equation}
\end{thm}
\noindent{\it Proof}.
See, e.g., \cite[Proposition 1.2.3]{nunubook}.
\qed

\bigskip

{\bf Ornstein-Uhlenbeck semigroup}.
We now introduce the extension of (\ref{pt}) in our infinite-dimensional setting.

\begin{defi}
The Ornstein-Uhlenbeck semigroup is the family of linear operators $(P_t)_{t\geq 0}$ defined on $L^2(\Omega)$ by
\begin{equation}\label{Ptinfini}
P_tF=\sum_{q=0}^{\infty }e^{-qt}I^B_{q}(f_q),
\end{equation}
where the symmetric kernels $f_q$ are given by (\ref{E}).
\end{defi}

A crucial property of $(P_t)_{t\geq 0}$ is the Mehler formula, that gives an alternative and often useful
representation formula for $P_t$. To be able to state it, we need to introduce a further notation.
Let $(B,B')$ be a two-dimensional Brownian motion
defined on the product probability space $({\bf \Omega},{\bf \mathcal{F}},{\bf P})=(\Omega\times\Omega',\mathscr{F}\otimes
\mathscr{F}',P\times P')$. Let $F\in L^2(\Omega)$. Since $F$ is measurable with respect to the Brownian motion $B$,
we can write $F=\Psi_F(B)$ with $\Psi_F$ a measurable mapping determined
$P\circ B^{-1}$ a.s.. As a consequence, for any $t\geq 0$ the random variable
$\Psi_F(e^{-t} B+\sqrt{1-e^{-2t}}B')$ is well-defined $P\times P'$ a.s.
(note indeed that $e^{-t}B+\sqrt{1-e^{-2t}}B'$ is again a Brownian motion for any $t>0$). We then have the following formula.

\begin{thm}[Mehler's formula] \label{T : Mehler}
For every $F=F(B)\in L^2(\Omega)$ and every $t\geq 0$, we have
\begin{equation}\label{mehler}
P_t(F)=E'\big[\Psi_F(e^{-t}B+\sqrt{1-e^{-2t}}B')\big],
\end{equation}
where $E'$ denotes the expectation with respect to $P'$.
\end{thm}
{\it Proof}.
By using standard arguments, one may show that the linear span of random variables $F$ having the form $F=\exp{\int_0^\infty h(s)dB_s}$
with $h\in L^2(\R_+)$ is dense
in $L^2(\Omega)$. Therefore, it suffices to consider the case where $F$ has this particular form.
On the other hand, we have the following identity, see, e.g., \cite[Proposition 1.4.2$(vi)$]{book-malliavinstein}: for all $c,x\in\R$,
\[
e^{cx-c^2/2} = \sum_{q=0}^\infty \frac{c^q}{q!}H_q(x),
\]
with $H_q$ the $q$th Hermite polynomial.
By setting $c=\|h\|_{L^2(\R_+)}=\|h\|$ and
$x=\int_0^\infty \frac{h(s)}{\|h\|}dB_s$, we deduce that
\[
\exp{\int_0^\infty h(s)dB_s}=e^{\frac12\|h\|^2}\sum_{q=0}^\infty \frac{\|h\|^q}{q!} H_q\left(\int_0^\infty \frac{h(s)}{\|h\|}dB_s\right),
\]
implying in turn, using (\ref{linkhermite}), that
\begin{equation}\label{tgv}
\exp{\int_0^\infty h(s)dB_s}=e^{\frac12\|h\|^2}\sum_{q=0}^\infty \frac{1}{q!} I^B_q\left(h^{\otimes q}\right).
\end{equation}
Thus, for $F=\exp{\int_0^\infty h(s)dB_s}$,
\[
P_t F  = e^{\frac12\|h\|^2}\sum_{q=0}^\infty \frac{e^{-qt}}{q!} I^B_q\left(h^{\otimes q}\right).
\]
On the other hand,
\begin{eqnarray*}
&&E'\big[\Psi_F(e^{-t}B+\sqrt{1-e^{-2t}}B')\big]=
E'\left[\exp\int_0^\infty h(s)(e^{-t}dB_s+\sqrt{1-e^{-2t}}dB'_s)\right]\\
&=&{\rm exp}\left(e^{-t}\int_0^\infty h(s)dB_s\right) \exp\left(\frac{1-e^{-2t}}{2}\,\|h\|^2\right)\\
&=&\exp\left(\frac{1-e^{-2t}}{2}\,\|h\|^2\right)\,e^{\frac{e^{-2t}}2\|h\|^2}\sum_{q=0}^\infty \frac{e^{-qt}}{q!} I^B_q\left(h^{\otimes q}\right)\quad\mbox{by (\ref{tgv})}\\
&=&P_tF.
\end{eqnarray*}
The desired conclusion follows.\qed

\bigskip

{\bf Generator of the Ornstein-Uhlenbeck semigroup}.
Recall the definition (\ref{dm2}) of the Sobolev-Watanabe spaces $\mathbb{D}^{m,2}$, $m\geq 1$, and that
the symmetric kernels $f_q\in L^2(\R_+^q)$ are
uniquely defined through (\ref{E}).

\begin{defi}
1. The generator of the Ornstein-Uhlenbeck semigroup is the linear operator $L$ defined on $\mathbb{D}^{2,2}$ by
\[
LF=-\sum_{q=0}^{\infty }qI^B_{q}(f_q).
\]
2. The pseudo-inverse of $L$ is the linear operator $L^{-1}$ defined on $L^2(\Omega)$ by
\[
L^{-1}F=-\sum_{q=1}^{\infty }\frac{1}q\,I^B_{q}(f_q).
\]
\end{defi}

It is obvious that, for any $F \in L^2(\Omega )$, we have that $L^{-1} F\in  \mathbb{D}^{2,2}$
and
\begin{equation}\label{Lmoins1}
LL^{-1} F = F - E[F].
\end{equation}
Our terminology for $L^{-1}$ is explained by the identity (\ref{Lmoins1}).
Another crucial property of $L$ is contained in the following result, which is the exact
generalization of Proposition \ref{P:Libp}.

\begin{prop}\label{P:Libp2}
Let $F\in\mathbb{D}^{2,2}$ and  $G\in \mathbb{D}^{1,2}$.
Then
\begin{equation}\label{Libp2}
E[LF\times G]=-E[\langle DF,DG\rangle_{L^2(\R_+)}].
\end{equation}
\end{prop}
{\it Proof}. By bilinearity and approximation, it is enough to show (\ref{Libp2})
for $F=I^B_p(f)$ and $G=I^B_q(g)$ with $p,q\geq 1$ and $f\in L^2(\R_+^p)$, $g\in L^2(\R_+^q)$ symmetric.
When $p\neq q$, we have
\[
E[LF\times G]=-pE[I^B_p(f)I^B_q(g)]=0\]
and \[
E[\langle DF,DG\rangle_{L^2(\R_+)}]=pq\int_0^\infty E[I^B_{p-1}(f(\cdot,t))I^B_{q-1}(g(\cdot,t))]dt=0\]
by (\ref{isom-dif-standard}), so the desired conclusion holds true in this case.
When $p=q$,  we have
\[
E[LF\times G]=-pE[I^B_p(f)I^B_p(g)]=-pp!\langle f,g\rangle_{L^2(\R_+^p)}\]
 and
 \begin{eqnarray*}
 E[\langle DF,DG\rangle_{L^2(\R_+)}]&=&p^2\int_0^\infty E[I^B_{p-1}(f(\cdot,t))I^B_{p-1}(g(\cdot,t))]dt\\
 &=&
 p^2(p-1)!\int_0^\infty\langle f(\cdot,t),g(\cdot,t)\rangle_{L^2(\R_+^{p-1})}dt=pp!\langle f,g\rangle_{L^2(\R_+^p)}
 \end{eqnarray*}
 by (\ref{isom}), so the desired conclusion holds true also in this case.
 \qed

 \bigskip

We are now in position to state and prove an integration by parts formula which will play a crucial role in the sequel.
\begin{thm}\label{ivangioformulatheorem}
Let $\varphi :\mathbb{R}\rightarrow \mathbb{R}$ be both of class $\mathcal{C}^1$ and Lipschitz, and let $F\in{\mathbb{D}}^{1,2}$ and $G\in L^2(\Omega)$.
Then
\begin{equation}
{\rm Cov}\big(G,\varphi(F)\big)=E\big[\varphi'(F)\langle DF,-DL^{-1}G\rangle_{L^2(\R_+)}\big].
\label{ivangioformula}
\end{equation}
\end{thm}
\noindent{\it Proof}.
Using the assumptions made on $F$ and $\varphi$, we can write:
\begin{eqnarray*}
{\rm Cov}\big(G,\varphi(F)\big)&=&E\big[L(L^{-1}G)\times \varphi(F)\big]\quad\mbox{(by (\ref{Lmoins1}))}\\
&=&E\left[\langle D\varphi(F),-DL^{-1}G\rangle_{L^2(\R_+)}\right]\quad\mbox{(by (\ref{Libp2}))}\\
&=&E\left[\varphi'(F)\langle D\varphi(F),-DL^{-1}G\rangle_{L^2(\R_+)}\right]\quad\mbox{(by (\ref{chainrule}))},
\end{eqnarray*}
which is the announced formula.
\qed

\bigskip

Theorem \ref{ivangioformulatheorem} admits a useful extension to indicator functions.
Before stating and proving it, we recall the following classical
result from measure theory.
\begin{prop}\label{lusin}
Let $C$ be a Borel set in $\R$, assume that $C\subset[-A,A]$ for some $A>0$, and let $\mu$ be a finite measure on $[-A,A]$.
Then, there exists a sequence $(h_n)$ of continuous functions
with support included in $[-A,A]$ and such that $h_n(x)\in[0,1]$ and ${\bf 1}_C(x)=\lim_{n\to\infty}h_n(x)$ $\mu$-a.e.
\end{prop}
\noindent{\it Proof}. This is an immediate corollary of Lusin's theorem, see e.g. \cite[page 56]{rudin}.
\qed

\bigskip

\begin{cor}\label{cor2}
Let $C$ be a Borel set in $\R$, assume that $C\subset[-A,A]$ for some $A>0$, and let
$F\in \mathbb{D}^{1,2}$ be such that $E[F]=0$.
Then
\[
E\left[F\int_{-\infty}^F {\bf 1}_C(x)dx\right]=E\big[{\bf 1}_{C}(F)\langle DF,-DL^{-1}F\rangle_{L^2(\R_+)}\big].
\]
\end{cor}
\noindent{\it Proof}.
Let $\lambda$ denote the Lebesgue measure and let $P_F$ denote the law of $F$.
By Proposition \ref{lusin} with $\mu=(\lambda+P_F)|_{[-A,A]}$ (that is, $\mu$ is the restriction of $\lambda+P_F$
to $[-A,A]$),
there is a sequence $(h_n)$ of continuous functions with support included in $[-A,A]$ and
such that $h_n(x)\in[0,1]$ and ${\bf 1}_C(x)=\lim_{n\to\infty}h_n(x)$ $\mu$-a.e.
In particular, ${\bf 1}_C(x)=\lim_{n\to\infty} h_n(x)$ $\lambda$-a.e. and $P_F$-a.e.
By Theorem \ref{ivangioformulatheorem}, we have moreover that
\[
E\left[F\int_{-\infty}^F h_n(x)dx\right]=E\big[h_n(F)\langle DF,-DL^{-1}F\rangle_{L^2(\R_+)}\big].
\]
The dominated convergence applies and yields the desired conclusion.
\qed

\bigskip

As a corollary of both Theorem \ref{ivangioformulatheorem} and Corollary \ref{cor2}, we shall prove that the law of any multiple
Wiener-It\^o integral is always absolutely continuous with respect to the Lebesgue measure except, of course, when its kernel is identically zero.
\begin{cor}[Shigekawa; see \cite{Shigekawa}]\label{shigekawa}
Let $q\geq 1$ be an integer and let $f$ be a non zero element of $L^2(\R_+^q)$.
Then the law of $F=I^B_q(f)$ is absolutely continuous with respect to the Lebesgue measure.
\end{cor}
\noindent{\it Proof}. Without loss of generality, we further assume that $f$ is symmetric. The proof is by induction on $q$.
When $q=1$, the desired property is readily checked because
$I^B_1(f) \sim \mathcal{N}(0, \|f\|_{L^2(\R_+)}^2)$, see (\ref{gaugaugau}). Now, let $q\geq 2$ and
assume that the statement of Corollary \ref{shigekawa} holds true for $q-1$, that is, assume that
the law of $I^B_{q-1}(g)$ is absolutely continuous
for any symmetric element $g$ of $L^2(\R_+^{q-1})$ such that
$\|g\|_{L^2(\R_+^{q-1})}>0$.
Let $f$ be a symmetric element of $L^2(\R_+^q)$ with
$\|f\|_{L^2(\R_+^q)}>0$. Let $h\in L^2(\R)$ be such that
$\left\|\int_0^\infty f(\cdot,s)h(s)ds\right\|_{L^2(\R_+^{q-1})} \neq 0$.
(Such an $h$ necessarily exists because, otherwise, we would have that $f(\cdot,s)=0$ for almost all $s\geq 0$ which, by symmetry, would imply
that $f\equiv 0$; this would be in contradiction with our assumption.)
Using the induction assumption, we have that the law of
\[
\langle DF,h\rangle_{L^2(\R_+)}=\int_0^\infty D_sF \,h(s)ds=qI^B_{q-1}\left(\int_0^\infty f(\cdot,s)h(s)ds\right)
\]
is absolutely continuous with respect to the Lebesgue measure.
In particular,
\[
P(\langle DF,h\rangle_{L^2(\R_+)} = 0)=0,
\]
implying in turn,
because $\{\|DF\|_{L^2(\R_+)}=0\}\subset \{\langle DF,h\rangle_{L^2(\R_+)} = 0\}$, that
\begin{equation}\label{bouleauhirsch}
P(\|DF\|_{L^2(\R_+)}>0)=1.
\end{equation}

Now, let $C$ be a Borel set in $\R$. Using Corollary \ref{cor2},
we can write, for every $n\geq 1$,
\begin{eqnarray*}
E\left[{\bf 1}_{C\cap[-n,n]}(F)\frac1q\|DF\|^2_{L^2(\R_+)}\right]&=&
E\left[{\bf 1}_{C\cap[-n,n]}(F)\langle DF,-DL^{-1}F\rangle_{L^2(\R_+)}\right]\\
&=&
E\left[F\int_{-\infty}^F{\bf 1}_{C\cap[-n,n]}(y)dy\right].
\end{eqnarray*}
Assume that the Lebesgue measure of $C$ is zero. The previous equality implies that
\[
E\left[{\bf 1}_{C\cap[-n,n]}(F)\frac1q\|DF\|^2_{L^2(\R_+)}\right]=0,\quad n\geq 1.
\]
But (\ref{bouleauhirsch}) holds as well, so $P(F\in C\cap[-n,n])=0$ for all $n\geq 1$. By monotone convergence, we actually get $P(F\in C)=0$.
This shows that the law of $F$ is absolutely continuous with respect to the Lebesgue measure. The proof of Corollary \ref{shigekawa} is concluded.
\qed

\bigskip

{\bf To go further}.
In the literature, the most quoted reference on Malliavin calculus is the excellent book \cite{nunubook} by Nualart.
It contains many applications of this theory (such as the study of the smoothness of probability laws or
the anticipating
stochastic calculus) and constitutes, as such, an unavoidable reference to go further.

\noindent
\section{Stein meets Malliavin}\label{stein+malliavin}

We are now in a position to prove the Fourth Moment Theorem \ref{NP}. As we will see, to do so we will combine the results of Section \ref{sec:stein} (Stein's method)
with those of Section \ref{malliavin} (Malliavin calculus), thus explaining the title of the current section!
It is a different strategy with respect to the original proof, which is based on the use of the Dambis-Dubins-Schwarz theorem.

We start by introducing the distance we shall use to measure the closeness of the laws of random variables.

\begin{defi}
The total variation distance between the laws of two real-valued random variables $Y$ and $Z$ is defined by
\begin{equation}\label{dtv1}
d_{TV}(Y,Z)=\sup_{C\in\mathscr{B}(\R)} \big| P(Y\in C)-P(Z \in C)\big|,
\end{equation}
where $\mathscr{B}(\R)$ stands for the set of  Borel sets in $\R$.
\end{defi}

When $C\in\mathscr{B}(\R)$, we have that $P(Y\in C\cap[-n,n])\to P(Y\in C)$ and
$P(Z\in C\cap[-n,n])\to P(Z\in C)$ as $n\to\infty$ by the monotone convergence theorem.
So, without loss we may restrict the supremum  in (\ref{dtv1}) to be taken over {\it bounded} Borel sets, that is,
\begin{equation}\label{dtv}
d_{TV}(Y,Z)=
\sup_{\substack{
C\in\mathscr{B}(\R)\\
C\,\mbox{\tiny bounded}}
}
\big| P(Y\in C)-P(Z \in C)\big|.
\end{equation}

We are now ready to derive a bound for the Gaussian approximation of any centered
element $F$ belonging to $\mathbb{D}^{1,2}$.

\begin{thm}[Nourdin, Peccati, 2009; see \cite{NP-PTRF}]\label{nou-pec}
Consider $F\in\mathbb{D}^{1,2}$ with $E[F]=0$.
Then, with $N\sim \mathcal{N}(0,1)$,
\begin{equation}\label{GioIvan}
d_{TV}(F,N)
\leq 2\,E\left[\big|1-\langle DF,-DL^{-1}F\rangle_{L^2(\R_+)}\big|\right].
\end{equation}
\end{thm}
\noindent{\it Proof}.
Let $C$ be a bounded Borel set in $\R$. Let $A>0$ be such that $C\subset [-A,A]$.
Let $\lambda$ denote the Lebesgue measure and let $P_F$ denote the law of $F$.
By Proposition \ref{lusin} with $\mu=(\lambda+P_F)|_{[-A-,A]}$ (the restriction of $\lambda+P_F$
to $[-A,A]$),
there is a sequence $(h_n)$ of continuous functions such that $h_n(x)\in[0,1]$ and ${\bf 1}_C(x)=\lim_{n\to\infty}h_n(x)$
$\mu$-a.e.
By the dominated convergence theorem,
$E[h_n(F)]\to P(F\in C)$
and $E[h_n(N)]\to P(N\in C)$ as $n\to\infty$.
On the other hand, using Lemma \ref{stein} (and denoting by $f_n$ the function associated with $h_n$) as well as (\ref{ivangioformula}) we can write, for each $n$,
\begin{eqnarray*}
\big|E[h_n(F)]-E[h_n(N)]\big|&=&
\big|E[f_{n}'(F)]-E[Ff_{n}(F)]\big|\\
&=&\big|E[f'_{n}(F)(1-\langle DF,-DL^{-1}F\rangle_{L^2(\R_+)}]\\
&\leq&2\,E\big[|1-\langle DF,-DL^{-1}F\rangle_{L^2(\R_+)}|\big].
\end{eqnarray*}
Letting $n$ goes to infinity yields
\[
\big|P(F\in C)-P(N\in C)\big|\leq
2\,E\big[|1-\langle DF,-DL^{-1}F\rangle_{L^2(\R_+)}|\big],
\]
which, together with (\ref{dtv}), implies the desired conclusion.
\qed

\bigskip

{\bf Wiener chaos and the Fourth Moment Theorem}.
In this section, we apply Theorem \ref{nou-pec} to a chaotic random variable $F$, that is, to a
random variable having the specific form of a multiple Wiener-It\^o integral.
We begin with a technical lemma which, among other, shows 
that the fourth moment of $F$ is necessarily greater than $3 E[F^2]^2$.
We recall from Definition \ref{def-con} the meaning of $f\widetilde{\otimes}_{r} f$.
\begin{lemma}\label{L : estimatesonvariance}
Let $q\geq 1$ be an integer and consider a symmetric function $f\in L^2(\R_+^q)$.
Set $F=I^B_q(f)$ and $\sigma^2=E[F^2]=q!\|f\|^2_{L^2(\R_+^q)}$.
The following two identities hold:
\begin{eqnarray}
\notag E\left[\left(\sigma^2-\frac1q\|DF\|^2_{L^2(\R_+)}\right)^2\right]=
\sum_{r=1}^{q-1}\frac{r^2}{q^2}\,r!^2\binom{q}{r}^4 (2q-2r)!\|f\widetilde{\otimes}_{r} f\|^2_{L^2(\R_+^{2q-2r})}\\
\label{lm-control-1d}
\end{eqnarray}
and
\begin{eqnarray}\label{aa3}
&&E[F^4]-3\sigma^4
=\frac{3}q\sum_{r=1}^{q-1}rr!^2\binom{q}{r}^4(2q-2r)!
\|f\widetilde{\otimes}_{r} f\|^2_{L^2(\R_+^{2q-2r})}\\
&=&\sum_{r=1}^{q-1} q!^2\binom{q}{r}^2\left\{
\|f\!\otimes_{r}\! f\|^2_{L^2(\R_+^{2q-2r})}
+\binom{2q-2r}{q-r}\|f\widetilde{\otimes}_{r} f\|^2_{L^2(\R_+^{2q-2r})}
\right\}.\notag\\
\label{aa3bis}
\end{eqnarray}
In particular,
\begin{equation}\label{momentdordre4}
E\left[\left(\sigma^2-\frac1q\|DF\|^2_{L^2(\R_+)}\right)^2\right]\leq \frac{q-1}{3q}\big(E[F^4]-3\sigma^4\big).
\end{equation}
\end{lemma}
\noindent{\it Proof}.
We follow \cite{surveyNP} for (\ref{lm-control-1d})-(\ref{aa3})
and \cite{nunugio} for (\ref{aa3bis}).
For any $t\geq 0$, we have
$D_tF=qI^B_{q-1}\big(f(\cdot,t)\big)$ so that, using (\ref{multiplication}),
\begin{eqnarray}
&&\frac1q\|DF\|^2_{L^2(\R_+)}
=q\int_0^\infty
I^B_{q-1}\big(f(\cdot,t)\big)^2dt\notag\\
&=&q\int_0^\infty\sum_{r=0}^{q-1}r!\binom{q-1}{r}^2 I^B_{2q-2-2r}\big(f(\cdot,t)\widetilde{\otimes}_r f(\cdot,t)\big)dt\notag\\
&=&q\int_0^\infty\sum_{r=0}^{q-1}r!\binom{q-1}{r}^2 I^B_{2q-2-2r}\big(f(\cdot,t)\otimes_r f(\cdot,t)\big)dt\notag\\
&=&q\sum_{r=0}^{q-1}r!\binom{q-1}{r}^2 I^B_{2q-2-2r}\left(\int_0^\infty f(\cdot,t)\otimes_r f(\cdot,t)dt
\right)\notag\\
&=&q\sum_{r=0}^{q-1}r!\binom{q-1}{r}^2 I^B_{2q-2-2r}(f\otimes_{r+1} f)\notag\\
&=&q\sum_{r=1}^{q}(r-1)!\binom{q-1}{r-1}^2 I^B_{2q-2r}(f\otimes_{r} f)\notag\\
&=&q!\|f\|^2_{L^2(\R_+^q)}+q\sum_{r=1}^{q-1}(r-1)!\binom{q-1}{r-1}^2 I^B_{2q-2r}(f\otimes_{r} f).\label{aa1}
\end{eqnarray}
Since $E[F^2]=q!\|f\|^2_{L^2(\R_+^q)}=\sigma^2$, the identity (\ref{lm-control-1d}) follows now from (\ref{aa1}) and the orthogonality properties of multiple Wiener-It\^o integrals.
Recall the hypercontractivity property (\ref{hypercontractivity}) of multiple Wiener-It\^o integrals, and
observe the relations $-L^{-1}F=\frac1q F$ and $D(F^3)=3F^2DF$.
By combining formula (\ref{ivangioformula}) with an approximation argument (the derivative of $\varphi(x)=x^3$ being not bounded),
we infer that
\begin{equation}
E[F^4]=E\big[F\times F^3\big]
=\frac3q E\big[F^2\|DF\|^2_{L^2(\R_+)}\big].\label{aa2}
\end{equation}
Moreover, the multiplication formula (\ref{multiplication}) yields
\begin{equation}
F^2=I^B_q(f)^2=\sum_{s=0}^{q}s!\binom{q}{s}^2 I^B_{2q-2s}(f\widetilde{\otimes}_s f).
\label{dollar2}
\end{equation}
By combining this last identity with (\ref{aa1}) and (\ref{aa2}),
we obtain (\ref{aa3}) and finally (\ref{momentdordre4}).
It remains to prove (\ref{aa3bis}).
Let $\sigma$ be a permutation of $\{1,\ldots,2q\}$
(this fact is written in symbols as $\sigma\in\mathfrak{S}_{2q}$).
If $r\in\{0,\ldots,q\}$ denotes the cardinality of $\{\sigma(1),\ldots,\sigma(q)\}\cap\{1,\ldots,q\}$
then it is readily checked that $r$ is also the cardinality of
$\{\sigma(q+1),\ldots,\sigma(2q)\}\cap\{q+1,\ldots,2q\}$ and that
\begin{eqnarray}
&&\int_{\R_+^{2q}}f(t_1,\ldots,t_q)f(t_{\sigma(1)},\ldots,t_{\sigma(q)})f(t_{q+1},\ldots,t_{2q})\notag\\
&&\hskip4cm\times
f(t_{\sigma(q+1)},\ldots,t_{\sigma(2q)})dt_1\ldots dt_{2q}\notag\\
&=&\int_{\R_+^{2q-2r}}(f\otimes_r f)(x_1,\ldots,x_{2q-2r})^2dx_1\ldots dx_{2q-2r} \notag\\
&=&
\|f\otimes _r f\|^2_{L^2(\R_+^{2q-2r})}.\label{ctr}
\end{eqnarray}
Moreover, for any fixed $r\in\{0,\ldots,q\}$, there are $\binom{q}{r}^2(q!)^2$
permutations $\sigma\in\mathfrak{S}_{2q}$ such that
$\#\{\sigma(1),\ldots,\sigma(q)\}\cap\{1,\ldots,q\}=r$.
(Indeed, such a permutation is completely determined by the choice of: $(a)$ $r$ distinct elements $y_1,\ldots,y_r$ of $\{1,\ldots,q\}$; $(b)$ $q-r$ distinct elements $y_{r+1},\ldots,y_q$ of $\{q+1,\ldots,2q\}$; $(c)$ a bijection between $\{1,\ldots,q\}$ and $\{y_1,\ldots,y_q\}$; $(d)$ a bijection between $\{q+1,\ldots,2q\}$ and $\{1,\ldots,2q\}\setminus \{y_1,\ldots,y_q\}$.)
Now, observe that the symmetrization of $f\otimes f$ is given by
\[
f\widetilde{\otimes} f(t_1,\ldots,t_{2q}) = \frac{1}{(2q)!}
\sum_{\sigma\in\mathfrak{S}_{2q}} f(t_{\sigma(1)},\ldots,t_{\sigma(q)})
f(t_{\sigma(q+1)},\ldots,t_{\sigma(2q)}).
\]
Therefore,
\begin{eqnarray*}
&&\|f\widetilde{\otimes} f\|^2_{L^2(\R_+^{2q})}\\
&=&
\frac{1}{(2q)!^2}
\sum_{\sigma,\sigma'\in\mathfrak{S}_{2q}}
\int_{\R_+^{2q}}
f(t_{\sigma(1)},\ldots,t_{\sigma(q)})f(t_{\sigma(q+1)},\ldots,t_{\sigma(2q)})\\
&&\hskip1cm\times
f(t_{\sigma'(1)},\ldots,t_{\sigma'(q)})f(t_{\sigma'(q+1)},\ldots,t_{\sigma'(2q)})
dt_1\ldots dt_{2q}\\
&=&
\frac{1}{(2q)!}
\sum_{\sigma\in\mathfrak{S}_{2q}}
\int_{\R_+^{2q}}
f(t_{1},\ldots,t_{q})f(t_{q+1},\ldots,t_{2q})\\
&&\hskip1cm\times
f(t_{\sigma(1)},\ldots,t_{\sigma(q)})f(t_{\sigma(q+1)},\ldots,t_{\sigma(2q)})
dt_1\ldots dt_{2q}\\
&=&\frac{1}{(2q)!}\sum_{r=0}^q
\sum_{\substack{\sigma\in\mathfrak{S}_{2q}\\
\{\sigma(1),\ldots,\sigma(q)\}\cap\{1,\ldots,q\}=r
}}
\int_{\R_+^{2q}}
f(t_{1},\ldots,t_{q})f(t_{q+1},\ldots,t_{2q})\\
&&\hskip1cm\times
f(t_{\sigma(1)},\ldots,t_{\sigma(q)})f(t_{\sigma(q+1)},\ldots,t_{\sigma(2q)})
dt_1\ldots dt_{2q}.
\end{eqnarray*}
Using (\ref{ctr}), we deduce that
\begin{equation}
(2q)!\|f\widetilde{\otimes} f\|^2_{L^2(\R_+^{2q})}
=2(q!)^2\|f\|_{L^2(\R_+^{q})}^4+(q!)^2\sum_{r=1}^{q-1}
\binom{q}{r}^2\|f\otimes_r f\|^2_{L^2(\R_+^{2q-2r})}.\label{beautyformula}
\end{equation}
Using the orthogonality and isometry properties of multiple Wiener-It\^o integrals,
the identity (\ref{dollar2})
yields
\begin{eqnarray*}
E[F^4] &=& \sum_{r=0}^{q} (r!)^2\binom{q}{r}^4 (2q-2r)!
\|f\widetilde{\otimes}_r f\|^2_{L^2(\R_+^{2q-2r})}\\
&=&(2q)! \|f\widetilde{\otimes} f\|^2_{L^2(\R_+^{2q})}
+(q!)^2\|f\|^4_{L^2(\R_+^{q})}\\
&&+\sum_{r=1}^{q-1} (r!)^2\binom{q}{r}^4 (2q-2r)!
\|f\widetilde{\otimes}_r f\|^2_{L^2(\R_+^{2q-2r})}.
\end{eqnarray*}
By inserting (\ref{beautyformula}) in the previous identity (and because
$(q!)^2\|f\|^4_{L^2(\R_+^q)}=E[F^2]^2=\sigma^4$), we
get (\ref{aa3bis}).
\qed

\medskip

As a consequence of Lemma \ref{L : estimatesonvariance}, we deduce the following bound on the
total variation distance for the Gaussian approximation of a normalized multiple Wiener-It\^o integral.
This is nothing but Theorem \ref{NP-PTRF} but we restate it for convenience.
\begin{thm}[Nourdin, Peccati, 2009; see \cite{NP-PTRF}]\label{mult-thm}
Let $q\geq 1$ be an integer and consider a symmetric function $f\in L^2(\R_+^q)$.
Set $F=I^B_q(f)$, assume that $E[F^2]=1$,  and let $N\sim\mathcal{N}(0,1)$.
Then
\begin{equation}\label{mult-ineq}
d_{TV}(F,N)\leq  2\sqrt{\frac{q-1}{3q}\big|E[F^4]-3\big|}.
\end{equation}
\end{thm}
\noindent{\it Proof}. Since
$L^{-1}F=-\frac1qF$, we have $\langle DF,-DL^{-1}F\rangle_{L^2(\R_+)}=\frac1q\|DF\|^2_{L^2(\R_+)}$.
So, we only need to apply Theorem \ref{nou-pec} and then formula (\ref{momentdordre4}) to conclude.  \qed

\medskip

The estimate (\ref{mult-ineq}) allows to deduce an easy proof of the following characterization of CLTs on Wiener chaos.
(This is the Fourth Moment Theorem \ref{NP} of Nualart and Peccati!). We note that our proof differs from the original one, which is based on the use of the Dambis-Dubins-Schwarz theorem. 

\begin{cor}[Nualart, Peccati, 2005; see \cite{nunugio}]\label{T : NPNOPTP}
Let $q\geq 1$ be an integer and consider a sequence $(f_n)$ of symmetric functions of $L^2(\R_+^q)$.
Set $F_n=I^B_q(f_n)$ and assume that $E[F_n^2]\to\sigma^2>0$ as $n\to \infty$.
Then,
as $n\to\infty$, the following three assertions are equivalent:
\begin{enumerate}
\item[\rm (i)] $F_n\overset{{\rm Law}}{\to}N\sim\mathcal{N}(0,\sigma^2)$;
\item[\rm (ii)] $E[F_n^4]\to E[N^4] = 3\sigma^4$;
\item[\rm (iii)] $\|f_n\widetilde{\otimes}_r f_n\|_{L^2(\R_+^{2q-2r})}\to 0$ for all $r=1,\ldots,q-1$.
\item[\rm (iv)] $\|f_n\otimes_r f_n\|_{L^2(\R_+^{2q-2r})}\to 0$ for all $r=1,\ldots,q-1$.
\end{enumerate}
\end{cor}
\noindent{\it Proof}.
Without loss of generality, we may and do assume that $\sigma^2=1$ and $E[F_n^2]=1$ for all $n$.
The implication (ii) $\to$ (i) is a direct application of Theorem \ref{mult-thm}.
The implication (i) $\to$ (ii) comes from the Continuous Mapping Theorem together with an approximation argument (observe
that $\sup_{n\geq 1}E[F_n^4]<\infty$ by the hypercontractivity relation (\ref{hypercontractivity})).
The equivalence between (ii) and (iii) is an immediate consequence of (\ref{aa3}).
The implication (iv) $\to$ (iii) is obvious (as  $\|f_n\widetilde{\otimes}_r f_n\|\leq  \|f_n\otimes_r f_n\|$)
whereas the implication (ii) $\to$ (iv) follows from (\ref{aa3bis}).
\qed

\bigskip

{\bf Quadratic variation of the fractional Brownian motion}.
In this section, we aim to illustrate Theorem \ref{nou-pec} in a concrete situation.
More precisely,
we shall use Theorem \ref{nou-pec} in order to derive an explicit bound for the
second-order approximation of the quadratic variation of a fractional Brownian motion on $[0,1]$.

Let $B^H=(B^H_t)_{t\geq 0}$ be a fractional Brownian motion with Hurst
index $H\in(0,1)$. This means that $B^H$ is a centered Gaussian process with covariance function given by
\[
E[B^H_tB^H_s]=\frac12\big(t^{2H}+s^{2H}-|t-s|^{2H}\big),\quad s,t\geq 0.
\]
It is easily checked that $B^H$ is selfsimilar of index $H$ and has stationary increments.

Fractional Brownian motion has been successfully used in order to model a variety of natural phenomena coming from different fields, including hydrology, biology, medicine, economics or traffic networks. 
A natural question is thus the identification of the Hurst parameter from real data. To do so, it is popular and classical to use the quadratic variation (on, say, $[0,1]$), which is
observable and given by
\[
S_n =\sum_{k=0}^{n-1} ( B^H_{(k+1)/n}-B^H_{k/n})^2,\quad n\geq 1.
\]
One may prove (see, e.g., \cite[(2.12)]{nourdinfbm}) that
\begin{equation}\label{imcsq2}
n^{2H-1}S_n\overset{\rm proba}{\to} 1\quad\mbox{as $n\to\infty$}.
\end{equation}
We deduce that the estimator $\widehat{H}_n$, defined as
\[
\widehat{H}_n = \frac12 - \frac{\log S_n}{2\log n},
\]
satisfies $\widehat{H}_n\overset{\rm proba}{\to} 1$ as $n\to\infty$. To study the asymptotic normality, consider
\[
F_n=\frac{n^{2H}}{\sigma_n}\sum_{k=0}^{n -1}\big[(B^H_{(k+1)/n}-B^H_{k/n})^2-n^{-2H}\big]
\,\overset{\rm (law)}{=}\,
\frac{1}{\sigma_n}\sum_{k=0}^{n-1}\big[(B^H_{k+1}-B^H_{k})^2-1\big],
\]
where $\sigma_n>0$ is so that $E[F_n^2]=1$. We then have the following result.
\begin{thm}\label{fBM}
Let $N\sim \mathcal{N}(0,1)$ and assume that $H\leq 3/4$. Then,
$\lim_{n\to\infty}\sigma_n^2/n=2\sum_{r\in\mathbb{Z}}\rho^2(r)$ if $H\in(0,\frac34)$,
with $\rho:\mathbb{Z}\to\R$ given by
\begin{equation}\label{rho}
\rho(r)=\frac12\big(|r+1|^{2H}+|r-1|^{2H}-2|r|^{2H}\big),
\end{equation}
and
$\lim_{n\to\infty}\sigma_n^2/(n\log n)=\frac{9}{16}$ if $H=\frac34$.
Moreover,
there exists a
constant $c_H>0$ (depending only on $H$) such that, for every
$n\geq 1$,
\begin{equation}\label{bmstat}
d_{TV}(F_n,N)\leq c_H\times
\left\{\begin{array}{lll}
\frac1{\sqrt{n}}&\,\,\mbox{if $H\in (0,\frac58)$}\\
\\
\frac{(\log n)^{3/2}}{\sqrt{n}}&\,\,\mbox{if $H=\frac58$}\\
\\
n^{4H-3}&\,\,\mbox{if $H\in (\frac58,\frac34)$}\\
\\
\frac1{\log n}&\,\,\mbox{if $H=\frac34$}\\
\end{array}\right..
\end{equation}
\end{thm}

As an immediate consequence of Theorem \ref{fBM}, provided
$H<3/4$ we obtain that
\begin{equation}\label{imcsq}
\sqrt{n}\big(n^{2H-1}S_n-1\big)\overset{\rm law}{\to}
\mathcal{N}\big(0,2\sum_{r\in \mathbb{Z}}\rho^2(r)\big)
\quad\mbox{as $n\to\infty$},
\end{equation}
implying in turn
\begin{equation}\label{cvest}
\sqrt{n}\log n\big(\widehat{H}_n-H\big)\overset{\rm law}{\to}
\mathcal{N}\big(0,\frac12\sum_{r\in \mathbb{Z}}\rho^2(r)\big)\quad\mbox{as $n\to\infty$}.
\end{equation}
Indeed, we can write
\[
\log x = x-1 - \int_1^xdu\int_1^u \frac{dv}{v^2}\quad\mbox{for all $x>0$},
\]
so that (by considering $x\geq 1$ and $0<x<1$)
\[
\big| \log x + 1-x\big|\leq \frac{(x-1)^2}{2}\left\{
1+\frac{1}{x^2}\right\}\quad\mbox{for all $x>0$}.
\]
As a result,
\[
\sqrt{n}\log n \big(\widehat{H}_n-H\big)=-\frac{\sqrt{n}}{2}
\log(n^{2H-1}S_n)=-\frac{\sqrt{n}}{2}(n^{2H-1}S_n-1)+R_n
\]
with
\[
|R_n|\leq \frac{\big(\sqrt{n}(n^{2H-1}S_n-1)\big)^2}{4\sqrt{n}}\left\{1+\frac{1}{(n^{2H-1}S_n)^2}\right\}.
\]
Using (\ref{imcsq2}) and (\ref{imcsq}), it is clear that
$R_n\overset{\rm proba}{\to} 0$ as $n\to\infty$
and then that (\ref{cvest}) holds true. 

Now we have motivated it, let us go back to the proof of Theorem \ref{fBM}. To perform our calculations, we will mainly follow ideas taken from \cite{BBL}. We first need the following ancillary result.
\begin{lemma}\label{lm-rho}
\begin{itemize}
\item[\rm 1.] For any $r\in\mathbb{Z}$, let $\rho(r)$ be defined by (\ref{rho}).
If $H\neq \frac12$, one has $\rho(r)\sim H(2H-1)|r|^{2H-2}$ as $|r|\to\infty$.
If $H=\frac12$ and $|r|\geq 1$, one has $\rho(r)=0$.
Consequently, $\sum_{r\in\mathbb{Z}} \rho^2(r)<\infty$ if and only if $H<3/4$.
\item[\rm 2.] For all $\alpha>-1$,
we have $\sum_{r=1}^{n-1} r^\alpha \sim \frac{n^{\alpha+1}}{\alpha+1}$ as $n\to\infty$.
\end{itemize}
\end{lemma}
\noindent{\it Proof}.
1. The sequence $\rho$ is symmetric, that is, one has $\rho(n)=\rho(-n)$. When $r\to\infty$,
\[
\rho(r)=H(2H-1)r^{2H-2} + o(r^{2H-2}).
\]
Using the usual criterion for convergence of Riemann sums, we deduce that
$\sum_{r\in\mathbb{Z}} \rho^2(r)<\infty$ if and only if $4H-4<-1$ if and only if $H<\frac34$.\\
2. For $\alpha>-1$, we have:
\[
\frac1n\sum_{r=1}^n \left(\frac{r}n\right)^\alpha \to
\int_{0}^1 x^\alpha dx = \frac{1}{\alpha+1}\quad\mbox{as $n\to\infty$.}
\]
We deduce that $\sum_{r=1}^n r^\alpha \sim \frac{n^{\alpha+1}}{\alpha+1}$ as $n\to\infty$.
\qed

\medskip

We are now in position to prove Theorem \ref{fBM}.

\bigskip

\noindent{\it Proof of Theorem \ref{fBM}}.
Without loss of generality, we will rather use the second expression of $F_n$:
\[
F_n=\frac{1}{\sigma_n}\sum_{k=0}^{n-1}\big[(B^H_{k+1}-B^H_{k})^2-1\big].
\]
Consider the linear span $\mathcal{H}$ of $(B^H_{k})_{k\in\mathbb{N}}$, that is, $\mathcal{H}$ is the closed linear subspace of $L^2(\Omega)$
generated by $(B^H_{k})_{k\in\mathbb{N}}$.
It is a real separable Hilbert space and, consequently,
there exists an isometry $\Phi:\mathcal{H}\to L^2(\R_+)$. For any $k\in\N$, set $e_k = \Phi(B^H_{k+1}-B^H_k)$;
we then have, for all $k,l\in\mathbb{N}$,
\begin{equation}\label{e:kerneldiscrete}
\int_0^\infty e_k(s)e_l(s)ds=E[(B^H_{k+1}-B^H_k)(B^H_{l+1}-B^H_l)]=\rho(k-l)
\end{equation}
with $\rho$ given by (\ref{rho}). Therefore,
\[
\{B^H_{k+1}-B^H_k:\,k\in\N\} \,\overset{\rm law}{=}\,\left\{ \int_0^\infty e_k(s)dB_s:\,k\in\N\right\}=\left\{ I^B_1(e_k):\,k\in\N\right\},
\]
where $B$ is a  Brownian motion and $I^B_p(\cdot)$, $p\geq 1$, stands for the $p$th multiple Wiener-It\^o integral associated to $B$.
As a consequence we can, without loss of generality, replace $F_n$ by
\[
F_n=\frac{1}{\sigma_n}\sum_{k=0}^{n-1}\left[\left(I^B_1(e_k)\right)^2-1\right].
\]
Now, using the multiplication formula (\ref{multiplication}), we deduce that
\[
F_n=I^B_2(f_n),\quad\mbox{with }f_n=\frac{1}{\sigma_n}\sum_{k=0}^{n-1}e_k\otimes e_k.
\]
By using the same arguments as in the proof of Theorem \ref{BM}, we obtain the exact value of $\sigma_n$:
\begin{eqnarray*}
\sigma_n ^2= 2\sum_{k,l=0}^{n-1}\rho^2(k-l)
=2\sum_{|r|<n} (n-|r|)\rho^2(r) .
\end{eqnarray*}
Assume that $H<\frac34$ and write
\[
\frac{\sigma_n^2}{n}=2\sum_{r\in\mathbb{Z}} \rho^2(r)\left(1-\frac{|r|}{n}\right)
{\bf 1}_{\{|r|<n\}}.
\]
Since $\sum_{r\in\mathbb{Z}} \rho^2(r)<\infty$ by Lemma \ref{lm-rho}, we obtain by dominated convergence that, when $H<\frac34$,
\begin{equation}\label{sigman}
\lim_{n\to\infty}
\frac{\sigma_n^2}{n}=2\sum_{r\in\mathbb{Z}} \rho^2(r).
\end{equation}
Assume now that $H=\frac34$. We then have
$
\rho^2(r)\sim \frac{9}{64|r|}
$ as $|r|\to\infty$, implying in turn
\[
n\sum_{|r|<n} \rho^2(r) \sim
\frac{9n}{64}\sum_{0<|r|<n}\frac{1}{|r|}
\sim \frac{9n\log n}{32}
\]
and
\[
\sum_{|r|<n} |r|\rho^2(r) \sim
\frac{9}{64}\sum_{|r|<n} 1
\sim \frac{9n}{32}
\]
as $n\to\infty$.
Hence, when $H=\frac34$,
\begin{equation}\label{sigman2}
\lim_{n\to\infty}\frac{\sigma_n^2}{n\log n}=\frac{9}{16}.
\end{equation}
On the other hand, recall that the convolution of two sequences $\{u(n)\}_{n\in\mathbb{Z}}$
and $\{v(n)\}_{n\in\mathbb{Z}}$ is the sequence $u*v$ defined as $(u*v)(j)=\sum_{n\in\mathbb{Z}}u(n)v(j-n)$,
and observe that $(u*v)(l-i)=\sum_{k\in\mathbb{Z}}u(k-l)v(k-i)$
whenever $u(n)=u(-n)$ and $v(n)=v(-n)$ for all $n\in\mathbb{Z}$.
Set \[
\rho_n(k)=|\rho(k)|{\bf 1}_{\{|k|\leq n-1\}},\quad k\in\mathbb{Z},\,n\geq 1.
\]
We then have (using (\ref{lm-control-1d}) for the first equality, and noticing that $f_n\otimes_1 f_n=f_n\widetilde{\otimes}_1 f_n$),
\begin{eqnarray*}
&&E\left[\left(1-
\frac12\|D[I^B_2(f_n)]
\|^2_{L^2(\R_+)}\right)^2\right]\\
&=&
8\,\|f_n\otimes_1 f_n\|^2_{L^2(\R_+^2)}
=\frac8{\sigma_n^4}\sum_{i,j,k,l=0}^{n-1}\rho(k-l)\rho(i-j)\rho(k-i)\rho(l-j)\\
&\leq&\frac8{\sigma_n^4}\sum_{i,l=0}^{n-1}\sum_{j,k\in\mathbb{Z}}\rho_n(k-l)\rho_n(i-j)\rho_n(k-i)\rho_n(l-j)\\
&=&\frac8{\sigma_n^4}\sum_{i,l=0}^{n-1}(\rho_n*\rho_n)(l-i)^2
\leq \frac{8n}{\sigma_n^4}
\sum_{k\in\mathbb{Z}}(\rho_n*\rho_n)(k)^2
= \frac{8n}{\sigma_n^4} \|\rho_n*\rho_n\|^2_{\ell^2(\mathbb{Z})}.
\end{eqnarray*}
Recall Young's inequality: if $s,p,q\geq 1$ are such that $\frac1p+\frac1q=1+\frac1s$, then
\begin{equation}\label{young}
\|u*v\|_{\ell^s(\mathbb{Z})}\leq \|u\|_{\ell^p(\mathbb{Z})}
\|v\|_{\ell^q(\mathbb{Z})}.
\end{equation}
Let us apply (\ref{young}) with $u=v=\rho_n$, $s=2$ and $p=\frac43$. We get
$\|\rho_n*\rho_n\|^2_{\ell^2(\mathbb{Z})}\leq
 \|\rho_n\|^4_{\ell^\frac43(\mathbb{Z})}$,
so that
\begin{equation}\label{bmst}
E\left[\left(1-
\frac12\|D[I^B_2(f_n)]
\|^2_{L^2(\R_+)}\right)^2\right]
\leq \frac{8n}{\sigma_n^4}\left(\sum_{|k|<n}|\rho(k)|^{\frac43}\right)^3.
\end{equation}
Recall the asymptotic behavior of $\rho(k)$ as $|k|\to\infty$ from Lemma \ref{lm-rho}(1).
Hence
\begin{equation}\label{siza}
\sum_{|k|<n} |\rho(k)|^\frac43=
\left\{\begin{array}{lll}
O(1)&\,\,\mbox{if $H\in (0,\frac58)$}\\
O(\log n)&\,\,\mbox{if $H=\frac58$}\\
O(n^{(8H-5)/3})&\,\,\mbox{if $H\in (\frac58,1)$}.
\end{array}\right.
\end{equation}
Assume first that $H<\frac34$ and recall (\ref{sigman}).
This, together with  (\ref{bmst}) and (\ref{siza}), imply
that
\begin{eqnarray*}
E\left[\left|1-
\frac12\|D[I^B_2(f_n)]
\|^2_{L^2(\R_+)}\right|\right]
&\leq&
\sqrt{E\left[\left(1-
\frac12\|D[I^B_2(f_n)]
\|^2_{L^2(\R_+)}\right)^2\right]}
\\
&\leq& c_H\times
\left\{\begin{array}{lll}
\frac1{\sqrt{n}}&\,\,\mbox{if $H\in (0,\frac58)$}\\
\\
\frac{(\log n)^{3/2}}{\sqrt{n}}&\,\,\mbox{if $H=\frac58$}\\
\\
n^{4H-3}&\,\,\mbox{if $H\in (\frac58,\frac34)$}
\end{array}\right..
\end{eqnarray*}
Therefore, the desired conclusion holds for $H\in(0,\frac34)$ by applying Theorem \ref{nou-pec}.
Assume now that $H=\frac34$ and recall (\ref{sigman2}).
This, together with  (\ref{bmst}) and (\ref{siza}), imply
that
\begin{eqnarray*}
E\left[\left|1-
\frac12\|D[I^B_2(f_n)]
\|^2_{L^2(\R_+)}\right|\right]
&\leq&
\sqrt{E\left[\left(1-
\frac12\|D[I^B_2(f_n)]
\|^2_{L^2(\R_+)}\right)^2\right]}
=O(1/\log n),
\end{eqnarray*}
and leads to the desired conclusion for $H=\frac34$ as well.
\qed

\bigskip

{\bf To go further}.
In \cite{NP-PTRF}, one may find a version of Theorem \ref{nou-pec} where $N$ is replaced by a centered Gamma law (see also \cite{noncentral}).
In \cite{bercu}, one associate to Corollary \ref{T : NPNOPTP} an almost sure central limit theorem.
In \cite{breton-nourdin}, the case where $H$ is bigger than $3/4$ in Theorem \ref{fBM} is analyzed.

\noindent
\section{The smart path method}\label{s:smartpath}

The aim of this section is to prove Theorem \ref{PT}
(that is, the multidimensional counterpart of the Fourth Moment Theorem), and even a more general version of it. Following the approach developed in the previous section for the one-dimensional case, a possible way for achieving this goal would have
consisted in extending Stein's method to the multivariate setting, so to combine them with the tools of Malliavin calculus.
This is indeed the approach developed in \cite{NPRev} and it works well.
In this survey, we will actually proceed differently (we follow \cite{surveyNP}), by using the so-called `smart path method' (which is a popular method in spin glasses
theory, see, e.g., Talagrand \cite{talagrand}).

\medskip

Let us first illustrate this approach in dimension one.
Let $F\in\mathbb{D}^{1,2}$ with $E[F]=0$, let $N\sim\mathcal{N}(0,1)$ and
let $h:\R\to\R$ be a $\mathcal{C}^2$ function satisfying $\|\varphi''\|_\infty<\infty$.
Imagine we want to estimate $E[h(F)]-E[h(N)]$. Without loss of generality, we may assume that $N$ and $F$
are stochastically independent. We further have:
\begin{eqnarray*}
&&E[h(F)]-E[h(N)]
=\int_0^1 \frac{d}{dt}E[h(\sqrt{t}F+\sqrt{1-t}N)] dt\\
&=&\int_0^1 \left(\frac{1}{2\sqrt{t}}E[h'(\sqrt{t}F+\sqrt{1-t}N)F]
-\frac{1}{2\sqrt{1-t}}E[h'(\sqrt{t}F+\sqrt{1-t}N)N]\right)dt.
\end{eqnarray*}
For any $x\in\R$ and $t\in[0,1]$, Theorem \ref{ivangioformulatheorem} implies that
\[
E[h'(\sqrt{t}F+\sqrt{1-t}x)F]=\sqrt{t}\,E[h''(\sqrt{t}F+\sqrt{1-t}x)\langle DF,-DL^{-1}F\rangle_{L^2(\R_+)}],
\]
whereas a classical integration by parts yields
\[
E[h'(\sqrt{t}x+\sqrt{1-t}N)N]=\sqrt{1-t}\,E[h''(\sqrt{t}x+\sqrt{1-t}N)].
\]
We deduce, since $N$ and $F$ are independent, that
\begin{equation}\label{smart1}
E[h(F)]-E[h(N)]=\frac12 \int_0^1 E[h''(\sqrt{t}x+\sqrt{1-t}N)(\langle DF,-DL^{-1}F\rangle_{L^2(\R_+)}-1)]dt,
\end{equation}
implying in turn
\begin{equation}\label{smart1d}
\big|E[h(F)]-E[h(N)]\big|\leq \frac12\|h''\|_\infty E\left[\big|1-\langle DF,-DL^{-1}F\rangle_{L^2(\R_+)}\big|\right],
\end{equation}
compare with (\ref{GioIvan}).

It happens that this approach extends easily to the multivariate setting.
To see why, we will adopt the following short-hand notation: for every $h:\R^d\to\R$ of
class $\mathcal{C}^2$, we set
\[
\|h''\|_\infty = \max_{i,j=1,\ldots,d}\sup_{x\in\R^d}\left|\frac{\partial^2 h}{\partial x_i\partial x_j}(x)\right|.
\]
Theorem \ref{theo:majDist2} below is a first step towards Theorem \ref{PT}, and is nothing but the multivariate counterpart of (\ref{smart1})-(\ref{smart1d}).
\begin{thm}\label{theo:majDist2}
Fix $d\geq 2$ and let $F=(F_1,\ldots,F_d)$ be such that $F_i\in\mathbb{D}^{1,2}$ with $E[F_i]=0$ for any $i$.
Let $C\in\mathcal{M}_d(\R)$ be a symmetric and positive matrix,
and let $N$ be a centered Gaussian vector with covariance $C$.
Then, for any $h:\R^d\to\R$ belonging to $\mathcal{C}^2$ and such that $\|h''\|_\infty<\infty$,
we have
\begin{equation}\label{waser2}
\big|E[h(F)]-E[h(N)]\big|\leq \frac12\|h''\|_\infty
\sum_{i,j=1}^d E\left[\left|C(i,j)-\langle DF_j,-DL^{-1}F_i\rangle_{L^2(\R_+)}\right|\right].
\end{equation}
\end{thm}
{\it Proof}.
Without loss of generality, we assume that $N$ is independent of
the underlying Brownian motion $B$.
Let $h$ be as in the statement of the theorem.
For any $t\in[0,1]$, set
$
\Psi(t)=E\big[h\big(\sqrt{1-t}F+\sqrt{t}N\big)\big],
$
so that
\[
E[h(N)]-E[h(F)]=\Psi(1)
-\Psi(0)= \int_0^1 \Psi'(t)dt.
\]
We easily see that $\Psi$ is differentiable on $(0,1)$ with
\[
\Psi'(t)=
\sum_{i=1}^d E\left[\frac{\partial h}{\partial x_i}\big(\sqrt{1-t}F
+\sqrt{t}N\big)\left(
\frac1{2\sqrt{t}}N_i-\frac{1}{2\sqrt{1-t}}F_i\right)
\right].
\]
By integrating by parts, we can write
\begin{eqnarray*}
&&E\left[
\frac{\partial h}{\partial x_i}\big(\sqrt{1-t}F
+\sqrt{t}N\big)N_i
\right]
=E\left\{
E\left[
\frac{\partial h}{\partial x_i}\big(\sqrt{1-t}x
+\sqrt{t}N\big)N_i
\right]_{|x=F}\right\}\\
&=&\sqrt{t}\sum_{j=1}^d C(i,j)\,E\left\{
E\left[
\frac{\partial^2 h}{\partial x_i\partial x_j}\big(\sqrt{1-t}x
+\sqrt{t}N\big)
\right]_{|x=F}\right\}\\
&=&\sqrt{t}\sum_{j=1}^d C(i,j)\,
E\left[
\frac{\partial^2h}{\partial x_i\partial x_j}\big(\sqrt{1-t}F
+\sqrt{t}N\big)
\right].
\end{eqnarray*}
By using Theorem \ref{ivangioformulatheorem} in order to perform the integration by parts, we can also write
\begin{eqnarray*}
&&E\left[
\frac{\partial h}{\partial x_i}\big(\sqrt{1-t}F
+\sqrt{t}N\big)F_i
\right]
=E\left\{
E\left[
\frac{\partial h}{\partial x_i}\big(\sqrt{1-t}F
+\sqrt{t}x\big)F_i
\right]_{|x=N}\right\}\\
&=&\sqrt{1-t}\sum_{j=1}^d E\left\{
E\left[
\frac{\partial^2 h}{\partial x_i\partial x_j}\big(\sqrt{1-t}F
+\sqrt{t}x\big)
\langle DF_j,-DL^{-1}F_i\rangle_{L^2(\R_+)}
\right]_{|x=N}\right\}\\
&=&\sqrt{1-t}\sum_{j=1}^d
E\left[
\frac{\partial^2 h}{\partial x_i\partial x_j}\big(\sqrt{1-t}F
+\sqrt{t}N\big)
\langle DF_j,-DL^{-1}F_i\rangle_{L^2(\R_+)}
\right].
\end{eqnarray*}
Hence
\[
\Psi'(t)=\frac1{2}
\sum_{i,j=1}^d
E\left[
\frac{\partial^2 h}{\partial x_i\partial x_j}\big(\sqrt{1-t}F
+\sqrt{t}N\big)\left(C(i,j)-
\langle DF_j,-DL^{-1}F_j\rangle_{L^2(\R_+)}\right)
\right],
\]
and the desired conclusion follows.
\qed

\bigskip

We are now in position to prove Theorem \ref{NP}
(using a different approach compared to the original proof; here, we rather follow \cite{NOL}). We will actually even show the following more general version.

\begin{thm}[Peccati, Tudor, 2005; see \cite{PTu04}]\label{PecTud}
Let $d\geq 2$ and $q_d, \ldots, q_1\geq 1$ be some fixed
integers. Consider vectors
\[
F_n=(F_{1,n},\ldots,F_{d,n})=
(I^B_{q_1}(f_{1,n}),\ldots,I^B_{q_d}(f_{d,n})), \quad n\geq 1,
\]
with
$f_{i,n}\in L^2(\R_+^{q_i})$ symmetric.
Let $C\in\mathcal{M}_d(\R)$ be a symmetric and positive matrix,
and let $N$ be a centered Gaussian vector with covariance $C$.
Assume that
\begin{equation}
\label{eq:asympcov}
\lim_{n\to\infty}
E[F_{i,n}F_{j,n}]=C(i,j),\quad 1\leq i,j\leq d.
\end{equation} Then, as $n\to\infty$,  the following two
conditions are equivalent:
\begin{itemize}
\item[(a)] $F_n$
converges in law to $N$;
\item[(b)] for every $1\leq i\leq d$, $F_{i,n}$
converges in law to $\mathcal{N}(0,C(i,i))$.
\end{itemize}
\end{thm}
{\it Proof}.
By symmetry, we assume without loss of generality that
$q_1\leq \ldots\leq q_d$.
The implication $(a)\Rightarrow (b)$ being trivial, we only concentrate on $(b)\Rightarrow (a)$.
So, assume $(b)$ and let us show that $(a)$ holds true.
Thanks to (\ref{waser2}), we are left to show that, for each $i,j=1,\ldots,d$,
\begin{equation}\label{luxo}
\langle DF_{j,n},-DL^{-1}F_{i,n}\rangle_{L^2(\R_+)}
=\frac{1}{q_i}\langle DF_{j,n},DF_{i,n}\rangle_{L^2(\R_+)}
\overset{L^2(\Omega)}\to C(i,j)\quad \mbox{as $n\to\infty$}.
\end{equation}
Observe first that, using the product formula (\ref{multiplication}),
\begin{eqnarray}
\frac{1}{q_i}\langle DF_{j,n},DF_{i,n}\rangle_{L^2(\R_+)}
&=&q_j\int_0^\infty I^B_{q_i-1}(f_{i,n}(\cdot,t))I^B_{q_j-1}(f_{j,n}(\cdot,t))dt\notag\\
&=&q_j\sum_{r=0}^{q_i\wedge q_j-1} r!\binom{q_i-1}{r}\binom{q_j-1}{r} I^B_{q_i+q_j-2-2r}\left(\int_0^\infty
f_{i,n}(\cdot,t)\otimes_r f_{j,n}(\cdot,t)dt\right)\notag\\
&=&q_j\sum_{r=0}^{q_i\wedge q_j-1} r!\binom{q_i-1}{r}\binom{q_j-1}{r} I^B_{q_i+q_j-2-2r}\left(f_{i,n}\otimes_{r+1} f_{j,n}\right)\notag\\
&=&q_j\sum_{r=1}^{q_i\wedge q_j} (r-1)!\binom{q_i-1}{r-1}\binom{q_j-1}{r-1} I^B_{q_i+q_j-2r}(f_{i,n}\otimes_r f_{j,n}).\label{calculinter}
\end{eqnarray}
Now, let us consider all the possible cases for $q_i$ and $q_j$ with
$j\geq i$.

\medskip

{\it First case}: $q_i=q_j=1$. We have
$
\langle DF_{j,n},DF_{i,n}\rangle_{L^2(\R_+)}=\langle f_{i,n},f_{j,n}\rangle_{L^2(\R_+)}=E[F_{i,n}F_{j,n}].
$
But it is our assumption that $E[F_{i,n}F_{j,n}]\to C(i,j)$ so (\ref{luxo}) holds true in this case.

\medskip

{\it Second case}: $q_i=1$ and $q_j\geq 2$. We have
$
\langle DF_{j,n},DF_{i,n}\rangle_{L^2(\R_+)}=\langle f_{i,n},DF_{j,n}\rangle_{L^2(\R_+)}=I^B_{q_j-1}(f_{i,n}\otimes_1 f_{j,n}).
$
We deduce that
\begin{eqnarray*}
E[\langle DF_{j,n},DF_{i,n}\rangle_{L^2(\R_+)}^2]
&=&(q_j-1)!\|f_{i,n}\widetilde{\otimes}_1 f_{j,n}\|^2_{L^2(\R_+^{q_j-1})}
\leq(q_j-1)!\|f_{i,n}\otimes_1 f_{j,n}\|^2_{L^2(\R_+^{q_j-1})}\\
&=&(q_j-1)!\langle f_{i,n}\otimes f_{i,n},f_{j,n}\otimes_{q_j-1}f_{j,n}\rangle_{L^2(\R_+^2)}\\
&\leq&(q_j-1)! \|f_{i,n}\|^2_{L^2(\R_+)} \|f_{j,n}\otimes_{q_j-1}f_{j,n}\|_{L^2(\R_+^2)}\\
&=&(q_j-1)! E[F_{i,n}^2]\|f_{j,n}\otimes_{q_j-1}f_{j,n}\|_{L^2(\R_+^2)}.
\end{eqnarray*}
At this stage, observe the following two facts. First, because $q_i\neq q_j$, we have $C(i,j)=0$ necessarily.
Second, since $E[F_{j,n}^2]\to C(j,j)$ and $F_{j,n}\overset{\rm Law}{\to}\mathcal{N}(0,C(j,j))$, we have by Theorem \ref{T : NPNOPTP}
that $\|f_{j,n}\otimes_{q_j-1}f_{j,n}\|_{L^2(\R_+^2)}\to 0$. Hence,  (\ref{luxo}) holds true in this case as well.

\medskip

{\it Third case}: $q_i=q_j\geq 2$.
By (\ref{calculinter}), we can write
\[
\frac{1}{q_i}\langle DF_{j,n},DF_{i,n}\rangle_{L^2(\R_+)}
=E[F_{i,n}F_{j,n}]+q_i\sum_{r=1}^{q_i-1} (r-1)!\binom{q_i-1}{r-1}^2 I^B_{2q_i-2r}(f_{i,n}\otimes_r f_{j,n}).
\]
We deduce that
\begin{eqnarray*}
&&E\left[\left(\frac{1}{q_i}\langle DF_{j,n},DF_{i,n}\rangle_{L^2(\R_+)}-C(i,j)\right)^2\right]\\
&=&\big(E[F_{i,n}F_{j,n}]-C(i,j)\big)^2
+q_i^2\sum_{r=1}^{q_i-1} (r-1)!^2\binom{q_i-1}{r-1}^4 (2q_i-2r)!\|f_{i,n}\widetilde{\otimes}_r f_{j,n}\|^2_{L^2(\R_+^{2q_i-2r})}.
\end{eqnarray*}
The first term of the right-hand side tends to zero by assumption.
For the second term, we can write, whenever $r\in\{1,\ldots,q_i-1\}$,
\begin{eqnarray*}
\|f_{i,n}\widetilde{\otimes}_r f_{j,n}\|^2_{L^2(\R_+^{2q_i-2r})}&\leq&
\|f_{i,n}\otimes_r f_{j,n}\|^2_{L^2(\R_+^{2q_i-2r})}\\
&=&\langle f_{i,n}\otimes_{q_i-r}f_{i,n},f_{j,n}\otimes_{q_i-r}f_{j,n}\rangle_{L^2(\R_+^{2r})}\\
&\leq&\| f_{i,n}\otimes_{q_i-r}f_{i,n}\|_{L^2(\R_+^{2r})}\|f_{j,n}\otimes_{q_i-r}f_{j,n}\|_{L^2(\R_+^{2r})}.
\end{eqnarray*}
Since $F_{i,n}\overset{\rm Law}{\to}\mathcal{N}(0,C(i,i))$ and  $F_{j,n}\overset{\rm Law}{\to}\mathcal{N}(0,C(j,j))$,
by Theorem \ref{T : NPNOPTP} we have that
$\| f_{i,n}\otimes_{q_i-r}f_{i,n}\|_{L^2(\R_+^{2r})}\|f_{j,n}\otimes_{q_i-r}f_{j,n}\|_{L^2(\R_+^{2r})}\to 0$, thereby
showing that (\ref{luxo}) holds true in our third case.

\medskip

{\it Fourth case}: $q_j>q_i\geq 2$.
By (\ref{calculinter}), we have
\[
\frac{1}{q_i}\langle DF_{j,n},DF_{i,n}\rangle_{L^2(\R_+)}
=q_j\sum_{r=1}^{q_i} (r-1)!\binom{q_i-1}{r-1}\binom{q_j-1}{r-1} I^B_{q_i+q_j-2r}(f_{i,n}\otimes_r f_{j,n}).
\]
We deduce that
\begin{eqnarray*}
&&E\left[\frac{1}{q_i}\langle DF_{j,n},DF_{i,n}\rangle_{L^2(\R_+)}^2\right]\\
&=&q_j^2\sum_{r=1}^{q_i} (r-1)!^2\binom{q_i-1}{r-1}^2\binom{q_j-1}{r-1}^2(q_i+q_j-2r)!\|f_{i,n}\widetilde{\otimes}_r f_{j,n}\|^2_{L^2(\R_+^{q_i+q_j-2r})}.
\end{eqnarray*}
For any $r\in\{1,\ldots,q_i\}$, we have
\begin{eqnarray*}
\|f_{i,n}\widetilde{\otimes}_r f_{j,n}\|^2_{L^2(\R_+^{q_i+q_j-2r})}&\leq&
\|f_{i,n}\otimes_r f_{j,n}\|^2_{L^2(\R_+^{q_i+q_j-2r})}\\
&=&\langle f_{i,n}\otimes_{q_i-r}f_{i,n},f_{j,n}\otimes_{q_j-r}f_{j,n}\rangle_{L^2(\R_+^{2r})}\\
&\leq&\| f_{i,n}\otimes_{q_i-r}f_{i,n}\|_{L^2(\R_+^{2r})}\|f_{j,n}\otimes_{q_j-r}f_{j,n}\|_{L^2(\R_+^{2r})}\\
&\leq&\| f_{i,n}\|_{L^2(\R_+^{q_i})}^2\|f_{j,n}\otimes_{q_j-r}f_{j,n}\|_{L^2(\R_+^{2r})}
\end{eqnarray*}
Since $F_{j,n}\overset{\rm Law}{\to}\mathcal{N}(0,C(j,j))$ and $q_j-r\in\{1,\ldots,q_j-1\}$, by  Theorem \ref{T : NPNOPTP}
we have that
$\|f_{j,n}\otimes_{q_j-r}f_{j,n}\|_{L^2(\R_+^{2r})}\to 0$.
We deduce that (\ref{luxo}) holds true in our fourth case.

\medskip

Summarizing, we have that (\ref{luxo}) is true for any $i$ and $j$, and the proof of the theorem is done.
\qed

\bigskip

When the integers $q_d,\ldots,q_1$ are pairwise disjoint in Theorem \ref{PecTud}, notice that $(\ref{eq:asympcov})$ is automatically
verified with $C(i,j)=0$ for all $i\neq j$, see indeed (\ref{isom-dif-standard}). As such, we recover the version of
Theorem \ref{PecTud} (that is, Theorem \ref{PT}) which was stated and used in Lecture 1 to prove Breuer-Major theorem.

\bigskip

{\bf To go further}.
In \cite{NPRev}, Stein's method is combined with Malliavin calculus in a multivariate setting to provide
bounds for the Wasserstein distance between
the laws of $N\sim\mathcal{N}_d(0,C)$ and $F=(F_1,\ldots,F_d)$ where each $F_i\in\mathbb{D}^{1,2}$ verifies $E[F_i]=0$. Compare with Theorem \ref{theo:majDist2}.

\noindent
\section{Cumulants on the Wiener space}

In this section, following \cite{np-jfa} our aim is to analyze the cumulants
of a given element $F$ of $\mathbb{D}^{1,2}$ and to show how
the formula we shall obtain allows us to give yet another proof of the Fourth Moment
Theorem \ref{NP}.

Let $F$ be a random variable with, say, all the moments (to simplify the exposition). Let $\phi_F$
denote its characteristic function, that is, $\phi_F(t)=E[e^{itF}]$, $t\in\R$.
Then, it is well-known that we may recover the moments of $F$ from
$\phi_F$ through the identity
\[
E[F^j]=(-i)^j\,\frac{d^j}{dt^j}|_{t=0}\,\phi_F(t).
\]
The cumulants of $F$, denoted by $\{\kappa_j(F)\}_{j\geq 1}$, are defined in a similar way, just by replacing
$\phi_F$ by $\log \phi_F$ in the previous expression:
\[
\kappa_j(F)=(-i)^j\,\frac{d^j}{dt^j}|_{t=0}\,\log\phi_F(t).
\]
The first few cumulants are
\begin{eqnarray*}
\kappa_1(F)&=&E[F],\\
\kappa_2(F)&=&E[F^2]-E[F]^2={\rm Var}(F),\\
\kappa_3(F)&=&E[F^3]-3E[F^2]E[F]+2E[F]^3.
\end{eqnarray*}
It is immediate that
\begin{equation}\label{cumeasy}
\kappa_j(F+G)=\kappa_j(F)+\kappa_j(G)\quad\mbox{and}\quad
\kappa_j(\lambda F)=\lambda^j \kappa_j(F)
 \end{equation}
 for all $j\geq 1$,
 when $\lambda\in\R$ and $F$ and $G$ are {\it independent} random variables
(with all the moments). Also, it is easy to express moments in terms of cumulants and vice-versa. Finally, let us observe that the cumulants of $F\sim
\mathcal{N}(0,\sigma^2)$ are all zero, except for the second one which is
$\sigma^2$. This fact, together with the two properties (\ref{cumeasy}), gives a quick proof
of the classical CLT and illustrates that cumulants are often relevant
when wanting to decide whether a given random variable is approximately normally distributed.

The following simple lemma is a useful link between moments and cumulants.
\begin{lemma}\label{decom-cum}
Let $F$ be a random variable (in a given probability space $(\Omega,\mathcal{F},P)$) having all the moments.
Then, for all $m\in\N$,
\[
E[F^{m+1}]=\sum_{s=0}^m \binom{m}{s}\kappa_{s+1}(F)E[F^{m-s}].
\]
\end{lemma}
{\it Proof}.
We can write
\begin{eqnarray*}
E[F^{m+1}]&=&(-i)^{m+1}\frac{d^{m+1}}{dt^{m+1}}|_{t=0}\,\phi_F(t)
=(-i)^{m+1}\frac{d^m}{dt^m}|_{t=0}\,\left(\phi_F(t)\frac{d}{dt}\log\phi_F(t)\right)\\
&=&(-i)^{m+1}\sum_{s=0}^m \binom{m}{s}\left(\frac{d^{s+1}}{dt^{s+1}}|_{t=0}\,\log\phi_F(t)
\right)\left(\frac{d^{m-s}}{dt^{m-s}}|_{t=0}\,\phi_F(t)\right)\quad\mbox{by Leibniz rule}\\
&=&\sum_{s=0}^m \binom{m}{s}\kappa_{s+1}(F)E[F^{m-s}].
\end{eqnarray*}
\qed

From now on, we will deal with a random variable $F$ with all moments that is further measurable with respect to the Brownian motion $(B_t)_{t\geq 0}$. We let the notation of Section \ref{malliavin} prevail
and we consider the chaotic expansion (\ref{E}) of $F$.
We further assume (only to avoid technical issues)
that $F$ belongs to $\mathbb{D}^\infty$, meaning that
$F\in\mathbb{D}^{m,2}$ for all $m\geq 1$ and that
$E[\|D^mF\|^p_{L^2(\R_+^m)}]<\infty$ for all $m\geq 1$ and all $p\geq 2$.
This assumption allows us to introduce recursively the following (well-defined) sequence of random variables related to $F$. Namely, set $\Gamma_0(F)=F$ and
\[
\Gamma_{j+1}(F)=\langle DF,-DL^{-1}\Gamma_j(F)\rangle_{L^2(\R_+)}.
\]
The following result contains a neat expression of the cumulants of $F$
in terms of the family $\{\Gamma_s(F)\}_{s\in\N}$.
\begin{thm}[Nourdin, Peccati, 2010; see \cite{np-jfa}]\label{cumnoupec}
Let $F\in\mathbb{D}^\infty$. Then, for any $s\in\N$,
\[
\kappa_{s+1}(F)=s!E[\Gamma_s(F)].
\]
\end{thm}
{\it Proof}.
The proof is by induction. It consists in computing $\kappa_{s+1}(F)$
 using the induction hypothesis, together with Lemma \ref{decom-cum} and (\ref{ivangioformula}).
 First, the result holds true for $s=0$, as it only says that
 $\kappa_1(F)=E[\Gamma_0(F)]=E[F]$.
 Assume now that $m\geq 1$ is given and that $\kappa_{s+1}(F)=s!E[\Gamma_s(F)]$ for all $s\leq m-1$.
 We can then write
 \begin{eqnarray*}
 \kappa_{m+1}(F)&=&E[F^{m+1}]-\sum_{s=0}^{m-1}\binom{m}{s}\kappa_{s+1}(F)E[F^{m-s}]\quad\mbox{by Lemma \ref{decom-cum}}\\
 &=&E[F^{m+1}]-\sum_{s=0}^{m-1}s!\binom{m}{s}E[\Gamma_s(F)]E[F^{m-s}]\quad\mbox{by the induction hypothesis}.
 \end{eqnarray*}
 On the other hand, by applying (\ref{ivangioformula}) repeatedly, we get
 \begin{eqnarray*}
 E[F^{m+1}]&=&E[F^m]E[\Gamma_0(F)]+{\rm Cov}(F^m,\Gamma_0(F))=E[F^m]E[\Gamma_0(F)]+
 mE[F^{m-1}\Gamma_1(F)]\\
 &=&E[F^m]E[\Gamma_0(F)]+
 mE[F^{m-1}]E[\Gamma_1(F)]+ m{\rm Cov}(F^{m-1},\Gamma_1(F))\\ &=&E[F^m]E[\Gamma_0(F)]+
 mE[F^{m-1}]E[\Gamma_1(F)]+ m(m-1)E[F^{m-2}\Gamma_2(F)] \\
 &=&\ldots\\
 &=&\sum_{s=0}^m s!\binom{m}{s}\,E[F^{m-s}]E[\Gamma_s(F)].
 \end{eqnarray*}
 Thus
 \[
\kappa_{m+1}(F)= E[F^{m+1}]-\sum_{s=0}^{m-1}s!\binom{m}{s}E[\Gamma_s(F)]E[F^{m-s}]
 =m!E[\Gamma_m(F)],
 \]
 and the desired conclusion follows.
\qed

\bigskip

 Let us now focus on the computation of cumulants associated to random variables having the form of a multiple Wiener-It\^{o} integral. The following statement provides a compact representation for the cumulants of such random variables.

\begin{thm}\label{thm-pasmal}
Let $q\geq 2$ and assume that $F=I^B_q(f)$, where $f\in L^2(\R_+^q)$. We have $\kappa_1(F)=0$, $\kappa_2(F)=q!\|f\|^2_{L^2(\R_+^q)}$ and,
for every $s\geq 3$,
\begin{equation}\label{formula-cumulants}
\kappa_s(F)= q!(s-1)!\sum c_q(r_1,\ldots,r_{s-2})
\big\langle (...((f\widetilde{\otimes}_{r_1} f) \widetilde{\otimes}_{r_2} f)\ldots
\widetilde{\otimes}_{r_{s-3}}f)\widetilde{\otimes}_{r_{s-2}}f,f\big\rangle_{L^2(\R_+^q)},
\end{equation}
where the sum $\sum$ runs over all collections of integers $r_1,\ldots,r_{s-2}$ such that:
\begin{enumerate}
\item[(i)] $1\leq r_1,\ldots, r_{s-2}\leq q$;
\item[(ii)] $r_1+\ldots+r_{s-2}=\frac{(s-2)q}{2}$;
\item[(iii)] $r_1<q$, $r_1+r_2 <\frac{3q}{2}$, $\ldots$, $r_1+\ldots +r_{s-3}< \frac{(s-2)q}2$;
\item[(iv)] $r_2\leq 2q-2r_1$, $\ldots$, $r_{s-2}\leq (s-2)q-2r_1-\ldots-2r_{s-3}$;
\end{enumerate}
and where the combinatorial constants $c_q(r_1,\ldots,r_{s-2})$ are recursively defined by the
relations
\[
c_q(r)=q(r-1)!\binom{q-1}{r-1}^2,
\]
and, for $a\geq 2$,
\[
c_q(r_1,\ldots,r_{a})=q(r_{a}-1)!\binom{aq-2r_1-\ldots - 2r_{a-1}-1}{r_{a}-1}
\binom{q-1}{r_{a}-1}c_q(r_1,\ldots,r_{a-1}).
\]
\end{thm}

\begin{remark}
{\rm
\begin{enumerate}
\item
If $sq$ is odd, then $\kappa_s(F)=0$, see indeed condition $(ii)$. This fact is easy to see in any case: use that $\kappa_s(-F)=(-1)^s\kappa_s(F)$ and observe that,
when $q$ is odd, then $F\overset{\rm (law)}{=}-F$ (since
$B\overset{\rm (law)}{=}-B$).
\item
If $q=2$ and $F=I^B_2(f)$ with $f\in L^2(\R_+^2)$, then the only possible integers $r_1,\ldots,r_{s-2}$ verifying $(i)-(iv)$ in the previous statement are $r_1=\ldots=r_{s-2}=1$.
On the other hand, we immediately compute that $c_2(1)=2$, $c_2(1,1)=4$, $c_2(1,1,1)=8$, and so on.
Therefore,
\begin{equation}\label{cumX2}
\kappa_s(I^B_2(f))= 2^{s-1}(s-1)!\big\langle (...(f\otimes_1 f)\ldots f)\otimes_1 f, f\big\rangle_{L^2(\R_+^2)},
\end{equation}
and we recover the classical expression of the cumulants of a double integral.
\item If $q\geq 2$ and $F=I^B_q(f)$, $f\in L^2(\R_+^q)$, then (\ref{formula-cumulants})
for $s=4$ reads
\begin{eqnarray}
\kappa_4(I^B_q(f))&=&6q!\sum_{r=1}^{q-1}c_q(r,q-r)\big\langle (f\widetilde{\otimes}_rf)\widetilde{\otimes}_{q-r}
f,f\big\rangle_{L^2(\R_+^q)}\notag\\
&=&\frac{3}{q}\sum_{r=1}^{q-1}rr!^2\binom{q}{r}^4(2q-2r)!\big\langle (f\widetilde{\otimes}_rf)\otimes_{q-r}
f,f\big\rangle_{L^2(\R_+^q)}\notag\\
&=&\frac{3}{q}\sum_{r=1}^{q-1}rr!^2\binom{q}{r}^4(2q-2r)!\big\langle f\widetilde{\otimes}_rf,f\otimes_r f
\big\rangle_{L^2(\R_+^{2q-2r})}\notag\\
&=&\frac{3}{q}\sum_{r=1}^{q-1}rr!^2\binom{q}{r}^4(2q-2r)!
\| f\widetilde{\otimes}_rf\|^2_{L^2(\R_+^{2q-2r})},
\label{bob}
\end{eqnarray}
and we recover the expression for $\kappa_4(F)$ given in (\ref{aa3})  by a different route.
\end{enumerate}
}
\end{remark}
{\it Proof of Theorem \ref{thm-pasmal}}.
Let us first show the following formula: for any $s\geq 2$, we claim that
\begin{eqnarray}
\Gamma_{s-1}(F)&=&
\sum_{r_1=1}^{q} \ldots\sum_{r_{s-1}=1}^{[(s-1)q-2r_1-\ldots-2r_{s-2}]\wedge q}c_q(r_1,\ldots,r_{s-1})
{\bf 1}_{\{r_1< q\}}
\ldots {\bf 1}_{\{
r_1+\ldots+r_{s-2}< \frac{(s-1)q}2
\}}\notag\\
&&
\hskip5cm\times I^B_{sq-2r_1-\ldots-2r_{s-1}}\big(
(...(f\widetilde{\otimes}_{r_1} f) \widetilde{\otimes}_{r_2} f)\ldots
f)\widetilde{\otimes}_{r_{s-1}}f
\big).\notag\\
\label{for}
\end{eqnarray}
We shall prove (\ref{for}) by induction.
When $s=2$, identity (\ref{for}) simply reads
$\Gamma_{1}(F)=
\sum_{r=1}^{q} c_q(r)
I^B_{2q-2r}(
f\widetilde{\otimes}_{r} f)$
and is nothing but (\ref{aa1}).
Assume now that (\ref{for}) holds for $\Gamma_{s-1}(F)$, and let us prove that it continues to hold for $\Gamma_s(F)$.
We have, using the product formula (\ref{multiplication}) and following the same line of reasoning
as in the proof of (\ref{aa1}),
\begin{eqnarray*}
\Gamma_{s}(F)&=&\langle DF,-DL^{-1}\Gamma_{s-1}F\rangle_{L^2(\R_+)}\\
&=&\sum_{r_1=1}^{q} \ldots\sum_{r_{s-1}=1}^{[(s-1)q-2r_1-\ldots-2r_{s-2}]\wedge q}q
c_q(r_1,\ldots,r_{s-1})
{\bf 1}_{\{r_1< q\}}
\ldots {\bf 1}_{\{
r_1+\ldots+r_{s-2}< \frac{(s-1)q}2
\}}\notag\\
&&
\hskip.5cm\times
{\bf 1}_{\{
r_1+\ldots+r_{s-1}< \frac{sq}2
\}} \big\langle I^B_{q-1}(f),I^B_{sq-2r_1-\ldots-2r_{s-1}-1}\big(
(...(f\widetilde{\otimes}_{r_1} f) \widetilde{\otimes}_{r_2} f)\ldots
f)\widetilde{\otimes}_{r_{s-1}}f
\big)\big\rangle_{L^2(\R_+)}\\
&=&\sum_{r_1=1}^{q} \ldots\sum_{r_{s-1}=1}^{[(s-1)q-2r_1-\ldots-2r_{s-2}]\wedge q}\,\,\,
\sum_{r_{s}=1}^{[sq-2r_1-\ldots-2r_{s-1}]\wedge q}
c_q(r_1,\ldots,r_{s-1})\times q
(r_s-1)!\\
&&\hskip1cm\times\binom{sq-2r_1-\ldots-2r_{s-1}-1}{r_s-1}\binom{q-1}{r_s-1}
{\bf 1}_{\{r_1< q\}}
\ldots {\bf 1}_{\{
r_1+\ldots+r_{s-2}< \frac{(s-1)q}2
\}}\notag\\
&&
\hskip1cm\times
{\bf 1}_{\{
r_1+\ldots+r_{s-1}< \frac{sq}2
\}} I^B_{(s+1)q-2r_1-\ldots-2r_{s}}\big(
(...(f\widetilde{\otimes}_{r_1} f) \widetilde{\otimes}_{r_2} f)\ldots
f)\widetilde{\otimes}_{r_{s}}f
\big),
\end{eqnarray*}
which is the desired formula for $\Gamma_s(F)$.
The proof of (\ref{for}) for all $s\geq 1$ is thus finished.
Now, let us take the expectation on both sides of (\ref{for}). We get
\begin{eqnarray*}
\kappa_{s}(F)&=&(s-1)!E[\Gamma_{s-1}(F)]\\
&=&(s-1)!
\sum_{r_1=1}^{q} \ldots\sum_{r_{s-1}=1}^{[(s-1)q-2r_1-\ldots-2r_{s-2}]\wedge q}c_q(r_1,\ldots,r_{s-1})
{\bf 1}_{\{r_1< q\}}
\ldots {\bf 1}_{\{
r_1+\ldots+r_{s-2}< \frac{(s-1)q}2
\}}
\\
&&
\hskip4cm\times {\bf 1}_{\{
r_1+\ldots+r_{s-1}= \frac{sq}2
\}}\times
(...(f\widetilde{\otimes}_{r_1} f) \widetilde{\otimes}_{r_2} f)\ldots
f)\widetilde{\otimes}_{r_{s-1}}f.
\end{eqnarray*}
Observe that, if $2r_1+\ldots+2r_{s-1}= sq$ and $r_{s-1}\leq
(s-1)q-2r_1-\ldots-2r_{s-2}$ then
$2r_{s-1}=q+(s-1)q-2r_1-\ldots-2r_{s-2}\geq q+r_{s-1}$, so that $r_{s-1}\geq q$.
Therefore,
\begin{eqnarray*}
\kappa_{s}(F)&=&(s-1)!
\sum_{r_1=1}^{q} \ldots\sum_{r_{s-2}=1}^{[(s-2)q-2r_1-\ldots-2r_{s-3}]\wedge q}c_q(r_1,\ldots,r_{s-2},q)
{\bf 1}_{\{r_1< q\}}
\ldots {\bf 1}_{\{
r_1+\ldots+r_{s-3}< \frac{(s-2)q}2
\}}
\\
&&
\hskip4cm\times {\bf 1}_{\{
r_1+\ldots+r_{s-2}= \frac{(s-2)q}2
\}}
\big\langle(...( f\widetilde{\otimes}_{r_1} f) \widetilde{\otimes}_{r_2} f)\ldots
f)\widetilde{\otimes}_{r_{s-2}}f,f\big\rangle_{L^2(\R_+^q)},
\end{eqnarray*}
which is the announced result,
since $c_q(r_1,\ldots,r_{s-2},q)=q!c_q(r_1,\ldots,r_{s-2})$.
\qed

\bigskip

We conclude this section by providing yet another proof (based on our new formula (\ref{formula-cumulants})) of the Fourth Moment Theorem
\ref{NP}. More precisely, let us show by another route that,
if $q\geq 2$ is fixed and if $(F_n)_{n\geq 1}$ is a sequence of the form $F_n=I^B_q(f_n)$
with $f_n\in L^2(\R_+^q)$ such that
$E[F_n^2]=q!\|f_n\|^2_{L^2(\R_+^q)}=1$ for all $n\geq 1$ and $E[F_n^4]\to 3$ as $n\to\infty$,
then $F_n\to \mathcal{N}(0,1)$ in law as $n\to\infty$.

To this end, observe that $\kappa_1(F_n)=0$ and $\kappa_2(F_n)=1$.
To estimate
$\kappa_s(F_n)$, $s\geq 3$, we consider the expression (\ref{formula-cumulants}).
Let $r_1,\ldots,r_{s-2}$ be some integers such that
$(i)$--$(iv)$ in Theorem \ref{thm-pasmal} are satisfied.
Using Cauchy-Schwarz and then successively
\[
\|g\widetilde{\otimes}_r h\|_{L^2(\R_+^{p+q-2r})}\leq
\|g\otimes_r h\|_{L^2(\R_+^{p+q-2r})}\leq
\|g\|_{L^2(\R_+^p)}\|h\|_{L^2(\R_+^q)}
\]
whenever $g\in L^2(\R_+^p)$, $h \in L^2(\R_+^q)$
and $r=1,\ldots,p\wedge q$, we get
that
\begin{eqnarray}
&&\big|\langle (...(f_n\widetilde{\otimes}_{r_1} f_n) \widetilde{\otimes}_{r_2} f_n)\ldots
f_n)\widetilde{\otimes}_{r_{s-2}}f_n,f_n\rangle_{L^2(\R_+^q)}\big|\notag\\
&&\quad\quad\quad\quad\quad\quad\quad\quad\quad \quad\quad\quad\leq
\| (...(f_n\widetilde{\otimes}_{r_1} f_n) \widetilde{\otimes}_{r_2} f_n)\ldots
f_n)\widetilde{\otimes}_{r_{s-2}}f_n\|_{L^2(\R_+^q)}
\|f_n\|_{L^2(\R_+^q)}\notag\\
&& \quad\quad\quad\quad\quad\quad\quad\quad\quad\quad\quad\quad  \leq\|f_n\widetilde{\otimes}_{r_1} f_n\|_{L^2(\R_+^{2q-2r_1})}
\|f_n\|_{L^2(\R_+^q)}^{s-2}\notag\\
&& \quad\quad\quad \quad\quad\quad\quad\quad\quad \quad\quad\quad =(q!)^{1-\frac{s}2}\,\|f_n\widetilde{\otimes}_{r_1} f_n\|_{L^2(\R_+^{2q-2r_1})}.
\label{majoration}
\end{eqnarray}
Since $E[F_n^4]-3=\kappa_4(F_n)\to 0$, we deduce from (\ref{bob}) that
$\|f_n\widetilde{\otimes}_r f_n\|_{L^2(\R_+^{2q-2r})}\to 0$
for all $r=1,\ldots,q-1$. Consequently, by combining (\ref{formula-cumulants})
with (\ref{majoration}), we get that
$\kappa_s(F_n)\to 0$ as $n\to\infty$ for all $s\geq 3$, implying in turn that
$F_n\to \mathcal{N}(0,1)$ in law.
\qed

\bigskip

{\bf To go further}.
The multivariate version of Theorem \ref{cumnoupec} may be found in \cite{NorNour}.

\noindent
\section{A new density formula}

In this section, following \cite{NV} we shall explain how the quantity $\langle DF,-DL^{-1}F\rangle_{L^2(\R_+)}$ is related to the density of
$F\in\mathbb{D}^{1,2}$ (provided it exists). More specifically, when $F\in\mathbb{D}^{1,2}$ is
 such that $E[F]=0$, let us
  introduce the function $g_F:\R\to\R$, defined by means of the following identity:
 \begin{equation}\label{gf}
 g_F(F)=E[\langle DF,-DL^{-1}F\rangle_{L^2(\R_+)}|F].
 \end{equation}
 A key property of the random variable $g_F(F)$ is as follows.
 \begin{prop}
If $F\in\mathbb{D}^{1,2}$ satisfies $E[F]=0$, then $P(g_F(F)\geq 0)=1$.
\end{prop}
 {\it Proof}.
 Let $C$ be a Borel set of $\R$ and set $\phi_n(x)=\int_0^x {\bf 1}_{C\cap[-n,n]}(t)dt$, $n\geq 1$ (with the usual convention
$\int_0^x=-\int_x^0$ for $x<0$). Since $\phi_n$ is increasing and vanishing at zero, we have $x\phi_n(x)\geq 0$
for all $x\in\R$. In particular,
\[
0\leq E[F\phi_n(F)]=E\left[F\int_0^F {\bf 1}_{C\cap[-n,n]}(t)dt\right]
=E\left[F\int_{-\infty}^F {\bf 1}_{C\cap[-n,n]}(t)dt\right].
\]
Therefore, we deduce from Corollary \ref{cor2} that
$E\left[g_F(F) {\bf 1}_{C\cap[-n,n]}(F)\right]\geq 0$.
By dominated convergence, this yields
$E\left[g_F(F) {\bf 1}_{C}(F)\right]\geq 0$,
implying in turn that $P(g_F(F)\geq 0)=1$.\qed

\bigskip

The following theorem gives a new density formula for $F$ in terms of the function $g_F$.
We will then study some of its consequences.

\begin{thm}[Nourdin, Viens, 2009; see \cite{NV}]
\label{firstkey}
Let $F\in\mathbb{D}^{1,2}$ with $E[F]=0$. Then,
the law of $F$ admits a density with respect to Lebesgue measure (say, $\rho:\R\to\R$) if and only if
$P(g_F(F)>0)=1$.
In this case, the
support of $\rho$, denoted by $\mathrm{supp\,}\rho$, is a closed interval of
$\mathbb{R}$ containing zero and we have, for (almost) all $x\in\mathrm{supp\,}%
\rho$:
\begin{equation}
\rho(x)=\frac{E[|F|]}{2g_{F}(x)}\,\,\mathrm{exp}\left(
-\int_{0}^{x}\frac{y\,dy}{g_{F}(y)}\right)  . \label{ohquellebelleformulebis}%
\end{equation}
\end{thm}
\textit{Proof}.
Assume that $P(g_{F}(F)>0)=1$ and
let $C$ be a Borel set. Let $n\geq 1$.
Corollary \ref{cor2} yields
\begin{equation}\label{ffff}
E\left[F\int_{-\infty}^F{\bf 1}_{C\cap[-n,n]}(t)dt\right] =E\big[{\bf 1}_{C\cap[-n,n]}(F)g_F(F)\big].
\end{equation}
Suppose that the Lebesgue measure of $C$ is zero.
Then $\int_{-\infty}^F {\bf 1}_{C\cap[-n,n]}(t)dt = 0$, so that
$E\big[{\bf 1}_{C\cap[-n,n]}(F)g_F(F)\big]=0$
by (\ref{ffff}).
But, since $P(g_F(F)>0)=1$, we get that $P(F\in C\cap[-n,n])=0$ and, by letting $n\to\infty$, that
$P(F\in C)=0$. Therefore, the
Radon-Nikodym criterion is verified, hence implying that the law of $F$ has a density.

Conversely, assume that the law of $F$ has a density, say $\rho$.
Let $\phi:\mathbb{R}\rightarrow\mathbb{R}$ be a
continuous function with compact support, and let $\Phi$ denote any
antiderivative of $\phi$. Note that $\Phi$ is necessarily bounded. We can write:
\begin{eqnarray*}
E\big[\phi(F)g_{F}(F)\big]  &   =&E\big[\Phi(F)F\big]\quad\mbox{by (\ref{ivangioformula})}\\
&  =&\int_{\mathbb{R}}\Phi(x)\,x\,\rho(x)dx \underset{(\ast)}{=}\int_{\mathbb{R}%
}\phi(x)\left(  \int_{x}^{\infty}y\rho(y)dy\right)  dx
 =E\left[  \phi(F)\frac{\int_{F}^{\infty}y\rho(y)dy}{\rho(F)}\right].
\end{eqnarray*}
Equation $(\ast)$ was obtained by integrating by parts, after observing that
\[
\int_{x}^{\infty}y\rho(y)dy\to0\,\mbox{ as $|x|\to\infty$}
\]
(for $x\rightarrow+\infty$, this is because $F\in L^{1}(\Omega)$; for
$x\rightarrow-\infty$, this is because $F$ has mean zero). Therefore, we have
shown that, $P$-a.s.,
\begin{equation}
g_{F}(F) =\frac{\int_{F}^{\infty}y\rho(y)dy}{\rho(F)}.
\label{ohlabelleformule}%
\end{equation}
(Notice that $P(\rho(F)>0)=\int_\R {\bf 1}_{\{\rho(x)>0\}}\rho(x)dx = \int_\R \rho(x)dx=1$, so that identity (\ref{ohlabelleformule}) always
makes sense.)
Since $F\in\mathbb{D}^{1,2}$, one has (see, e.g., \cite[Proposition 2.1.7]{nunubook}) that $\mathrm{supp\,}%
\rho=[\alpha,\beta]$ with $-\infty\leq\alpha<\beta\leq+\infty$.
Since $F$ has zero mean, note that $\alpha<0$ and $\beta>0$ necessarily. For
every $x\in(\alpha,\beta)$, define
\begin{equation}
\varphi\left(  x\right)  =\int_{x}^{\infty}y\rho\left(  y\right)  dy.
\label{phi}%
\end{equation}
The function $\varphi$ is differentiable almost everywhere on $(\alpha,\beta
)$, and its derivative is $-x\rho\left(  x\right)  $. In particular, since
$\varphi(\alpha)=\varphi(\beta)=0$ and $\varphi$ is strictly increasing before
$0$ and strictly decreasing afterwards, we have $\varphi(x)>0$ for all
$x\in(\alpha,\beta)$. Hence, (\ref{ohlabelleformule}) implies that $P(g_{F}(F)>0)=1$.

Finally, let us prove
 (\ref{ohquellebelleformulebis}). Let $\varphi$ still
be defined by (\ref{phi}). On the one hand, we have $\varphi^{\prime
}(x)=-x\rho(x)$ for almost all $x\in\mathrm{supp}\,\rho$. On the other hand,
by (\ref{ohlabelleformule}), we have, for almost all $x\in\mathrm{supp}\,\rho
$,
\begin{equation}
\varphi(x)=\rho(x)g_{F}(x). \label{stop1}%
\end{equation}
By putting these two facts together, we get the following ordinary
differential equation satisfied by $\varphi$:
\[
\frac{\varphi^{\prime}(x)}{\varphi(x)}=-\frac{x}{g_{F}(x)}\quad
\mbox{for almost all $x\in{\rm supp}\,\rho$.}
\]
Integrating this relation over the interval $[0,x]$ yields
\[
\log\varphi(x)=\log\varphi(0)-\int_{0}^{x}\frac{y\,dy}{g_{F}(y)}.
\]
Taking the exponential and using $0=E(F)=E(F_{+})-E(F_{-})$ so that
$E|F|=E(F_{+})+E(F_{-})=2E(F_{+})=2\varphi(0)$, we get
\[
\varphi(x)=\frac{1}{2}\,E[|F|]\,\mathrm{exp}\left(  -\int_{0}^{x}\frac
{y\,dy}{g_{F}(y)}\right)  .
\]
Finally, the desired conclusion comes from (\ref{stop1}). \qed

\bigskip

As a consequence of Theorem \ref{firstkey}, we have the following statement, yielding sufficient conditions in order for the law of $F$ to have a support equal to the real line.

\begin{cor}
\label{key-thm}
Let $F\in\mathbb{D}^{1,2}$ with $E[F]=0$.
Assume that there exists $\sigma_{\min}>0$ such that
\begin{equation}
g_{F}(F)\geq\sigma_{\min}^{2},\quad\mbox{$P$-a.s.} \label{sigmamin}%
\end{equation}
Then the law of $F$, which has a density $\rho$ by Theorem \ref{firstkey}, has $\R$ for support  and
(\ref{ohquellebelleformulebis}) holds almost everywhere in $\mathbb{R}$.
\end{cor}
\textit{Proof}.
It is an immediate consequence of Theorem \ref{firstkey}, except for the fact that
$\mathrm{supp\,}\rho=\R$. For the moment, we just know
that $\mathrm{supp\,}\rho=[\alpha,\beta]$ with $-\infty\leq\alpha
<0<\beta\leq+\infty$. Identity (\ref{ohlabelleformule}) yields
\begin{equation}
\int_{x}^{\infty}y\rho\left(  y\right)  dy\geq\sigma_{\min}^{2}%
\,\rho\left(  x\right)  \quad\mbox{for almost all $x\in(\alpha,\beta)$}.
\label{fgr}%
\end{equation}
Let $\varphi$ be defined by (\ref{phi}), and recall that $\varphi(x)>0$ for
all $x\in(\alpha,\beta)$. When multiplied by $x\in\lbrack{0,\beta)}$, the
inequality (\ref{fgr}) gives $\frac{\varphi^{\prime}\left(  x\right)
}{\varphi\left(  x\right)  }\geq-\frac{x}{\sigma_{\min}^{2}}$.
Integrating this relation over the interval $[0,x]$ yields $\log\varphi\left(
x\right)  -\log\varphi\left(  0\right)  \geq-\frac{x^{2}}{2\,\sigma
_{\min}^{2}}$, i.e., since $\varphi(0)=\frac12E|F|$,
\begin{equation}
\varphi\left(  x\right)  =\int_{x}^{\infty}y\rho\left(  y\right)
dy\geq\frac{1}{2}E|F|e^{-\frac{x^{2}}{2\,\sigma_{\min}^{2}}}. \label{Flb}%
\end{equation}
Similarly, when multiplied by $x\in(\alpha,0]$, inequality (\ref{fgr}) gives
$\frac{\varphi^{\prime}\left(  x\right)  }{\varphi\left(  x\right)  }%
\leq-\frac{x}{\sigma_{\min}^{2}}.$ Integrating this relation over the
interval $[x,0]$ yields $\log\varphi\left(  0\right)  -\log\varphi\left(
x\right)  \leq\frac{x^{2}}{2\,\sigma_{\min}^{2}}$, i.e. (\ref{Flb}) still
holds for $x\in(\alpha,0]$. Now, let us prove that $\beta=+\infty$. If this
were not the case, by definition, we would have $\varphi\left(  \beta\right)
=0$; on the other hand, by letting $x$ tend to $\beta$ in the above
inequality, because $\varphi$ is continuous, we would have $\varphi\left(
\beta\right)  \geq\frac{1}{2}E|F|e^{-\frac{\beta^{2}}{2\sigma_{\min}^{2}%
}}>0$, which contradicts $\beta<+\infty$. The proof of $\alpha=-\infty$ is
similar. In conclusion, we have shown that $\mathrm{supp\,}\rho=\mathbb{R}$.
\qed

\

Using Corollary \ref{key-thm}, we deduce a neat criterion for normality.

\begin{cor}\label{levy}
Let $F\in\mathbb{D}^{1,2}$ with $E[F]=0$ and assume that $F$ is not identically zero. Then $F$ is Gaussian if
and only if $\mathrm{Var}(g_{F}(F) )=0.$
\end{cor}
\textit{Proof}: By (\ref{ivangioformula}) (choose $\varphi(x)=x$, $G=F$ and recall that $E[F]=0$), we have
\begin{equation}\label{varianceformula}
E[\langle DF,-DL^{-1}F\rangle_{\EuFrak H}]=E[F^{2}]=\mathrm{Var}F.
\end{equation}
Therefore, the condition $\mathrm{Var}(g_{F}(F) )=0$ is equivalent to $P(g_{F}(F) =\mathrm{Var}F)=1$.
Let $F\sim\mathcal{N}(0,\sigma^{2})$ with $\sigma>0$. Using
(\ref{ohlabelleformule}), we immediately check that $g_{F}(F) =\sigma^{2}$,
$P$-a.s. Conversely, if $g_{F}(F) =\sigma^{2}>0$ $P$-a.s., then Corollary
\ref{key-thm} implies that the law of $F$ has a density $\rho$, given by
$\rho(x)=\frac{E|F|}{2\sigma^{2}}e^{-\frac{x^{2}}{2\,\sigma^{2}}}$ for almost
all $x\in\mathbb{R}$, from which we immediately deduce that $F\sim
\mathcal{N}(0,\sigma^{2})$.\qed

\

Observe that if $F\sim\mathcal{N}(0,\sigma^{2})$ with $\sigma>0$, then
$E|F|=\sqrt{2/\pi}\,\sigma$, so that the formula
(\ref{ohquellebelleformulebis}) for $\rho$ agrees, of course, with the usual
one in this case.

\

As a `concrete' application of (\ref{ohquellebelleformulebis}), let us consider the following situation.
Let $K:[0,1]^2\to\R$ be a square-integrable kernel such that $K(t,s)=0$ for $s>t$, and consider the centered
Gaussian process $X=(X_t)_{t\in[0,1]}$ defined as
\begin{equation}\label{X}
X_t=\int_0^1 K(t,s)dB_s=\int_0^t K(t,s)dB_s,\quad t\in[0,1].
\end{equation}
Fractional Brownian motion is an instance of such a process, see, e.g., \cite[Section 2.3]{nourdinfbm}.
Consider the maximum
\begin{equation}\label{supX}
Z=\sup_{t\in[0,1]}X_t.
\end{equation}
Assume further that the kernel $K$ satisfies
\begin{equation}\label{K}
\exists c,\alpha>0,\quad \forall s,t\in[0,1]^2,\,s\neq t,\quad 0<\int_0^1 (K(t,u)-K(s,u))^2du \leq c|t-s|^\alpha.
\end{equation}
This latter assumption ensures (see, e.g., \cite{decreununu}) that: $(i)$ $Z\in\mathbb{D}^{1,2}$; $(ii)$ the law of $Z$ has a density with respect to Lebesgue measure; $(iii)$ there exists
a (a.s.) unique random point $\tau\in[0,1]$ where the supremum is attained, that is, such that $Z=X_\tau=\int_0^1 K(\tau,s)dB_s$; and
$(iv)$ $D_t Z = K(\tau,t)$, $t\in[0,1]$.
We claim the following formula.
\begin{prop}
Let $Z$ be given by (\ref{supX}), $X$ be defined as (\ref{X})
and $K\in L^2([0,1]^2)$ be satisfying (\ref{K}).
Then, the law of $Z$ has a density $\rho$ whose support is $\R_+$, given by
\[
\rho(x)=\frac{E|Z-E[Z]|}{2h_Z(x)}{\rm exp}\left(-\int_{E[Z]}^x \frac{(y-E[Z])dy}{h_Z(y)}\right),\quad x\geq 0.
\]
Here,
\[
h_Z(x)=\int_0^\infty e^{-u}{\bf E}\left[R(\tau_0,\tau_u)|Z=x\right]du,
\]
where $R(s,t)=E[X_sX_t]$, $s,t\in[0,1]$, and $\tau_{u}$
is the (almost surely) unique random point where
\[
X^{(u)}_t=\int_0^1 K(t,s)(e^{-u}dB_s+\sqrt{1-e^{-2u}}dB_s^{\prime})
\]
attains its maximum on $[0,1]$, with $(B,B')$ a two-dimensional Brownian motion defined on
the product probability space $({\bf \Omega},{\bf \mathcal{F}},{\bf P})=(\Omega
\times\Omega^{\prime},\mathscr{F}\otimes\mathscr{F}^{\prime},P\times
P^{\prime})$.
\end{prop}
{\it Proof}.
Set $F=Z-E[Z]$. We have
$-D_{t}L^{-1}F=\sum_{q=1}^{\infty}I^B_{q-1}(f_{q}(\cdot,t))$
and $D_{t}F=\sum_{q=1}^{\infty}qI^B_{q-1}(f_{q}(\cdot,t))$. Thus
\begin{align*}
\int_{0}^{\infty}e^{-u}P_{u}(D_{t}F)du  &  =\sum_{q=1}^{\infty}I^B_{q-1}(f_{q}(\cdot,t))
\int_{0}^{\infty}e^{-u}qe^{-(q-1)u} du  =\sum_{q=1}^{\infty}I^B_{q-1}(f_{q}(\cdot,t)).
\end{align*}
Consequently,
\[
-D_tL^{-1}F=\int_{0}^{\infty}e^{-u}P_{u}(D_tF)du,\quad t\in[0,1].
\]
By Mehler's formula (\ref{mehler}), and since $DF=DZ=K(\tau,\cdot)$ with $\tau={\rm argmax}_{t\in[0,1]}\int_0^1 K(t,s)dB_s$,
we deduce that
\[
-D_tL^{-1}F=\int_{0}^{\infty}e^{-u}E^{\prime}[K(\tau_u,t)]du,
\]
implying in turn
\begin{eqnarray*}
g_F(F)&=&E[\langle DF,-DL^{-1}F\rangle_{L^2([0,1])}|F]=\int_0^1 dt\int_0^\infty du\, e^{-u}K(\tau_0,t)E[E'[K(\tau_u,t)|F]]\\
&=&\int_0^\infty e^{-u}E\left[E'\left[\int_0^1 K(\tau_0,t)K(\tau_u,t)dt|F\right]\right]du
=\int_0^\infty e^{-u}E\left[E'\left[R(\tau_0,\tau_u)|F\right]\right]du\\
&=&\int_0^\infty e^{-u}{\bf E}\left[R(\tau_0,\tau_u)|F\right]du.
\end{eqnarray*}
The desired conclusion follows now from Theorem \ref{firstkey}
and the fact that $F=Z-E[Z]$.\qed

\bigskip

{\bf To go further}.
The reference \cite{NV} contains concentration inequalities for centered
random variables $F\in\mathbb{D}^{1,2}$ satisfying $g_F(F)\leq \alpha F +\beta$.
The paper \cite{NQ2} shows how Theorem \ref{firstkey} can lead to
optimal Gaussian density estimates for a class of stochastic equations
with additive noise.

\noindent
\section{Exact rates of convergence}

In this section, we follow \cite{exact}. Let $\{F_n\}_{n\geq 1}$ be a sequence of
random variables in $\mathbb{D}^{1,2}$ such that $E[F_n]=0$,
${\rm Var}(F_n)=1$ and $F_n\overset{\rm law}{\to} N\sim \mathcal{N}(0,1)$ as $n\to\infty$.
Our aim is to develop tools for computing
the exact asymptotic expression of the (suitably normalized) sequence
\begin{equation*}\label{e:rimando}
P(F_n\leq x)-P(N\leq x),\quad n\geq 1,
\end{equation*}
when $x\in\R$ is fixed. This will complement the content of Theorem \ref{nou-pec}.

\bigskip

\noindent
{\bf A technical computation}.
For every fixed $x$, we denote by $f_x:\R\to\R$ the function
\begin{eqnarray}\notag
f_x(u)&=&
e^{u^2/2}\int_{-\infty}^u \big({\bf 1}_{(-\infty,x]}(a)-\Phi(x)\big)e^{-a^2/2}da\\
&=&
\sqrt{2\pi}e^{u^2/2}\times \left\{
\begin{array}{ll}
\Phi(u)(1-\Phi(x))&\quad\mbox{if $u\leq x$}\\
\Phi(x)(1-\Phi(u))&\quad\mbox{if $u\geq x$}
\end{array}
\right.,\label{fx}
\end{eqnarray}
where $\Phi(x)=\frac{1}{\sqrt{2\pi}}\int_{-\infty}^x e^{-a^2/2}da$.
We have the following result.
\begin{prop}
Let $N\sim\mathcal{N}(0,1)$. We have, for every $x\in\R$,
\begin{equation}\label{useful}
E[f'_x(N)N] = \frac{1}{3}(x^2-1)\frac{e^{-x^2/2}}{\sqrt{2\pi}}.
\end{equation}
\end{prop}
{\it Proof}. Integrating by parts (the bracket term is easily shown to vanish), we first obtain
that
\begin{eqnarray*}
E[f'_x(N) N]&=&\int_{-\infty}^{+\infty}f'_x(u)u\frac{e^{-u^2/2}}{\sqrt{2\pi}}du=
\int_{-\infty}^{+\infty}f_x(u)(u^2-1)\frac{e^{-u^2/2}}{\sqrt{2\pi}}du \notag\\
&=&\frac1{\sqrt{2\pi}}\int_{-\infty}^{+\infty}
(u^2-1)\left(\int_{-\infty}^u \big[{\bf 1}_{(-\infty,x]}(a)-\Phi(x)\big]e^{-a^2/2}da\right)du.
\end{eqnarray*}
Integrating by parts once again, this time using the relation $u^2-1=\frac{1}{3}(u^3-3u)'$, we deduce that
\begin{eqnarray*}
&&\int_{-\infty}^{+\infty}
(u^2-1)\left(\int_{-\infty}^u \big[{\bf 1}_{(-\infty,x]}(a)-\Phi(x)\big]e^{-a^2/2}da\right)du
\\
&=&-\frac{1}{3}\int_{-\infty}^{+\infty} (u^3-3u)\big[
{\bf 1}_{(-\infty,x]}(u)-\Phi(x)\big]e^{-u^2/2}du\\
&=&-\frac{1}{3}\left(\int_{-\infty}^{x} (u^3-3u)e^{-u^2/2}du
-\Phi(x)\int_{-\infty}^{+\infty} (u^3-3u)e^{-u^2/2}du\right)\\
&=&\frac{1}{3}\,(x^2-1)e^{-x^2/2},\quad\mbox{since $[(u^2-1)e^{-u^2/2}]'=-(u^3-3u)e^{-u^2/2}$}.
\end{eqnarray*}
\qed

\bigskip

\noindent
{\bf A general result}.
Assume that $\{F_n\}_{n\geq 1}$ is a sequence of (sufficiently regular) centered random variables with unitary variance such that the sequence
\begin{equation}\label{IBPstein3}
\varphi(n):=\sqrt{E[(1-\langle DF_n,-DL^{-1}F_n\rangle_{L^2(\R_+)})^2]},\quad n\geq 1,
\end{equation}
converges to zero as $n\to \infty$.
According to Theorem \ref{nou-pec} one has that, for any $x\in\R$ and as $n\to \infty$,
\begin{equation}\label{e:shura}
P(F_n\leq x) - P(N\leq x)\leq d_{TV}(F_n,N)\leq 2\varphi(n) \to 0,
\end{equation}
where $N\sim \mathcal{N}(0,1)$. The forthcoming result provides a useful criterion in order to
compute an exact asymptotic expression (as $n\to \infty$) for the quantity
\[
\frac{P(F_n\leq x) - P(N\leq x)}{\varphi(n)}, \quad n\geq 1.
\]

\begin{thm}[Nourdin, Peccati, 2010; see \cite{exact}]\label{T : UberMain}
Let $\{F_n\}_{n\geq 1}$ be a sequence of random variables belonging to $\mathbb{D}^{1,2}$,
and such that $E[F_n]=0$, ${\rm Var}[F_n]=1$.
Suppose moreover that the following three conditions hold:
\begin{enumerate}
\item[(i)] we have $0<\varphi(n)<\infty$ for every $n$ and $\varphi(n)\to 0$ as $n\to\infty$;
\item[(ii)] the law of $F_n$ has a density with respect to Lebesgue measure for every $n$;
\item[(iii)] as $n\rightarrow\infty$, the two-dimensional vector
$\left(F_n,\frac{\langle DF_n,-DL^{-1}F_n\rangle_{L^2(\R_+)}-1}{\varphi(n)}\right)$
converges in distribution to a centered
two-dimensional Gaussian vector $(N_1,N_2)$, such that
$E[N_1^2]=E[N_2^2]=1$ and $E[N_1 N_2]=\rho$.
\end{enumerate}
Then, as $n\to\infty$, one has for every $x\in\R$,
\begin{equation}\label{Uber}
\frac{P(F_n\leq x) - P(N\leq x)}{\varphi(n)}
\to
\frac{\rho}{3}(1-x^2)\frac{
e^{-x^2/2}}{\sqrt{2\pi}}.
\end{equation}
\end{thm}
{\it Proof.}  For any integer $n$ and any $\mathcal{C}^1$-function $f$ with a bounded derivative,
we know by Theorem \ref{ivangioformulatheorem} that
\[
E[F_nf(F_n)]=E[f'(F_n)\langle DF_n,-DL^{-1}F_n\rangle_{L^2(\R_+)}].
\]
Fix $x\in\R$ and observe
that the function $f_x$ defined by (\ref{fx}) is not $\mathcal{C}^1$ due to the singularity in $x$. However,
by using a regularization argument given assumption $(ii)$, one can show that the identity
\[
E[F_nf_x(F_n)]=E[f_x'(F_n)\langle DF_n,-DL^{-1}F_n\rangle_{L^2(\R_+)}]
\]
is true for any $n$.
Therefore,
since $P(F_n\leq x)- P(N\leq x)=  E[f_x'(F_n)] - E[F_nf_x(F_n)]$,
we get
\[
\frac{P(F_n\leq x) -
P(N\leq x)}{\varphi(n)}=E\left[f'_x(F_n)\times\frac{1-\langle
DF_n,-DL^{-1}F_n\rangle_{L^2(\R_+)}}{\varphi(n)}\right].
\]
Reasoning as in Lemma \ref{stein}, one may show that $f_x$ is Lipschitz with constant 2. Since
$\varphi(n)^{-1}(1-\langle DF_n,-DL^{-1}F_n\rangle_{L^2(\R_+)})$ has
variance 1 by definition of $\varphi(n)$, we deduce that the sequence
\[
f'_x(F_n)\times\frac{1-\langle DF_n,-DL^{-1}F_n\rangle_{L^2(\R_+)}}{\varphi(n)},
\quad n\geq 1,
\]
is uniformly integrable. Definition (\ref{fx}) shows
that $u\to f'_x(u)$ is continuous at every $u\neq x$. This yields
that, as $n\rightarrow \infty$ and due to assumption $(iii)$,
\[
E\left[f'_x(F_n)\times\frac{1-\langle DF_n,-DL^{-1}F_n\rangle_{L^2(\R_+)}}{\varphi(n)}\right]
\to -E[f'_x(N_1)N_2] =-\rho\, E[f'_x(N_1)N_1].
\]
Consequently, relation (\ref{Uber}) now follows from formula (\ref{useful}).
\qed

\bigskip

\noindent
{\bf The double integrals case and a concrete application}.
When applying Theorem \ref{T : UberMain} in concrete situations,
the main issue is often to check that condition $(ii)$ therein holds true.
In the particular case of sequences belonging to the second Wiener chaos,
we can go further in the analysis, leading to the following result.

\begin{prop}\label{roc}
Let
$N\sim\mathcal{N}(0,1)$ and
let $F_n=I^B_2(f_n)$ be such that $f_n\in L^2(\R_+^2)$ is symmetric for all $n\geq 1$.
Write $\kappa_p(F_n)$, $p\geq 1$, to indicate the sequence of the cumulants of $F_n$. Assume that
$\kappa_2(F_n)=E[F_n^2]=1$ for all $n\geq 1$ and that
$\kappa_4(F_n)=E[F_n^4]-3\to 0$ as $n\to\infty$. If we have in addition that
\begin{equation}\label{EZZ}
\frac{\kappa_3(F_n)}{\sqrt{\kappa_4(F_n)}}\to
\alpha  \quad \mbox{and}\quad
\frac{\kappa_8(F_n)}{\big(\kappa_4(F_n)\big)^2}\to
0,
\end{equation}
then, for all $x\in\R$,
\begin{equation}\label{exact-rate}
\frac{P(F_n\leq x)-P(N\leq x)}{
\sqrt{\kappa_4(F_n)}
}\to \frac{\alpha}{6\sqrt{2\pi}}
\left(1-x^2\right)e^{-\frac{x^2}2}\quad\mbox{as $n\to\infty$}.
\end{equation}
\end{prop}

\begin{rem}\label{roc-rem}
{\rm
Due to (\ref{EZZ}), we see that (\ref{exact-rate}) is equivalent to
\[
\frac{P(F_n\leq x)-P(N\leq x)}{\kappa_3(F_n)} \to
 \frac{1}{6\sqrt{2\pi}}
\left(1-x^2\right)e^{-\frac{x^2}2}\quad\mbox{as $n\to\infty$}.
\]
Since each $F_n$ is centered, one also has that $\kappa_3(F_n) = E[F_n^3]$.}
\end{rem}
{\it Proof}.
We shall apply Theorem \ref{T : UberMain}.
Thanks to (\ref{aa3bis}), we get that
\[
\frac{\kappa_4(F_n)}{6}=\frac{E[F_n^4]-3}{6}=8\,\|f_n\otimes_1
f_n\|_{L^2(\R_+^2)}^2.
\]
By combining this identity with (\ref{lm-control-1d}) (here, it is worth observing
that $f_n\otimes_1 f_n$ is symmetric, so that the symmetrization $f_n\widetilde{\otimes}_1 f_n$
is immaterial), we see that
the quantity $\varphi(n)$ appearing in (\ref{IBPstein3})
is given by $\sqrt{\kappa_4(F_n)/{6}}$.
In particular, condition $(i)$ in Theorem \ref{T : UberMain} is met
(here, let us stress that one may show that $\kappa_4(F_n)>0$ for all $n$ by means of (\ref{aa3bis})).
On the other hand, since $F_n$ is a non-zero double integral, its law has a density
with respect to Lebesgue measure, according to Theorem \ref{shigekawa}.
This means that condition $(ii)$ in Theorem \ref{T : UberMain} is also in order.
Hence, it remains to check condition $(iii)$.
Assume that (\ref{EZZ}) holds.
Using (\ref{cumX2}) in the cases $p=3$
and $p=8$, we deduce that
\[
\frac{\kappa_3(F_n)}{\sqrt{\kappa_4(F_n)}}
= \frac{8\,\langle f_n, f_n\otimes_1
f_n\rangle_{L^2(\R_+^2)}}{\sqrt6\,\varphi(n)}
\]
and
\[
\frac{\kappa_8(F_n)}{\left(\kappa_4(F_n)\right)^2}
=\frac{ 17920 \|(f_n\otimes_1 f_n)\otimes_1
(f_n\otimes_1 f_n)\|^2_{L^2(\R_+^2)}}{\varphi(n)^4}.
\]
On the other hand, set
\[
Y_n=\frac{\frac12\|DF_n\|^2_{L^2(\R_+)} -1 }{\varphi(n)}.
\]
By (\ref{aa1}), we have $\frac12\|DY_n\|_{L^2(\R_+)}^2 -1 = 2\,I^B_2(f_n\otimes_1 f_n)$.
Therefore, by (\ref{lm-control-1d}), we get that
\begin{eqnarray*}
E\left[\left(\frac12\|DY_n\|_{L^2(\R_+)}^2 -1\right)^2\right]&=&\frac{128}{\varphi(n)^4}\|(f_n\otimes_1 f_n)\otimes_1(f_n\otimes_1 f_n)\|_{L^2(\R_+^2)}\\
&=&\frac{\kappa_8(F_n)}{140\left(\kappa_4(F_n)\right)^2}
\to 0 \quad\mbox{as $n\to\infty$}.
\end{eqnarray*}
Hence, by Theorem \ref{T : NPNOPTP}, we deduce that
$Y_n\overset{{\rm Law}}{\to}\mathcal{N}(0,1)$. We also have
\[
E[Y_nF_n]=\frac{4}{\varphi(n)}\langle f_n\otimes_1 f_n,f_n\rangle_{L^2(\R_+^2)}=
\frac{\sqrt{6}\,\kappa_3(F_n)}{2\sqrt{\kappa_4(F_n)}}\to
\frac{\alpha\sqrt6}{2}=:\rho\quad\mbox{as $n\to\infty$}.
\]
Therefore, to conclude that condition $(iii)$ in Theorem \ref{T : UberMain} holds true, it suffices to apply
Theorem \ref{PecTud}.
\qed

\bigskip

To give a concrete application of Proposition \ref{roc}, let us go back to the quadratic variation of fractional Brownian motion.
Let $B^H=(B^H_t)_{t\geq 0}$ be a fractional Brownian motion with Hurst
index $H\in(0,\frac12)$ and let
\[
F_n=
\frac{1}{\sigma_n}\sum_{k=0}^{n-1}\big[(B^H_{k+1}-B^H_{k})^2-1\big],
\]
where $\sigma_n>0$ is so that $E[F_n^2]=1$.
Recall from Theorem \ref{fBM} that
$\lim_{n\to\infty}\sigma_n^2/n=2\sum_{r\in\mathbb{Z}}\rho^2(r)<\infty$,
with $\rho:\mathbb{Z}\to\R_+$ given by (\ref{rho}); moreover,
there exists a
constant $c_H>0$ (depending only on $H$) such that, with $N\sim \mathcal{N}(0,1)$,
\begin{equation}\label{dtvfn}
d_{TV}(F_n,N)\leq \frac{c_H}{\sqrt{n}},\quad n\geq 1.
\end{equation}

The next two results aim to show that one can associate a {\it lower} bound to (\ref{dtvfn}).
We start by the following auxiliary result.

\begin{prop}\label{prop-cum-bm}
Fix an integer $s\geq 2$, let $F_n$ be as above and let $\rho$ be given by (\ref{rho}).
Recall that $H<\frac12$, so that $\rho\in \ell^1(\Z)$.
Then, the $s${\rm th} cumulant of $F_n$ behaves asymptotically as
\begin{equation}\label{touslescumulantsdansbreuermajor}
\kappa_s(F_n)\sim n^{1-s/2}\,2^{s/2-1}(s-1)!\,\frac{\langle \rho^{*(s-1)},\rho\rangle_{\ell^2(\mathbb{Z})}}{\|\rho\|^s_{\ell^2(\Z)}} \quad\mbox{as $n\to\infty$}.
\end{equation}
\end{prop}
{\it Proof}.
As in the proof of Theorem \ref{fBM}, we have that
$F_n\overset{\rm law}{=}I^B_2(f_n)$ with $f_n=\frac{1}{\sigma_n}\sum_{k=0}^{n-1} e_k^{\otimes 2}.$
Now, let us proceed with the proof. It is divided into several steps.

\medskip

\noindent{\it First step}. Using the formula (\ref{cumX2}) giving the cumulants of $F_n=I^B_2(f_n)$
as well as the very definition of the contraction $\otimes_1$,
we immediately check that
\begin{eqnarray*}
\kappa_s(F_n) = \frac{2^{s-1}(s-1)!}{\sigma_n^s}
\sum_{k_1,\ldots,k_s=0}^{n-1}\rho(k_s-k_{s-1})\ldots \rho(k_2-k_1)\rho(k_1-k_s).
\end{eqnarray*}

\noindent{\it Second step}. Since $H<\frac12$, we have that
$\rho\in \ell^1(\Z)$.  Therefore, by applying Young inequality repeatedly, we have
\begin{eqnarray*}
\|\,|\rho|^{*(s-1)}\|_{\ell^\infty(\mathbb{Z})}&\leq&
\|\rho\|_{\ell^1(\mathbb{Z})}
\|\,|\rho|^{*(s-2)}\|_{\ell^\infty(\mathbb{Z})}
\leq
\ldots
\leq
\|\rho\|_{\ell^1(\mathbb{Z})}^{s-1}<\infty.
\end{eqnarray*}
In particular, we have that
$\langle|\rho|^{*(s-1)},|\rho|\rangle_{\ell^2(\Z)}\leq\|\rho\|_{\ell^1(\mathbb{Z})}^{s}<\infty$.

\medskip

\noindent{\it Third step}.
Thanks to the result shown in the previous step, observe first that
\begin{eqnarray*}
\sum_{k_2,\ldots,k_{s}\in\mathbb{Z}}
\left|
\rho(k_2)\rho(k_2-k_3)\rho(k_3-k_4)\ldots\rho(k_{s-1}-k_s)\rho(k_s)\right|
=\langle |\rho|^{*(s-1)},|\rho|\rangle_{\ell^2(\mathbb{Z})}<\infty.
\end{eqnarray*}
Hence, one can apply dominated convergence to get, as $n\to\infty$, that
\begin{eqnarray}
&&\frac{\sigma_n^s}{2^{s-1}(s-1)!\,n}\,\kappa_s(F_n)\notag\\
&=&
\frac1n
\sum_{k_1=0}^{n-1}\,\,\,
\sum_{k_2,\ldots,k_s=-k_1}^{n-1 -k_1}
\rho(k_2)\rho(k_2-k_3)\rho(k_3-k_4)\ldots\rho(k_{s-1}-k_s)\rho(k_s)\notag\\
&=&
\sum_{k_2,\ldots,k_{s}\in\mathbb{Z}}
\rho(k_2)\rho(k_2-k_3)\rho(k_3-k_4)\ldots\rho(k_{s-1}-k_s)\rho(k_s)\notag\\
&&\hskip1cm\times\left[1\wedge\left(1-\frac{\max\{k_2,\ldots, k_s\}}{n}\right)-0\vee\left(\frac{\min\{k_2,\ldots, k_s\}}{n}\right)\right]
{\bf 1}_{\{|k_2|< n,\ldots,|k_s|< n\}}\notag\\
&\to&
\sum_{k_2,\ldots,k_{s}\in\mathbb{Z}}
\rho(k_2)\rho(k_2-k_3)\rho(k_3-k_4)\ldots\rho(k_{s-1}-k_s)\rho(k_s)=
\langle \rho^{*(s-1)},\rho\rangle_{\ell^2(\mathbb{Z})}.\notag\\
\label{lilipuce}
\end{eqnarray}
Since $\sigma_n\sim \sqrt{2n}\,\|\rho\|_{\ell^2(\Z)}$ as $n\to\infty$, the desired conclusion follows.
\qed

\begin{cor}
Let $F_n$ be as above (with $H<\frac12$), let $N\sim\mathcal{N}(0,1)$, and let $\rho$ be given by (\ref{rho}).
Then, for all $x\in\R$, we have
\[
\sqrt{n}\big(P(F_n\leq x)-P(N\leq x)\big) \to
\frac{\langle \rho^{*2},\rho\rangle_{\ell^2(\mathbb{Z})}}{3\|\rho\|^2_{\ell^2(\Z)}}
\,(1-x^2)\,e^{-\frac{x^2}2}
\quad\mbox{as $n\to\infty$}.
\]
In particular, we deduce that there exists $d_H>0$ such that
\begin{equation}\label{lowerboundz}
\frac{d_H}{\sqrt{n}}\leq \big|P(F_n\leq 0)-P(N\leq 0)\big|\leq
d_{TV}(F_n,N), \quad n\geq 1.
\end{equation}
\end{cor}
{\it Proof}. The desired conclusion follows immediately by combining Propositions \ref{roc} and \ref{prop-cum-bm}.
\qed

\bigskip

By paying closer attention to the used estimates, one may actually show that (\ref{lowerboundz}) holds true for any $H<\frac58$
(not only $H<\frac12$). See \cite[Theorem 9.5.1]{book-malliavinstein} for the details.

\bigskip

{\bf To go further}.
The paper \cite{exact} contains several other examples of application of Theorem \ref{T : UberMain} and Proposition \ref{roc}.
In the reference \cite{BBNP}, one shows that the deterministic sequence
\[
\max\{E[F_n^3],E[F_n^4]-3\},\quad n\geq 1,
\]
completely characterizes the rate of convergence
(with respect to smooth distances) in CLTs involving chaotic random variables.

\noindent
\section{An extension to the Poisson space (following the invited talk by Giovanni Peccati)}\label{sec:gio}

Let $B=(B_t)_{t\geq 0}$ be a Brownian motion,
let $F$ be any centered element of $\mathbb{D}^{1,2}$ and let $N\sim\mathcal{N}(0,1)$.
We know from Theorem \ref{nou-pec} that
\begin{equation}\label{dtvaq}
d_{TV}(F,N)\leq 2\,E[|1-\langle DF,-DL^{-1}F\rangle_{L^2(\R_+)}|].
\end{equation}
The aim of this section, which follows \cite{peccati,PSTU}, is to explain how to deduce inequalities of the type (\ref{dtvaq}),
when $F$ is a regular functional of a {\it Poisson measure} $\eta$
and when the target law $N$ is either Gaussian or Poisson.

We first need to introduce the basic concepts in this framework.

\bigskip
\noindent
{\bf Poisson measure}.
In what follows, we shall use the symbol $Po(\lambda)$ to indicate the Poisson distribution of parameter $\lambda>0$ (that is,
$\mathcal{P}_\lambda\sim Po(\lambda)$ if and only if $P(\mathcal{P}_\lambda=k)=e^{-\lambda}\frac{\lambda^k}{k!}$ for
all $k\in\N$), with
the convention that $Po(0)=\delta_0$ (Dirac mass at 0).
Set $A=\R^d$ with $d\geq 1$, let $\mathcal{A}$ be the Borel $\sigma$-field on $A$, and
let $\mu$ be a positive, $\sigma$-finite and atomless measure over $(A,\mathcal{A})$.
We set $\mathcal{A}_\mu=\{B\in\mathcal{A}:\,\mu(B)<\infty\}$.
\begin{defi}
A Poisson measure $\eta$ with control $\mu$ is an object of the form $\{\eta(B)\}_{B\in\mathcal{A}_\mu}$ with the following
features:
\begin{enumerate}
\item[(1)] for all $B\in\mathcal{A}_\mu$, we have $\eta(B)\sim Po(\mu(B))$.
\item[(2)] for all $B,C\in\mathcal{A}_\mu$ with $B\cap C\neq\emptyset$, the random variables
$\eta(B)$ and $\eta(C)$ are independent.
\end{enumerate}
Also, we note $\widehat{\eta}(B)=\eta(B)-\mu(B)$.
\end{defi}

\begin{rem}{\rm
\begin{enumerate}

\item As a simple example, note that for $d=1$ and $\mu = \lambda\times {\rm Leb}$ (with `${\rm Leb}$' the Lebesgue measure) the process
$\{\eta([0,t])\}_{t\geq 0}$ is nothing but a Poisson process with intensity $\lambda$.

\item Let $\mu$ be a $\sigma$-finite atomless measure over $(A,\mathcal{A})$, and observe that this implies that there exists a sequence of disjoint sets $\{A_j : j\geq 1\}\subset \mathcal{A}_\mu$ such that $\cup_j A_j = A$. For every $j=1,2,...$ belonging to the set $J_0$ of those indices such that $\mu(A_j)>0$ consider the following objects: $X^{(j)} =\{X_k^{(j)} : k\geq 1\}$ is a sequence of i.i.d. random variables with values in $A_j$ and with common distribution $\frac{\mu_{|_ {A_j}}}{\mu(A_j)}$;  $P_j$ is a Poisson random variable with parameter $\mu(A_j)$. Assume moreover that : (i) $X^{(j)}$ is independent of $X^{(k)}$ for every $k\neq j$, (ii) $P_j$ is independent of $P_k$ for every $k\neq j$, and (iii) the classes $\{X^{(j)}\}$ and $\{P_j\}$ are independent. Then, it is a straightforward computation to verify that the random point measure
\[
\eta(\cdot) = \sum_{j\in J_0} \sum_{k=1}^{P_j} \delta_{X^{(j)}_k}(\cdot),
\]
where $\delta_x$ indicates the Dirac mass at $x$ and $\sum_{k=1}^0 =0$ by convention, is a a Poisson random measure with control $\mu$. See e.g. \cite[Section 1.7]{penrosebook}.

\end{enumerate}
}
\end{rem}

\bigskip
\noindent
{\bf Multiple integrals and chaotic expansion}.
As a preliminary remark, let us observe that $E[\widehat{\eta}(B)]=0$ and $E[\widehat{\eta}(B)^2]=\mu(B)$ for all $B\in\mathcal{A}_\mu$.
For any $q\geq 1$, set $L^2(\mu^q)=L^2(A^q,\mathcal{A}^q,\mu^q)$.
We want to appropriately define
\[
I^{\widehat{\eta}}_q(f)=\int_{A^q}f(x_1,\ldots,x_q)\widehat{\eta}(dx_1)\ldots\widehat{\eta}(dx_q)
\]
when $f\in L^2(\mu^q)$.
To reach our goal, we proceed in a classical way.
We first consider the subset $\mathcal{E}(\mu^q)$ of simple functions, which is defined as
\[
\mathcal{E}(\mu^q)={\rm span}\left\{
{\bf 1}_{B_1}\otimes\ldots\otimes{\bf 1}_{B_q},\,\mbox{with}\,
B_1,\ldots,B_q\in\mathcal{A}_\mu\,\mbox{such that}\,B_i\cap B_j=\emptyset\mbox{ for all $i\neq j$}
\right\}.
\]
When $f={\bf 1}_{B_1}\otimes\ldots\otimes{\bf 1}_{B_q}$ with $B_1,\ldots,B_q\in\mathcal{A}_\mu$ such that $B_i\cap B_j=\emptyset$ for all $i\neq j$, we naturally set
\[
I^{\widehat{\eta}}_q(f):=\widehat{\eta}(B_1)\ldots\widehat{\eta}(B_q)=\int_{A^q}f(x_1,\ldots,x_q)\widehat{\eta}(dx_1)\ldots
\widehat{\eta}(dx_q).
\]
(For such a simple function $f$, note that the right-hand side in the previous formula makes perfectly sense by considering $\widehat{\eta}$
as a signed measure.)
We can extend by linearity the definition of $I^{\widehat{\eta}}_q(f)$ to any $f\in \mathcal{E}(\mu^q)$.
It is then a simple computation to check that
\[
E[I^{\widehat{\eta}}_p(f)I^{\widehat{\eta}}_q(g)]=p!\delta_{p,q}\,\langle \widetilde{f},\widetilde{g}\rangle_{L^2(\mu^p)}
\]
for all $f\in \mathcal{E}(\mu^p)$ and $g\in\mathcal{E}(\mu^q)$, with $\widetilde{f}$ (resp. $\widetilde{g}$)
the symmetrization of $f$ (resp. $g$) and $\delta_{p,q}$ the Kronecker symbol.
Since $\mathcal{E}(\mu^q)$ is dense in $L^2(\mu^q)$ (it is precisely here that the fact that $\mu$ has no atom is crucial!), we can define
$I^{\widehat{\eta}}_q(f)$ by isometry to any $f\in L^2(\mu^q)$.
Relevant properties of $I^{\widehat{\eta}}_q(f)$ include $E[I^{\widehat{\eta}}_q(f)]=0$, $I^{\widehat{\eta}}_q(f)=I^{\widehat{\eta}}_q(\widetilde{f})$ and (importantly!) the fact that $I^{\widehat{\eta}}_q(f)$ is a {\it true}
multiple integral when $f\in\mathcal{E}(\mu^q)$.

\begin{defi}
Fix $q\geq 1$. The set of random variables of the form $I^{\widehat{\eta}}_q(f)$ is called the $q$th Poisson-Wiener chaos.
\end{defi}

In this framework, we have an analogue of the chaotic decomposition (\ref{E}) -- see e.g. \cite[Corollary 10.0.5]{PeTa} for a proof.
\begin{thm}
For all $F\in L^2(\sigma\{\eta\})$ (that is, for all random variable $F$ which is square integrable and measurable with respect to $\eta$),
we have
\begin{equation}\label{chaopoisson}
F=E[F]+\sum_{q=1}^{\infty} I^{\widehat{\eta}}_q(f_q),
\end{equation}
where the kernels $f_q$ are ($\mu^q$-a.e.) symmetric elements of $L^2(\mu^q)$ and are
uniquely determined by $F$.
\end{thm}

\bigskip
\noindent
{\bf Multiplication formula and contractions}.
When $f\in \mathcal{E}(\mu^p)$ and $g\in\mathcal{E}(\mu^q)$ are symmetric, we define, for all $r=0,\ldots,p\wedge q$ and
$l=0,\ldots,r$:
\begin{eqnarray*}
&&f\star_r^l g (x_1,\ldots,x_{p+q-r-l})\\
&=&\int_{A^l} f(y_1,\ldots,y_l,x_1,\ldots,x_{r-l},x_{r-l+1},\ldots,x_{p-l})g(y_1,\ldots,y_l,x_{1},\ldots,x_{r-l},x_{p-l+1},\ldots,x_{p+q-r-l})\\
&&\hskip12.5cm\times \mu(dy_1)\ldots \mu(dy_l).
\end{eqnarray*}
We then have the following product formula, compare with (\ref{multiplication}).
\begin{thm}[Product formula]
Let $p,q\geq 1$ and let $f\in\mathcal{E}(\mu^p)$ and $g\in\mathcal{E}(\mu^q)$ be symmetric. Then
\[
I^{\widehat{\eta}}_p(f)I^{\widehat{\eta}}_q(g)=\sum_{r=0}^{p\wedge q} r!\binom{p}{r}\binom{q}{r}\sum_{l=0}^r \binom{r}{l}I^{\widehat{\eta}}_{p+q-r-l}\left(\widetilde{f\star_r^l g}\right).
\]
\end{thm}
{\it Proof}.
Recall that, when dealing with functions in $\mathcal{E}(\mu^p)$,
$I^{\widehat{\eta}}_p(f)$ is a {\it true} multiple integral (by seeing $\widehat{\eta}$ as a signed measure). We deduce
\[
I^{\widehat{\eta}}_p(f)I^{\widehat{\eta}}_q(g)=\int_{A^{p+q}}f(x_1,\ldots,x_p)g(y_1,\ldots,y_q)\widehat{\eta}(dx_1)\ldots \widehat{\eta}(dx_p)
\widehat{\eta}(dy_1)\ldots \widehat{\eta}(dy_q).
\]
By definition of $f$ (the same applies for $g$), we have that $f(x_1,\ldots,x_p)=0$ when $x_i=x_j$ for some $i\neq j$.
Consider $r=0,\ldots,p\wedge q$, as well as pairwise disjoint indices $i_1,\ldots,i_r\in\{1,\ldots,p\}$ and
pairwise disjoint indices $j_1,\ldots,j_r\in\{1,\ldots,q\}$.
Set $\{k_1,\ldots,k_{p-r}\}=\{1,\ldots,p\}\setminus\{i_1,\ldots,i_r\}$
and $\{l_1,\ldots,l_{q-r}\}=\{1,\ldots,q\}\setminus\{j_1,\ldots,j_r\}$.
We have, since $\mu$ is atomless and using $\widehat{\eta}(dx)=\eta(dx)-\mu(dx)$,
\begin{eqnarray*}
&&\int_{A^{p+q}}f(x_1,\ldots,x_p)g(y_1,\ldots,y_q)
{\bf 1}_{\{
x_{i_1}=y_{j_1},\ldots,x_{i_r}=y_{j_r}
\}}
\widehat{\eta}(dx_1)\ldots \widehat{\eta}(dx_p)
\widehat{\eta}(dy_1)\ldots \widehat{\eta}(dy_q)\\
&=&\int_{A^{p+q-2r}}f(x_{k_1},\ldots,x_{k_{p-r}},x_{i_1},\ldots,x_{i_r})g(y_{l_1},\ldots,y_{l_{q-r}},x_{i_1},\ldots,x_{i_r})\\
&&\hskip4cm\times\widehat{\eta}(dx_{k_1})\ldots\widehat{\eta}(dx_{k_{p-r}})\widehat{\eta}(dy_{l_1})\ldots \widehat{\eta}(dy_{l_{q-r}})
\eta(dx_{i_1})\ldots \eta(dx_{i_r})\\
&=&\int_{A^{p+q-2r}}f(x_1,\ldots,x_{p-r},a_{1},\ldots,a_{r})g(y_{1},\ldots,y_{q-r},a_{1},\ldots,a_{r})\\
&&\hskip4cm\times\widehat{\eta}(dx_{1})\ldots\widehat{\eta}(dx_{p-r})\widehat{\eta}(dy_{1})\ldots \widehat{\eta}(dy_{q-r})
\eta(da_{1})\ldots \eta(da_{r}).
\end{eqnarray*}
By decomposing over the hyperdiagonals $\{x_i=y_j\}$, we deduce that
\begin{eqnarray*}
I^{\widehat{\eta}}_p(f)I^{\widehat{\eta}}_q(g)&=&\sum_{r=0}^{p\wedge q}r!\binom{p}{r}\binom{q}{r}\int_{A^{p+q-2r}}f(x_1,\ldots,x_{p-r},a_{1},\ldots,a_{r})g(y_{1},\ldots,y_{q-r},a_{1},\ldots,a_{r})\\
&&\hskip4cm\times\widehat{\eta}(dx_{1})\ldots\widehat{\eta}(dx_{p-r})\widehat{\eta}(dy_{1})\ldots \widehat{\eta}(dy_{q-r})
\eta(da_{1})\ldots \eta(da_{r}),
\end{eqnarray*}
and we get the desired conclusion by using the relationship
\[
\eta(da_{1})\ldots \eta(da_{r})=\big(\widehat{\eta}(da_1)+\mu(da_1)\big)\ldots\big(\widehat{\eta}(da_r)+\mu(da_r)\big).
\]
\qed

\bigskip
\noindent
{\bf Malliavin operators}. Each time we deal with a random element $F$ of $L^2(\{\sigma(\eta)\})$, in what follows we always consider its chaotic expansion (\ref{chaopoisson}).
\begin{defi}
1. Set ${\rm Dom}D=\{F\in L^2(\sigma\{\eta\}):\,\sum qq!\|f_q\|^2_{L^2(\mu^q)}<\infty\}$. If $F\in{\rm Dom}D$, we set
\[
D_t F =\sum_{q=1}^\infty q I^{\widehat{\eta}}_{q-1}(f_q(\cdot,t)),\quad t\in A.
\]
The operator $D$ is called the Malliavin derivative.\\
2. Set ${\rm Dom}L=\{F\in L^2(\sigma\{\eta\}):\,\sum q^2q!\|f_q\|^2_{L^2(\mu^q)}<\infty\}$.
If $F\in{\rm Dom}L$, we set
\[
LF=-\sum_{q=1}^\infty qI^{\widehat{\eta}}_q(f_q).
\]
The operator $L$ is called the generator of the Ornstein-Uhlenbeck semigroup.\\
3. If $F\in L^2(\sigma\{\eta\})$, we set
\[
L^{-1}F=-\sum_{q=1}^\infty \frac1q I^{\widehat{\eta}}_q(f_q).
\]
The operator $L^{-1}$ is called the pseudo-inverse of $L$.
\end{defi}
It is readily checked that $LL^{-1}F=F-E[F]$ for $F\in L^2(\sigma\{\eta\})$.
Moreover, proceeding {\sl mutatis mutandis} as in the proof of Theorem \ref{ivangioformulatheorem}, we get
the following result.
\begin{prop}
Let $F\in L^2(\sigma\{\eta\})$ and let $G\in{\rm Dom}D$.
Then
\begin{equation}\label{covfg}
{\rm Cov}(F,G)=E[\langle DG,-DL^{-1}F\rangle_{L^2(\mu)}].
\end{equation}
\end{prop}
The operator $D$ does not satisfy the chain rule. Instead, it admits an `add-one cost' representation which
plays an identical role.
\begin{thm}[Nualart, Vives, 1990; see \cite{NuVi}]
Let $F\in{\rm Dom}D$. Since $F$ is measurable with respect to $\eta$, we can view it as $F=F(\eta)$ with a slight abuse
of notation.
Then
\begin{equation}\label{nuvi}
D_t F = F(\eta+\delta_t) - F(\eta),\quad t\in A,
\end{equation}
where $\delta_t$ stands for the Dirac mass at $t$.
\end{thm}
{\it Proof}. By linearity and approximation, it suffices to prove the claim for $F=I^{\widehat{\eta}}_q(f)$, with $q\geq 1$ and $f\in\mathcal{E}(\mu^q)$ symmetric.
In this case, we have
\[
F(\eta+\delta_t)=\int_{A^q} f(x_1,\ldots,x_q)\big(\widehat{\eta}(dx_1)+\delta_t(dx_1)\big)\ldots
\big(\widehat{\eta}(dx_q)+\delta_t(dx_q)\big).
\]
Let us expand the integrator. Each member of such an expansion such that there is strictly more than one Dirac mass in the resulting expression gives a contribution equal to zero,
since $f$ vanishes on diagonals. We therefore deduce that
\begin{eqnarray*}
F(\eta+\delta_t)&=&F(\eta)+\sum_{l=1}^q \int_{A^q} f(x_1,\ldots,x_{l-1},t,x_{l+1},\ldots,x_q)\widehat{\eta}(dx_1)\ldots
\widehat{\eta}(dx_{l-1})\widehat{\eta}(dx_{l+1})\ldots\widehat{\eta}(dx_q)\\
&=&F(\eta)+qI^{\widehat{\eta}}_{q-1}(f(t,\cdot))\quad\mbox{by symmetry of $f$}\\
&=&F(\eta)+D_tF.
\end{eqnarray*}
\qed
As an immediate corollary of the previous theorem, we get the formula
\[
D_t (F^2)= (F+D_tF)^2-F^2 = 2F\,D_tF + (D_tF)^2,\quad t\in A,
\]
which shows how $D$ is far from satisfying the chain rule (\ref{chainrule}).

\bigskip
\noindent
{\bf Gaussian approximation}.
It happens that it is the following distance which is appropriate in our framework.

\begin{defi}
The Wasserstein distance between the laws of two real-valued random variables $Y$ and $Z$ is defined by
\begin{equation}\label{wasser}
d_{W}(Y,Z)=\sup_{h\in {\rm Lip}(1)} \big| E[h(Y)]-E[h(Z)]\big|,
\end{equation}
where ${\rm Lip}(1)$ stands for the set of Lipschitz functions $h:\R\to\R$ with constant 1.
\end{defi}

Since we are here dealing with Lipschitz functions $h$, we need a suitable version of Stein's lemma.
Compare with Lemma \ref{stein}.
\begin{lemma}[Stein, 1972; see \cite{Stein_orig}]\label{stein2}
Suppose $h:\R\to\R$ is a Lipschitz constant with constant 1. Let $N\sim\mathcal{N}(0,1)$. Then, there exists
a solution to the equation
\[
f'(x)-xf(x)=h(x)-E[h(N)],\quad x\in\R,
\]
that satisfies $\|f'\|_\infty\leq \sqrt{\frac2{\pi}}$ and $\|f''\|_\infty\leq 2$.
\end{lemma}
{\it Proof}.
Let us recall that, according to Rademacher's theorem, a function which is Lipschitz continuous on $\R$ is almost everywhere differentiable.
Let $f:\R\to\R$ be the (well-defined!) function given by
\begin{equation}\label{godspell}
f(x)=
-\int_0^\infty \frac{e^{-t}}{\sqrt{1-e^{-2t}}}
E[h(e^{-t}x+\sqrt{1-e^{-2t}}N)N]dt.
\end{equation}
By dominated convergence we have that $f_h\in\mathcal{C}^1$ with
\begin{eqnarray*}
f'(x)&=&-\int_0^\infty \frac{e^{-2t}}{\sqrt{1-e^{-2t}}}
E[h'(e^{-t}x+\sqrt{1-e^{-2t}}N)N]dt.
\end{eqnarray*}
We deduce, for any $x\in\R$,
\begin{equation}\label{bounn}
|f'(x)|\leq E|N|\int_0^\infty \frac{e^{-2t}}{\sqrt{1-e^{-2t}}}dt = \sqrt{\frac{2}{\pi}}.
\end{equation}
Now, let $F:\R\to\R$ be the function given by
\[
F(x)=\int_0^\infty
E[h(N)-h(e^{-t}x+\sqrt{1-e^{-2t}}N)]
dt,\quad x\in\R.\]
Observe that $F$ is well-defined since
$h(N)-h(e^{-t}x+\sqrt{1-e^{-2t}}N)$ is integrable
due to
\begin{eqnarray*}
\big|h(N)-h(e^{-t}x+\sqrt{1-e^{-2t}}N)\big|&\leq& e^{-t}|x|+\big(1-\sqrt{1-e^{-2t}}\big)|N|\\
&\leq& e^{-t} |x|+e^{-2t}|N|,
\end{eqnarray*}
where the last inequality follows from $1-\sqrt{1-u}=u/(\sqrt{1-u}+1)\leq
u$
if $u\in[0,1]$.
By dominated convergence, we immediately see that
$F$ is differentiable with
\[
F'(x)=-\int_0^\infty e^{-t}\,E[h'(e^{-t}x+\sqrt{1-e^{-2t}}N)]dt.
\]
By integrating by parts, we see that
$F'(x)=f(x)$.
Moreover, by using the notation introduced in Section \ref{malliavin}, we can write
\begin{eqnarray*}
&&f'(x)-xf(x)\\
&=&LF(x),\quad\mbox{by decomposing in Hermite polynomials, since $LH_q=-qH_q=H_q''-XH_q'$}\\
&=&-\int_0^\infty LP_th(x)dt,\quad\mbox{since $F(x)=\int_0^\infty \big(E[h(N)]-P_th(x)\big)dt$}\\
&=&-\int_0^\infty \frac{d}{dt}P_th(x)dt\\
&=&P_0h(x)-P_\infty h(x)=h(x)-E[h(N)].
\end{eqnarray*}
This proves the claim for $\|f'\|_\infty$. The claim for $\|f''\|_\infty$ is a bit more difficult to achieve;
we refer to Stein \cite[pp. 25-28]{Stein_orig} to keep the length of this survey within bounds.
\qed

\bigskip

We can now derive a bound for the Gaussian approximation of any centered
element $F$ belonging to ${\rm Dom}D$, compare with (\ref{dtvaq}).

\begin{thm}[Peccati, Sol\'e, Taqqu, Utzet, 2010; see \cite{PSTU}]\label{PSTU}
Consider $F\in{\rm Dom}D$ with $E[F]=0$.
Then, with $N\sim \mathcal{N}(0,1)$,
\[
d_{W}(F,N)
\leq \sqrt{\frac{2}{\pi}}\,E\left[\big|1-\langle DF,-DL^{-1}F\rangle_{L^2(\mu)}\big|\right]
+E\left[\int_A (D_tF)^2|D_tL^{-1}F|\mu(dt)\right].
\]
\end{thm}
\noindent{\it Proof}.
Let $h\in {\rm Lip}(1)$ and let $f$ be the function of Lemma \ref{stein2}.
Using (\ref{nuvi}) and a Taylor formula, we can write
\[
D_tf(F)=f(F+D_tF)-f(F)=f'(F)D_tF+R(t),
\]
with $|R(t)|\leq \frac12\|f''\|_\infty(D_tF)^2
\leq (D_tF)^2$.
We deduce, using (\ref{covfg}) as well,
\begin{eqnarray*}
E[h(F)]-E[h(N)]&=&E[f'(F)]-E[Ff(F)]=E[f'(F)]-E[\langle Df(F),-DL^{-1}F\rangle_{L^2(\mu)}]\\
&=&E[f'(F)(1-\langle DF,-DL^{-1}F\rangle_{L^2(\mu)})]+\int_A (-D_tL^{-1}F)R(t)\mu(dt).
\end{eqnarray*}
Consequently, since $\|f'\|_\infty\leq \sqrt{\frac{2}{\pi}}$,
\begin{eqnarray*}
d_W(F,N)&=&\sup_{h\in {\rm Lip}(1)} |E[h(F)]-E[h(N)]|\\
&\leq&
\sqrt{\frac{2}{\pi}}\,E\left[\big|1-\langle DF,-DL^{-1}F\rangle_{L^2(\mu)}\big|\right]
+E\left[\int_A (D_tF)^2|D_tL^{-1}F|\mu(dt)\right].
\end{eqnarray*}
\qed

\bigskip
\noindent
{\bf Poisson approximation}.
To conclude this section, we will prove a very interesting result, which may be seen as a Poisson counterpart of Theorem \ref{PSTU}.
\begin{thm}[Peccati, 2012; see \cite{peccati}]\label{peccatidelamathematique}
Let $F\in {\rm Dom}D$ with $E[F]=\lambda>0$ and $F$ taking its values in $\mathbb{N}$.
Let $\mathcal{P}_\lambda\sim Po(\lambda)$. Then,
\begin{eqnarray}
&&\sup_{C\subset\N} \big|P(F\in C)-P(\mathcal{P}_\lambda\in C)\big|\label{popopo}\\
&\leq&\frac{1-e^{-\lambda}}{\lambda}\,E|\lambda-\langle DF,-DL^{-1}F\rangle_{L^2(\mu)}|
+\frac{1-e^{-\lambda}}{\lambda^2}\,E\int |D_tF(D_tF-1)D_tL^{-1}F|\mu(dt).\notag
\end{eqnarray}
\end{thm}
Just as a mere illustration, consider the case where $F=\eta(B)=I^{\widehat{\eta}}_1({\bf 1}_B)$ with $B\in\mathcal{A}_\mu$.
We then have $DF=-DL^{-1}F={\bf 1}_B$, so that $\langle DF,-DL^{-1}F\rangle_{L^2(\mu)}=\int {\bf 1}_Bd\mu=\mu(B)$
and $DF(DF-1)=0$ a.e. The right-hand side of (\ref{popopo}) is therefore zero, as it was expected since $F\sim Po(\lambda)$.

\bigskip

During the proof of Theorem \ref{peccatidelamathematique}, we shall use an analogue of Lemma \ref{stein} in the Poisson context, which reads as follows.

\begin{lemma}[Chen, 1975; see \cite{chen}]\label{chenstein}
Let $C\subset\N$, let $\lambda>0$ and let $\mathcal{P}_\lambda\sim Po(\lambda)$.
The equation with unknown $f:\N\to\R$,
\begin{equation}\label{cchen}
\lambda\,f(k+1)-kf(k)={\bf 1}_C(k)-P(\mathcal{P}_\lambda\in C),\quad k\in\N,
\end{equation}
admits a unique solution such that $f(0)=0$, denoted by $f_C$. Moreover, by setting $\Delta f(k)=f(k+1)-f(k)$,
we have $\|\Delta f_C\|_\infty\leq \frac{1-e^{-\lambda}}{\lambda}$ and $\|\Delta^2f_C\|_\infty\leq\frac{2}{\lambda}\|\Delta f_C\|_\infty$.
\end{lemma}
{\it Proof}.  We only provide a proof for the bound on $\Delta f_C$; the estimate on $\Delta^2f_C$ is proved e.g. by Daly in \cite{daly}. Multiplying both sides of (\ref{cchen}) by $\lambda^k/k!$ and summing up yields that, for every $k\geq 1$,
\begin{eqnarray}
f_C(k) &=& \frac{(k-1)!}{\lambda^k} \sum_{r=0}^{k-1} \frac{\lambda^r}{r!}[ {\bf 1}_C(r) - P(\mathcal{P}_\lambda\in C)] \label{e:cs1}\\
&=& \sum_{j\in C} f_{\{j\}} (k) \label{e:cs2}\\
&=& -f_{C^c} (k) \label{e:cs3} \\
&=& - \frac{(k-1)!}{\lambda^k} \sum_{r=k}^{\infty} \frac{\lambda^r}{r!}[ {\bf 1}_C(r) - P(\mathcal{P}_\lambda\in C)] \label{e:cs4},
\end{eqnarray}
where $C^c$ denotes the complement of $C$ in $\mathbb{N}$.
(Identity (\ref{e:cs2}) comes from the additivity property of $C\mapsto f_C$, identity (\ref{e:cs3}) is because
$f_\N\equiv 0$ and identity (\ref{e:cs3}) is due to
\[
\sum_{r=0}^{\infty} \frac{\lambda^r}{r!}[ {\bf 1}_C(r) - P(\mathcal{P}_\lambda\in C)]=E[{\bf 1}_C(\mathcal{P}_\lambda)
-E[{\bf 1}_C(\mathcal{P}_\lambda)]]=0.\big)
  \]
  Since $f_C(k) - f_C(k+1) = f_{C^c}(k+1) - f_{C^c}(k)$ (due to (\ref{e:cs3})), it is sufficient to prove that, for every $k\geq 1$ and every $C\subset \mathbb{N}$, $f_{C}(k+1) -f_{C}(k)\leq (1 - e^{-\lambda})/\lambda$. One has the following fact: for every $j\geq 1$ the mapping $k\mapsto f_{\{j\}}(k)$ is negative and decreasing for $k=1,...,j$ and positive and decreasing for $k\geq j+1$. Indeed, we use (\ref{e:cs1}) to deduce that, if $1\leq k\leq j$,
\[
f_{\{j\}} (k) = -e^{-\lambda} \frac{\lambda^j}{j!} \sum_{r=1}^{k} \lambda^{-r} \frac{(k-1)!}{(k-r)!} \quad\text{(which is negative and decreasing in $k$)},
\]
whereas (\ref{e:cs4}) implies that, if $k\geq j+1$,
\[
f_{\{j\}} (k) = e^{-\lambda} \frac{\lambda^j}{j!} \sum_{r=0}^{\infty} \lambda^{r} \frac{(k-1)!}{(k+r)!}\quad\text{(which is positive and decreasing in $k$)}.
\]
Using (\ref{e:cs2}), one therefore infers that $f_{C}(k+1) -f_{C}(k)\leq f_{\{k\}}(k+1) -f_{\{k\}}(k)$, for every $k\geq 0$. Since
\begin{eqnarray*}
f_{\{k\}}(k+1) -f_{\{k\}}(k)& =& e^{-\lambda} \left[ \sum_{r=0}^{k-1} \frac{\lambda^r}{r!k} + \sum_{r=k+1}^{\infty} \frac{\lambda^{r-1}}{r!}  \right]
=\frac{e^{-\lambda}}{\lambda}
\left[ \sum_{r=1}^{k} \frac{\lambda^r}{r!}\times\frac{r}{k} + \sum_{r=k+1}^{\infty} \frac{\lambda^{r}}{r!}  \right]
\\
&\leq& \frac{1-e^{-\lambda}}{\lambda},
\end{eqnarray*}
the proof is concluded.
\qed

\bigskip

We are now in a position to prove Theorem \ref{peccatidelamathematique}.

\bigskip

\noindent
{\it Proof of Theorem \ref{peccatidelamathematique}}. The main ingredient is the following simple inequality, which is a kind
of Taylor formula: for all $k,a\in \N$,
\begin{equation}\label{ingremain}
\big|f(k)-f(a)-\Delta f(a)(k-a)\big|\leq \frac12 \|\Delta^2 f\|_\infty|(k-a)(k-a-1)|.
\end{equation}
Assume for the time being that (\ref{ingremain}) holds true and fix $C\subset\N$. We have, using Lemma \ref{chenstein}
and then (\ref{covfg})
\begin{eqnarray*}
\big|P(F\in C)-P(\mathcal{P}_\lambda\in C)\big|&=&\big|E[\lambda f_C(F+1)]-E[Ff_C(F)]\big| \\
&=&\big|\lambda E[\Delta f_C(F)]-E[(F-\lambda)f_C(F)]\big|\\
&=&\big|\lambda E[\Delta f_C(F)]-E[\langle Df_C(F),-DL^{-1}F\rangle_{L^2(\mu)}]\big|.
\end{eqnarray*}
Now, combining (\ref{nuvi}) with (\ref{ingremain}), we can write
\[
D_tf_C(F)=\Delta f_C(F)D_t F+ S(t),
\]
with $S(t)\leq \frac12\|\Delta^2 f_C\|_\infty |D_tF(D_tF-1)|\leq \frac{1-e^{-\lambda}}{\lambda^2}|D_tF(D_tF-1)|$,
see indeed Lemma \ref{chenstein} for the last inequality.
Putting all these bounds together and since $\|\Delta f_C\|_\infty\leq \frac{1-e^{-\lambda}}{\lambda}$
by
Lemma \ref{chenstein}, we get the desired conclusion.

So, to conclude the proof, it remains to show that (\ref{ingremain}) holds true. Let us first assume that $k\geq a+2$.
We then have
\begin{eqnarray*}
f(k)&=&f(a)+\sum_{j=a}^{k-1}\Delta f(j)=f(a)+\Delta f(a)(k-a)+\sum_{j=a}^{k-1}(\Delta f(j)-\Delta f(a))\\
&=&f(a)+\Delta f(a)(k-a)+\sum_{j=a}^{k-1}\sum_{l=a}^{j-1}\Delta^2 f(l)
=f(a)+\Delta f(a)(k-a)+\sum_{l=a}^{k-2}\Delta^2 f(l)(k-l-1),
\end{eqnarray*}
so that
\begin{eqnarray*}
\left|f(k)-f(a)-\Delta f(a)(k-a)\right|\leq \|\Delta^2f\|_\infty \sum_{l=a}^{k-2}(k-l-1)
=\frac12\|\Delta^2f\|_\infty (k-a)(k-a-1),
\end{eqnarray*}
that is, (\ref{ingremain}) holds true in this case. When $k=a$ or $k=a+1$, (\ref{ingremain}) is obviously true.
Finally, consider the case $k\leq a-1$.
We have
\begin{eqnarray*}
f(k)&=&f(a)-\sum_{j=k}^{a-1}\Delta f(j)=f(a)+\Delta f(a)(k-a)+\sum_{j=k}^{a-1}(\Delta f(a)-\Delta f(j))\\
&=&f(a)+\Delta f(a)(k-a)+\sum_{j=k}^{a-1}\sum_{l=j}^{a-1}\Delta^2 f(l)
=f(a)+\Delta f(a)(k-a)+\sum_{l=k}^{a-1}\Delta^2 f(l)(l-k+1),
\end{eqnarray*}
so that
\begin{eqnarray*}
\left|f(k)-f(a)-\Delta f(a)(k-a)\right|\leq \|\Delta^2f\|_\infty \sum_{l=k}^{a-1}(l-k+1)
=\frac12\|\Delta^2f\|_\infty (a-k)(a-k+1),
\end{eqnarray*}
that is, (\ref{ingremain}) holds true in this case as well. The proof of Theorem \ref{peccatidelamathematique} is done.
\qed

\bigskip

{\bf To go further}.
A multivariate extension of Theorem \ref{PSTU} can be found in \cite{PecZheng}.
The reference \cite{PecLac} contains several explicit applications of the tools developed in this section.

\noindent
\section{Fourth Moment Theorem and free probability}\label{sec:free}

To conclude this survey, we shall explain how the Fourth Moment Theorem \ref{NP} extends in the theory of free probability,
which provides a convenient framework for
investigating limits of random matrices.
We start with a short introduction to free probability.  We refer to \cite{nicaspeicher} for a systematic presentation and
to \cite{bianespeicher} for specific results on Wigner multiple integrals.

\bigskip

{\bf Free tracial probability space}.
A {\em non-commutative probability space} is a von Neumann algebra $\mathscr{A}$ (that is,
an algebra of operators on a complex
separable Hilbert space,
closed under adjoint and convergence in the weak operator topology) equipped with a {\em trace} $\ff$,
that is, a unital linear functional (meaning preserving the identity) which is weakly
continuous, positive (meaning $\ff(X)\ge 0$
whenever $X$ is a non-negative element of $\mathscr{A}$; i.e.\ whenever $X=YY^\ast$
for some $Y\in\mathscr{A}$), faithful (meaning that if
$\ff(YY^\ast)=0$ then $Y=0$), and tracial (meaning that $\ff(XY)=\ff(YX)$ for all
$X,Y\in\mathscr{A}$, even though in general $XY\ne YX$).

\bigskip

{\bf Random variables}.
In a non-commutative probability space, we refer to the self-adjoint elements of the
algebra as
{\em random variables}.  Any
random variable $X$ has
a {\em law}: this is the unique
probability measure $\mu$ on $\R$
with the same moments as $X$; in other words, $\mu$ is such that
\begin{equation}\label{mu}
\int_{\R} x^k d\mu(x) = \ff(X^k),\quad k\geq 1.
\end{equation}
(The existence and uniqueness of $\mu$ follow from the positivity of $\ff$,  see \cite[Proposition 3.13]{nicaspeicher}.)

\bigskip

{\bf Convergence in law}.
We say that a sequence $(X_{1,n},\ldots,X_{k,n})$, $n\geq 1$, of random vectors {\it converges in law}
to a random vector $(X_{1,\infty},\ldots,X_{k,\infty})$, and we write \[
(X_{1,n},\ldots,X_{k,n})\overset{\rm law}{\to} (X_{1,\infty},\ldots,X_{k,\infty}),
\]
to indicate the convergence in the sense of (joint) moments, that is,
\begin{equation}\label{cvlaw}
\lim_{n\to\infty}\ff\left(\poly(X_{1,n},\ldots,X_{k,n})\right) = \ff\left(\poly(X_{1,\infty},\ldots,X_{k,\infty})\right),
\end{equation}
for any polynomial $\poly$ in $k$ non-commuting variables.

We say that a sequence $(F_n)$ of {\it non-commutative stochastic processes} (that is, each $F_n$ is a
one-parameter family of
self-adjoint operators $F_n(t)$ in  $(\mathscr{A},\ff)$) {\it converges in the sense of finite-dimensional
distributions} to a non-commutative stochastic process $F_\infty$,
and we write \[
F_{n}\overset{\rm f.d.d.}{\to} F_\infty,
\]
to indicate that, for any $k\geq 1$ and any $t_1,\ldots,t_k\geq 0$,
\[
(F_{n}(t_1),\ldots,F_{n}(t_k))\overset{\rm law}{\to} (F_{\infty}(t_1),\ldots,F_{\infty}(t_k)).
\]

\bigskip

{\bf Free independence}.
In the free probability setting, the notion of {\em independence} (introduced by
Voiculescu in \cite{Voiculescu})
goes as follows. Let $\mathscr{A}_1,\ldots,\mathscr{A}_p$ be unital subalgebras of $\mathscr{A}$.  Let $X_1,\ldots, X_m$ be elements
chosen from among the $\mathscr{A}_i$'s such that, for $1\le j<m$,
two consecutive elements $X_j$ and $X_{j+1}$ do not come from the same $\mathscr{A}_i$ and are
such that $\ff(X_j)=0$ for each $j$.  The subalgebras $\mathscr{A}_1,\ldots,\mathscr{A}_p$ are said to be {\em free} or {\em freely
independent} if, in this circumstance,
\begin{equation}\label{free-def}
\ff(X_1X_2\cdots X_m) = 0.
\end{equation}
Random variables are called freely independent if the unital algebras
they generate are freely independent.  Freeness is in general much more complicated than classical independence. For example, if
$X,Y$ are free and $m,n\geq 1$, then by (\ref{free-def}),
\[
\ff\big((X^m-\ff(X^m)1)(Y^n-\ff(Y^n)1)\big)=0.
\]
By expanding (and using the linear property of $\ff$), we get
\begin{equation}\label{covcov}
\ff(X^mY^n)=\ff(X^m)\ff(Y^n),
\end{equation}
which is what we would expect under classical independence. But, by setting
$X_1=X_3=X-\ff(X)1$ and $X_2=X_4=Y-\ff(Y)$
in (\ref{free-def}), we also have
\[
\ff\big((X-\varphi(X)1)(Y-\varphi(Y)1)(X-\varphi(X)1)(Y-\varphi(Y)1)\big)=0.
\]
By expanding, using (\ref{covcov}) and the tracial property of $\ff$ (for instance $\ff(XYX)=\ff(X^2Y)$)
we get
\begin{eqnarray*}
\ff(XYXY)
&=&\ff(Y)^2\ff(X^2)+\ff(X)^2\ff(Y^2)-\ff(X)^2\ff(Y)^2,
\end{eqnarray*}
which is different from $\ff(X^2)\ff(Y^2)$, which is what one would have obtained if $X$ and $Y$ were
classical independent random variables.
Nevertheless, if $X,Y$ are freely independent, then their joint moments
are determined by the moments of $X$ and $Y$ separately,
exactly as in the classical case.

\bigskip

{\bf Semicircular distribution}.
The {\em semicircular distribution}
$\mathcal{S}(m,\sigma^2)$
with mean $m\in\R$ and variance $\sigma^2>0$
 is the probability distribution
\begin{equation} \label{eq semicircle}
\mathcal{S}(m,\sigma^2)(dx) = \frac{1}{2\pi \sigma^2} \sqrt{4\sigma^2-(x-m)^2}\,{\bf 1}_{\{|x-m|\le 2\sigma\}}\,dx.
\end{equation}
If $m=0$, this distribution is symmetric around $0$,
and therefore its odd moments are all $0$. A simple calculation shows that the even centered moments are given by
(scaled) {\em Catalan numbers}: for non-negative integers $k$,
\[ \int_{m-2\sigma}^{m+2\sigma} (x-m)^{2k} \mathcal{S}(m,\sigma^2)(dx) = C_k \sigma^{2k}, \]
where $C_k = \frac{1}{k+1}\binom{2k}{k}$
(see, e.g., \cite[Lecture 2]{nicaspeicher}).
In particular, the variance is $\sigma^2$ while the centered fourth moment is $2\sigma^4$.
The semicircular distribution plays here the role of the Gaussian distribution.
It has the following similar properties:
\begin{enumerate}
\item If $S\sim\mathcal{S}(m,\sigma^2)$ and $a,b\in\R$, then $aS+b\sim \mathcal{S}(am+b,a^2\sigma^2)$.
\item If $S_1\sim\mathcal{S}(m_1,\sigma_1^2)$ and $S_2\sim\mathcal{S}(m_2,\sigma_2^2)$ are freely independent,
then $S_1+S_2\sim\mathcal{S}(m_1+m_2,\sigma_1^2+\sigma_2^2)$.
\end{enumerate}

\bigskip

{\bf Free Brownian Motion}.
A {\it free Brownian motion} $S=\{S(t)\}_{t\geq 0}$ is a non-commutative stochastic process
with the following defining characteristics:
\begin{itemize}
\item[(1)] $S(0) = 0$.
\item[(2)] For $t_2>t_1\geq 0$, the law of $S(t_2)-S(t_1)$ is the semicircular distribution of mean 0 and variance $t_2-t_1$.
\item[(3)] For all $n$ and $t_n>\cdots>t_2>t_1>0$, the increments $S(t_1)$, $S(t_2)-S(t_1)$, \ldots, $S(t_n)-S(t_{n-1})$ are
freely independent.
\end{itemize}

We may think of free Brownian motion as `infinite-dimensional matrix-valued Brownian motion'.

\bigskip

{\bf Wigner integral}.
Let $S=\{S(t)\}_{t\geq 0}$ be a free Brownian motion.
Let us quickly sketch out the construction of the {\it Wigner integral} of $f$ with respect to $S$.
For an indicator function $f={\bf 1}_{[u,v]}$, the Wigner integral of $f$ is defined
by
\[
\int_0^\infty {\bf 1}_{[u,v]}(x)dS(x)=S(v)-S(u).
\]
We then extend this definition by linearity to simple functions of the form
$
f=\sum_{i=1}^k \alpha_i {\bf 1}_{[u_i,v_i]},
$
where $[u_i,v_i]$ are disjoint intervals of $\R_+$.
Simple computations show that
\begin{eqnarray}
\ff\left(\int_0^\infty f(x)dS(x)\right)&=&0\label{wigner1}\\
\ff\left(\int_0^\infty f(x)dS(x)\times\int_0^\infty g(x)dS(x)\right)&=&\langle f,g\rangle_{L^2(\R_+)}.\label{wigner2}
\end{eqnarray}
By isometry, the definition of $\int_0^\infty f(x)dS(x)$ is extended to all $f\in L^2(\R_+)$,
and (\ref{wigner1})-(\ref{wigner2}) continue to hold in this more general setting.

\bigskip

{\bf Multiple Wigner integral}.
Let $S=\{S(t)\}_{t\geq 0}$ be a free Brownian motion, and let $q\geq 1$ be an integer.
When $f\in L^2(\R_+^q)$ is real-valued, we write $f^*$ to indicate the function of $L^2(\R_+^q)$ given by
$f^*(t_1,\ldots,t_q)=f(t_q,\ldots,t_1)$.

Following \cite{bianespeicher}, let us quickly sketch out the construction of the {\it multiple Wigner
integral} of $f$ with respect to $S$.
Let $D^q\subset\R_+^q$ be the collection of all diagonals, i.e.
\begin{equation}\label{dp}
D^q=\{(t_1,\ldots,t_q)\in\R_+^q:\,t_i=t_j\mbox{ for some $i\neq j$}\}.
\end{equation}
For an indicator function $f={\bf 1}_A$, where $A\subset\R_+^q$ has the form
$
A=[u_1,v_1]\times\ldots\times [u_q,v_q]
$
with $A\cap D^q=\emptyset$, the $q$th multiple Wigner integral of $f$
is defined
by
\[
I^S_q(f)=(S(v_1)-S(u_1))\ldots (S(v_q)-S(u_q)).
\]
We then extend this definition by linearity to simple functions of the form
$
f=\sum_{i=1}^k \alpha_i {\bf 1}_{A_i},
$
where
$
A_i=[u^i_1,v^i_1]\times\ldots\times [u^i_q,v^i_q]
$
are disjoint $q$-dimensional rectangles as above which do not meet the diagonals.
Simple computations show that
\begin{eqnarray}
\ff(I^S_q(f))&=&0\label{isomfree1}\\
\ff(I^S_q(f)I^S_q(g))&=&\langle f,g^*\rangle_{L^2(\R_+^q)}.\label{isomfree}
\end{eqnarray}
By isometry, the definition of $I^S_q(f)$ is extended to all $f\in L^2(\R_+^q)$,
and (\ref{isomfree1})-(\ref{isomfree}) continue to hold in this more general setting.
If one wants $I^S_q(f)$ to be a random variable,
it is necessary for $f$ to be {\it mirror symmetric}, that is, $f=f^*$
(see \cite{knps}).
Observe that $I^S_1(f)=\int_0^\infty f(x)dS(x)$ when $q=1$.
We have moreover
\begin{equation}\label{isom-dif-free}
\ff(I^S_p(f)I^S_q(g))=0\,\,\mbox{ when $p\neq q$, $f\in L^2(\R_+^p)$ and $g\in L^2(\R_+^q)$}.
\end{equation}

When $r\in\{1,\ldots,p\wedge q\}$, $f\in L^2(\R_+^p)$ and
$g\in L^2(\R_+^q)$,
let us  write $f\overset{r}\frown g$ to indicate the $r$th {\it contraction} of $f$ and $g$, defined as
being the element of $L^2(\R_+^{p+q-2r})$ given by
\begin{eqnarray}
&&f\overset{r}\frown g (t_1,\ldots,t_{p+q-2r})\label{contr}\\
&=&\int_{\R_+^r}f(t_1,\ldots,t_{p-r},x_1,\ldots,x_r)g(x_r,\ldots,x_1,t_{p-r+1},\ldots,t_{p+q-2r})dx_1\ldots dx_r.\notag
\end{eqnarray}
By convention, set $f\overset{0}{\frown} g= f\otimes g$ as
being the tensor product of $f$ and $g$.
Since $f$ and $g$ are not necessarily symmetric functions, the position
of the identified variables $x_1,\ldots,x_r$ in (\ref{contr}) is important, in contrast
to what happens in classical probability.
Observe moreover that
\begin{equation}\label{cauchy}
\|f\overset{r}{\frown} g\|_{L^2(\R_+^{p+q-2r})}\leq \|f\|_{L^2(\R_+^p)}\|g\|_{L^2(\R_+^q)}
\end{equation}
by Cauchy-Schwarz, and also that
$f\overset{q}{\frown} g=\langle f,g^*\rangle_{L^2(\R_+^q)}$ when $p=q$.

We have the following {\it product formula} (see \cite[Proposition 5.3.3]{bianespeicher}), valid for any $f\in L^2(\R_+^p)$
and $g\in L^2(\R_+^q)$:
\begin{equation}\label{prodfree}
I^S_p(f)I^S_q(g)=\sum_{r=0}^{p\wedge q} I^S_{p+q-2r}(f\overset{r}{\frown}g).
\end{equation}
We deduce (by a straightforward induction) that, for any $e\in L^2(\R_+)$ and any $q\geq 1$,
\begin{equation}\label{lien-tch}
U_q\left(\int_0^\infty e(x)dS_x\right) = I^S_q(e^{\otimes q}),
\end{equation}
where $U_0=1$, $U_1=X$, $U_2=X^2-1$, $U_3=X^3-2X$, $\ldots$, is the
sequence of Tchebycheff polynomials of second kind
(determined by the recursion $XU_k=U_{k+1}+U_{k-1}$),
$\int_0^\infty e(x)dS(x)$ is understood as a Wigner integral,
and $e^{\otimes q}$ is the $q$th tensor product of $e$. This is the exact analogue of (\ref{linkhermite})
in our context.

\bigskip

We are now in a position to offer a free analogue of the Fourth Moment Theorem \ref{T : NPNOPTP}, which reads as follows.

\begin{thm}[Kemp, Nourdin, Peccati, Speicher, 2011; see \cite{knps}]\label{thm-knps}
Fix an integer $q\geq 2$ and let $\{S_t\}_{t\geq 0}$ be a free Brownian motion.
Whenever $f\in L^2(\R_+^q)$, set $I^S_q(f)$ to denote the $q$th multiple Wigner integrals of $f$
with respect to $S$.
Let $\{F_n\}_{n\geq 1}$
be a sequence of Wigner multiple integrals of the form
\[
F_n=I^S_q(f_n),
\]
where each $f_n\in L^2(\R_+)$ is mirror-symmetric, that is, is such that $f_n=f_n^*$.
Suppose moreover that $\varphi(F_n^2)\to 1$ as $n\to\infty$.
Then, as $n\to\infty$, the following two assertions are equivalent:
\begin{enumerate}
\item[(i)] $F_n\overset{\rm Law}\to S_1\sim \mathcal{S}(0,1)$;
\item[(ii)] $\ff(F_n^4)\to 2=\ff(S_1^4)$.
\end{enumerate}
\end{thm}
{\it Proof} (following \cite{yet}).
Without loss of generality and for sake of simplicity, we suppose that $\ff(F_n^2)=1$ for all $n$
(instead of $\ff(F_n^2)\to 1$ as $n\to\infty$).
The proof of the implication $(i)\Rightarrow (ii)$ being trivial by the very definition
of the convergence in law in a free tracial probability space, we only concentrate on
the proof of $(ii)\Rightarrow (i)$.

Fix $k\geq 3$. By iterating the product formula (\ref{prodfree}), we can write
\[
F_n^k=I^S_q(f_n)^k=\sum_{(r_1,\ldots,r_{k-1})\in A_{k,q}}
I^S_{kq-2r_1-\ldots-2r_{k-1}}\big(
(\ldots((
f_n\overset{r_1}{\frown}f_n
)\overset{r_2}{\frown} f_n)
\ldots)\overset{r_{k-1}}{\frown}f_n
\big),
\]
where
\begin{eqnarray*}
A_{k,q}&=&\big\{
(r_1,\ldots,r_{k-1})\in \{0,1,\ldots,q\}^{k-1}:\,r_2\leq 2q-2r_1,\,\,r_3\leq 3q-2r_1-2r_2,\ldots,\\
&&\hskip7cm r_{k-1}\leq(k-1)q-2r_1-\ldots-2r_{k-2}
\big\}.
\end{eqnarray*}
By taking the $\ff$-trace in the previous expression and taking into account that (\ref{isomfree1}) holds, we deduce that
\begin{equation}\label{ffFn}
\ff(F_n^k)=\ff(I^S_q(f_n)^k)=\sum_{(r_1,\ldots,r_{k-1})\in B_{k,q}}
(\ldots((
f_n\overset{r_1}{\frown}f_n
)\overset{r_2}{\frown} f_n)
\ldots)\overset{r_{k-1}}{\frown}f_n,
\end{equation}
with
\[
B_{k,q}=\big\{(r_1,\ldots,r_{k-1})\in A_{k,q}:\,2r_1+\ldots+2r_{k-1}=kq\big\}.
\]
Let us decompose $B_{k,q}$ into $C_{k,q}\cup E_{k,q}$, with $C_{k,q}=B_{k,q}\cap\{0,q\}^{k-1}$ and $E_{k,q}=B_{k,q}\setminus C_{k,q}$.
We then have
\begin{eqnarray*}
\ff(F_n^k)&=& \sum_{(r_1,\ldots,r_{k-1})\in C_{k,q}}
\big(
(\ldots((
f_n\overset{r_1}{\frown}f_n
)\overset{r_2}{\frown} f_n)
\ldots)\overset{r_{k-1}}{\frown}f_n
\big)\\
&&+\sum_{(r_1,\ldots,r_{k-1})\in E_{k,q}}
\big(
(\ldots((
f_n\overset{r_1}{\frown}f_n
)\overset{r_2}{\frown} f_n)
\ldots)\overset{r_{k-1}}{\frown}f_n
\big).
\end{eqnarray*}
Using the two relationships
$f_n\overset{0}{\frown} f_n = f_n\otimes f_n$ and
\[
f_n\overset{q}{\frown} f_n =
\int_{\R_+^q}f_n(t_1,\ldots,t_q)f_n(t_q,\ldots,t_1)dt_1\ldots dt_q=\|f_n\|^2_{L^2(\R_+^q)}=1,
\]
it is evident that
$(\ldots((
f_n\overset{r_1}{\frown}f_n
)\overset{r_2}{\frown} f_n)
\ldots)\overset{r_{k-1}}{\frown}f_n
=1$ for all $(r_1,\ldots,r_{k-1})\in C_{k,q}$. We deduce that
\begin{eqnarray*}
\ff(F_n^k)&=& \#C_{k,q} +\sum_{(r_1,\ldots,r_{k-1})\in E_{k,q}}
\big(
(\ldots((
f_n\overset{r_1}{\frown}f_n
)\overset{r_2}{\frown} f_n)
\ldots)\overset{r_{k-1}}{\frown}f_n
\big).
\end{eqnarray*}
On the other hand, by applying (\ref{ffFn}) with $q=1$, we get that
\begin{eqnarray*}
\ff(S_1^k)&=&\ff(I^S_1({\bf 1}_{[0,1]})^k)=\sum_{(r_1,\ldots,r_{k-1})\in B_{k,1}}
(\ldots((
{\bf 1}_{[0,1]}\overset{r_1}{\frown}{\bf 1}_{[0,1]}
)\overset{r_2}{\frown} {\bf 1}_{[0,1]})
\ldots)\overset{r_{k-1}}{\frown}{\bf 1}_{[0,1]}\\
&=&\sum_{(r_1,\ldots,r_{k-1})\in B_{k,1}} 1
= \# B_{k,1}.
\end{eqnarray*}
But it is clear
that $C_{k,q}$ is in bijection with $B_{k,1}$ (by dividing all the $r_i$'s in $C_{k,q}$ by $q$).
Consequently,
\begin{equation}\label{clekey}
\ff(F_n^k)= \ff(S_1^k) +\sum_{(r_1,\ldots,r_{k-1})\in E_{k,q}}
\big(
(\ldots((
f_n\overset{r_1}{\frown}f_n
)\overset{r_2}{\frown} f_n)
\ldots)\overset{r_{k-1}}{\frown}f_n
\big).
\end{equation}
Now, assume that $\varphi(F_n^4)\to \varphi(S_1^4)=2$ and let us show that
$\varphi(F_n^k)\to \varphi(S_1^k)$ for {\it all} $k\geq 3$.
Using that $f_n=f_n^*$, observe that
\begin{eqnarray*}
&&f_n\overset{r}{\frown} f_n(t_1,\ldots,t_{2q-2r})\\
&=&\int_{\R_+^r} f_n(t_1,\ldots,t_{q-r},s_1,\ldots,s_r)f_n(s_r,\ldots,s_1,t_{q-r+1},\ldots,t_{2q-2r})
ds_1\ldots ds_r
\\
&=&\int_{\R_+^{r}} f_n(s_r,\ldots,s_1,t_{q-r},\ldots,t_1)f_n(t_{2q-2r},\ldots,t_{q-r+1},s_1,\ldots,s_r)
ds_1\ldots ds_{r}
\\
&=&f_n\overset{r}{\frown} f_n(t_{2q-2r},\ldots,t_1)=(f_n\overset{r}{\frown} f_n)^*(t_1,\ldots,t_{2q-2r}),
\end{eqnarray*}
that is, $f_n\overset{r}{\frown} f_n = (f_n\overset{r}{\frown} f_n)^*$.
On the other hand, the  product formula (\ref{prodfree}) leads to
$F_n^2= \sum_{r=0}^q I^S_{2q-2r}(f_n\overset{r}{\frown} f_n)$.
Since two multiple integrals of different orders are orthogonal (see (\ref{isom-dif-free})), we deduce
that
\begin{eqnarray}
\ff(F_n^4) &=& \|f_n\otimes f_n\|^2_{L^2(\R_+^{2q})} + \big(\|f_n\|^2_{L^2(\R_+^q)}\big)^2+
\sum_{r=1}^{q-1}
\langle f_n\overset{r}{\frown} f_n, (f_n\overset{r}{\frown} f_n)^*\rangle_{L^2(\R_+^{2q-2r})}
  \notag\\
&=&  2\|f_n\|^4_{L^2([0,1]^q)}+\sum_{r=1}^{q-1}\|f_n\overset{r}{\frown} f_n\|^2_{L^2(\R_+^{2q-2r})}
=2+\sum_{r=1}^{q-1}\|f_n\overset{r}{\frown} f_n\|^2_{L^2(\R_+^{2q-2r})}.\label{fourthpower}
\end{eqnarray}
Using that $\varphi(F_n^4)\to 2$, we deduce that
\begin{equation}\label{step1}
\|f_n\overset{r}{\frown} f_n\|^2_{L^2(\R_+^{2q-2r})}\to 0\quad\mbox{
for all $r=1,\ldots,q-1$}.
\end{equation}
Fix $(r_1,\ldots,r_{k-1})\in E_{k,q}$ and let
$j\in\{1,\ldots,k-1\}$ be the smallest integer such that $r_j\in\{1,\ldots,q-1\}$.
Then:
\begin{eqnarray*}
&&\big| (\ldots((f_n\overset{r_1}{\frown}f_n
)\overset{r_2}{\frown} f_n)
\ldots)\overset{r_{k-1}}{\frown}f_n \big| \\
&=&\big|
(\ldots((f_n\overset{r_1}{\frown}f_n)\overset{r_2}{\frown}f_n)\ldots   \overset{r_{j-1}}{\frown}  f_n) \overset{r_{j}}{\frown} f_n)
\overset{r_{j+1}}{\frown} f_n)\ldots) \overset{r_{k-1}}{\frown} f_n\big|\\
&=&\big|
(\ldots((f_n\otimes\ldots \otimes f_n) \overset{r_{j}}{\frown} f_n)
\overset{r_{j+1}}{\frown} f_n)\ldots) \overset{r_{k-1}}{\frown} f_n\big|
\quad
\mbox{(since $f_n\overset{q}{\frown} f_n=1$)}\\
&=&\big|
(\ldots((f_n\otimes\ldots \otimes f_n)\otimes (f_n \overset{r_{j}}{\frown} f_n))
\overset{r_{j+1}}{\frown} f_n)\ldots) \overset{r_{k-1}}{\frown} f_n\big|\\
&\leq&
\|(f_n\otimes\ldots \otimes f_n)\otimes (f_n  \overset{r_{j}}{\frown} f_n)\|
\|f_n\|^{k-j-1}\quad
\mbox{\small (Cauchy-Schwarz)\normalsize}\\
&=&\|f_n  \overset{r_{j}}{\frown} f_n\|\quad\mbox{(since $\|f_n\|^2=1$)}\\
&\to& 0\quad\mbox{as $n\to\infty$}\quad\mbox{by (\ref{step1})}.
\end{eqnarray*}
Therefore, we deduce from (\ref{clekey}) that $\varphi(F_n^k)\to \varphi(S_1^k)$,
which is the desired conclusion and concludes the proof of the theorem.\qed

\bigskip

During the proof of Theorem \ref{thm-knps}, we actually showed (see indeed (\ref{fourthpower})) that the two assertions $(i)$-$(ii)$ are both equivalent to a third one, namely
\[
\mbox{$(iii)$:}\quad \quad
\|f_n\overset{r}{\frown} f_n\|^2_{L^2(\R_+^{2q-2r})}\to 0\,\,\mbox{for all $r=1,\ldots,q-1$}.
\]
Combining $(iii)$ with Corollary \ref{T : NPNOPTP}, we immediately deduce an interesting
transfer principle for translating results between the classical and free chaoses.
\begin{cor}
Fix an integer $q\geq 2$, let $\{B_t\}_{t\geq 0}$ be a standard Brownian motion and let $\{S_t\}_{t\geq 0}$ be a free Brownian motion.
Whenever $f\in L^2(\R_+^q)$, we write $I^B_q(f)$ (resp. $I^S_q(f)$) to indicate the $q$th multiple Wiener integrals of $f$
with respect to $B$ (resp. $S$).
Let $\{f_n\}_{n\geq 1}\subset L^2(\R_+^q)$ be a sequence of symmetric functions and let $\sigma>0$ be a finite constant.
Then, as $n\to\infty$, the following two assertions hold true.
\begin{enumerate}
\item[(i)] $E[I_q^B(f_n)]\to q!\sigma^2$ if and only if $\ff(I_q^S(f_n)^2)\to\sigma^2$.
\item[(ii)] If the asymptotic relations in $(i)$ are verified, then $I^B_q(f_n)\overset{\rm law}{\to}\mathcal{N}(0,q!\sigma^2)$
if and only if $I^S_q(f_n)\overset{\rm law}{\to}\mathcal{S}(0,\sigma^2)$.
\end{enumerate}
\end{cor}

\bigskip

{\bf To go further}.
A multivariate version of Theorem \ref{thm-knps} (free counterpart of Theorem \ref{PecTud}) can be found in \cite{NPS}.
In \cite{NP-poisson} (resp. \cite{tetilla}), one exhibits a version of Theorem \ref{thm-knps} in which the semicircular law in the limit is replaced by the free Poisson law (resp. the so-called tetilla law).
An extension of Theorem \ref{thm-knps} in the context of the $q$-Brownian motion (which is an interpolation
between the standard Brownian motion corresponding to $q=1$ and the free Brownian motion corresponding to $q=0$)
is given in \cite{qBM}.

\end{document}